                     \def\version{26.09.2021}                          %

\documentclass[reqno,11pt]{amsart}  

\usepackage{amsmath}
\usepackage{amssymb} 
\usepackage{amsthm}
\usepackage{mathtools}   
\usepackage{graphicx}
\usepackage{graphics}  
\usepackage[pdftex]{color}
\usepackage{color}
\usepackage{verbatim} 
\usepackage[normalem]{ulem} 
\usepackage[mathscr]{eucal}
\usepackage{upgreek}
\usepackage{enumerate} 
\usepackage{color}
\usepackage{dsfont}
\usepackage{tikz}
\usetikzlibrary{decorations.pathreplacing,angles,quotes}
\usepackage{subcaption}
 \usepackage{scalerel}
 
\usepackage[titletoc]{appendix}

\DeclareRobustCommand{\bigO}{%
  \text{\usefont{OMS}{cmsy}{m}{n}O}%
}

\numberwithin{equation}{section}  
\renewcommand{\theequation}{\thesection.\arabic{equation}} 
 
\newtheorem{theorem}{Theorem}[section] 
\newtheorem{lemma}[theorem]{Lemma} 
\newtheorem{prop}[theorem] {Proposition} 
\newtheorem{cor}[theorem]  {Corollary} 
\newtheorem{remark}[theorem]  {Remark}

\theoremstyle{definition}

 \usepackage{float}
\floatstyle{boxed}
\restylefloat{figure}

\DeclareMathAlphabet{\mathpzc}{OT1}{pzc}{m}{it}

\newcommand{\babs}[1]{{\bigl\lvert #1\bigr\rvert}}
\newcommand{\Babs}[1]{{\Bigl\lvert #1\Bigr\rvert}}
\newcommand{\bnorm}[1]{{\boldsymbol{\lvert} #1\boldsymbol{\rvert}}}
\newcommand{\tnorm}[1]{{|\hspace{-0.35mm}\lVert #1\rVert\hspace{-0.35mm}|}}

\DeclarePairedDelimiter{\abs}{\lvert}{\rvert}
\DeclarePairedDelimiter{\norm}{\lVert}{\rVert}

\newcommand{\bA} {\boldsymbol{A}} 
\newcommand{\bB} {\boldsymbol{B}} 
\newcommand{\bC} {\boldsymbol{C}} 
\newcommand{\bE} {\boldsymbol{E}}
\newcommand{\bF} {\boldsymbol{F}}
\newcommand{\bG} {\boldsymbol{G}}
\newcommand{\bH} {\boldsymbol{H}}
\newcommand{\bL} {\boldsymbol{L}} 
\newcommand{\bM} {\boldsymbol{M}} 
\newcommand{\bP} {\boldsymbol{P}} 
\newcommand{\bQ} {\boldsymbol{Q}} 
\newcommand{\bR} {\boldsymbol{R}} 
\newcommand{\bT} {\boldsymbol{T}} 
\newcommand{\bU} {\boldsymbol{U}} 
\newcommand{\bV} {\boldsymbol{V}} 
\newcommand{\bW} {\boldsymbol{W}}
\newcommand{\bX} {\boldsymbol{X}} 
\newcommand{\bY} {\boldsymbol{Y}}
\newcommand{\bZ} {\boldsymbol{Z}} 

\newcommand{\bb} {\boldsymbol{b}}
\newcommand{\bd} {\boldsymbol{d}}
\newcommand{\be} {\boldsymbol{e}}
\newcommand{\bff} {\boldsymbol{f}}
\newcommand{\bp} {\boldsymbol{p}}
\newcommand{\bq} {\boldsymbol{q}}
\newcommand{\bt} {\boldsymbol{t}}
\newcommand{\bu} {\boldsymbol{u}}
\newcommand{\bx} {\boldsymbol{x}}
\newcommand{\by} {\boldsymbol{y}}
\newcommand{\bz} {\boldsymbol{z}}
 
\newcommand{\cbH} {\boldsymbol{\mathcal H}}
\newcommand{\cbR} {\boldsymbol{\mathcal R}}
\newcommand{\cbX} {\boldsymbol{\mathcal X}}
\newcommand{\cbY} {\boldsymbol{\mathcal Y}}
\newcommand{\cbV} {\boldsymbol{\mathcal V}}
\newcommand{\cbT} {\boldsymbol{\mathcal T}}

\newcommand{\bbX} {\mathbb X}
\newcommand{\bbY} {\mathbb Y}
\newcommand{\bbZ} {\mathbb Z}
 
\newcommand{\bXi} {\boldsymbol{\Xi}}
\newcommand{\bgamma} {\boldsymbol{\gamma}}
\newcommand{\bs} {\boldsymbol{\sigma}}
\newcommand{\ba} {\boldsymbol{\alpha}}
\newcommand{\balpha} {\boldsymbol{\alpha}}
\newcommand{\bbeta} {\boldsymbol{\beta}}
\newcommand{\blambda}{\boldsymbol{\lambda}}
\newcommand{\bmu}{\boldsymbol{\mu}}
\newcommand{\bcdot} {\boldsymbol{\cdot}}
\newcommand{\one} {\boldsymbol{1}}
\newcommand{\zero} {\boldsymbol{0}}

\newcommand{\G} {\Gamma} 
\renewcommand{\L} {\Lambda} 
\renewcommand{\bL}{\boldsymbol{\Lambda}}
\newcommand{\brho}{\boldsymbol{\rho}}
\renewcommand{\O} {\Omega} 
\newcommand{\Phn} {\Phi} 
\renewcommand{\S} {\Sigma} 

\def\a{\alpha}
\def\b{\beta}
\def\d{\delta} 
\newcommand{\e} {\varepsilon} 
\newcommand{\eps}{\varepsilon} 
\newcommand{\ph} {\varphi} 
\newcommand{\phn} {\mathsf{\phi}} 
\newcommand{\g} {\gamma}
\def\l{\lambda} 
\newcommand{\n} {\nu} 
\def\r{\varrho} 
\newcommand{\s} {\sigma}
\def\t{\tau} 
\def\th{\theta} 
\def\u{\upsilon} 
\newcommand{\z} {\zeta} 
\newcommand{\lsc} {\rm \tiny l.s.c.}

\newfam\Bbbfam 
\font\tenBbb=msbm10 
\font\sevenBbb=msbm7 
\font\fiveBbb=msbm5 
\textfont\Bbbfam=\tenBbb 
\scriptfont\Bbbfam=\sevenBbb 
\scriptscriptfont\Bbbfam=\fiveBbb 
\def\Bbb{\fam\Bbbfam \tenBbb} 
\newcommand{\B}     {\mathbb{B}} 
 \newcommand{\C}     {\mathbb{C}} 
\newcommand{\R}     {\mathbb{R}} 
\newcommand{\Z}     {\mathbb{Z}} 
\newcommand{\N}     {\mathbb{N}} 
\renewcommand{\P}   {\mathbb{P}} 
\newcommand{\D}     {\mathbb{D}} 
\newcommand{\E}     {\mathbb{E}} 
\newcommand{\Q}     {\mathbb{Q}} 
\newcommand{\T}     {\mathbb{T}} 
 \newcommand{\floor}[1]{\left\lfloor #1 \right\rfloor}
\newcommand{\ceil}[1]{\left\lceil #1 \right\rceil}
\newcommand{\counting} {\mathfrak{c}} 
\newcommand{\smfrac}[2]{\textstyle{\frac {#1}{#2}}}
\def\cov{\mbox{\rm Cov}} 
\def\var{\mbox{\rm var}} 
\def\1{{\mathchoice {1\mskip-4mu\mathrm l}      
{1\mskip-4mu\mathrm l} 
{1\mskip-4.5mu\mathrm l} {1\mskip-5mu\mathrm l}}} 
\newcommand{\ssup}[1] {{\scriptscriptstyle{({#1}})}} 

\newcommand{\ssupp}[1]{\scaleto{({#1})}{2.5pt}}

\newcommand{\PMF} {\ssup{\scaleto{\scriptscriptstyle{PMF}}{2.55pt}}}

\newcommand{\PMFbc} {\ssup{\scaleto{\scriptscriptstyle{PMF}}{2.55pt}, \bc}}

\newcommand{\CMF} {\ssup{\scaleto{CMF}{2.55pt}}}

\newcommand{\HYL} {\ssup{\scaleto{HYL}{2.55pt}}}

\newcommand{\stup}[1] {{\scriptscriptstyle{({\scriptsize#1}})}} 
\newtheoremstyle{thm}{2ex}{2ex}{\itshape\rmfamily}{} 
{\bfseries\rmfamily}{}{1.7ex}{} 
 
\newtheoremstyle{rem}{1.3ex}{1.3ex}{\rmfamily}{} 
{\itshape\rmfamily}{}{1.5ex}{} 
 
\newenvironment{proofsect}[1] 
{\vskip0.1cm\noindent{\scshape #1.}\hskip0.5cm} 

 \def\e{{\rm e}}
\def\dx{{\rm d} x}
\def\dy{{\rm d} y}
\def\Poi{{\mathsf{Poi}}}
\def\Bin{{\mathsf{Bin}}}
\def\Geo{{\mathsf{Geo}}}
\def\NBin{{\mathsf{NBin}}}
\def\Ber{{\mathsf Ber}}
\def\Uni{{\mathsf U}}
\def\Norm{{\mathsf N}}
\def\Var{{\mathsf{Var}}}
\def\Cov{{\mathsf{Cov}}}
  
\def\ox{\by}
\def\oy{\bX}
 
\newcommand{\cP} {{\mathcal P}} 
\newcommand{\cR} {{\mathcal R}} 
\newcommand{\cS} {{\mathcal S}} 
\newcommand{\cT} {{\mathcal T}} 
\newcommand{\cW} {{\mathcal W}} 
\newcommand{\cU} {{\mathcal U}} 
\newcommand{\cX} {{\mathcal X}} 
\newcommand{\cA} {{\mathcal A}} 

\newcommand{\Acal}   {{\mathcal A }}
\newcommand{\Bcal}   {{\mathcal B }}
\newcommand{\Ccal}   {{\mathcal C }} 
\newcommand{\Dcal}   {{\mathcal D }} 
\newcommand{\Ecal}   {{\mathcal E }} 
\newcommand{\Fcal}   {{\mathcal F }} 
\newcommand{\Gcal}   {{\mathcal G }} 
\newcommand{\Hcal}   {{\mathcal H }} 
\newcommand{\Ical}   {{\mathcal I }} 
\newcommand{\Jcal}   {{\mathcal J }} 
\newcommand{\Kcal}   {{\mathcal K }} 
\newcommand{\Lcal}   {{\mathcal L }} 
\newcommand{\Mcal}   {{\mathcal M }} 
\newcommand{\Ncal}   {{\mathcal N }} 
\newcommand{\Ocal}   {{\mathcal O }} 
\newcommand{\Pcal}   {{\mathcal P }} 
\newcommand{\Qcal}   {{\mathcal Q }} 
\newcommand{\Rcal}   {{\mathcal R }} 
\newcommand{\Scal}   {{\mathcal S }} 
\newcommand{\Tcal}   {{\mathcal T }} 
\newcommand{\Ucal}   {{\mathcal U }} 
\newcommand{\Vcal}   {{\mathcal V }} 
\newcommand{\Wcal}   {{\mathcal W }} 
\newcommand{\Xcal}   {{\mathcal X }} 
\newcommand{\Ycal}   {{\mathcal Y }} 
\newcommand{\Zcal}   {{\mathcal Z }}

\newcommand{\Ascr} {\mathscr{A}}
\newcommand{\Bscr} {\mathscr{B}}
\newcommand{\Tscr} {\mathscr{T}}

\newcommand{\Cscr} {\mathscr{C}}
\newcommand{\Lscr} {\mathscr{L}}
\newcommand{\Escr}{\mathscr{E}}
\newcommand{\Fscr}{\mathscr{F}}
\newcommand{\Gscr}{\mathscr{G}}
\newcommand{\Hscr} {\mathscr{H}}
\newcommand{\Iscr}{\mathscr{I}}
\newcommand{\Jscr}{\mathscr{J}}
\newcommand{\Kscr}{\mathscr{K}}
\newcommand{\Nscr}{\mathscr{N}}
\newcommand{\Oscr}{\mathscr{O}}
\newcommand{\Pscr}{\mathscr{P}}
\newcommand{\Qscr}{\mathscr{Q}}
\newcommand{\Rscr} {\mathscr{R}}
\newcommand{\Sscr}{\mathscr{S}}

\newcommand{\Uscr}{\mathscr{U}}
\newcommand{\Wscr}{\mathscr{W}}

\newcommand{\m} {{\mathfrak m}}
\newcommand{\ip} {{\mathit p}}
\newcommand{\fX} {{\mathfrak X}}

\newcommand{\Abscr}{\boldsymbol{\mathscr{A}}}
 \newcommand{\ex}{{\rm e}} 
\newcommand{\com}{{\rm c}} 
\renewcommand{\d}{{\rm d}} 
\newcommand{\per}{{\text{\rm per}}} 
\newcommand{\opn}{\operatorname} 
\newcommand{\Li}{\operatorname{Li}\,}
\newcommand{\Leb}{{\rm Leb}} 
\newcommand{\Sym}{\mathfrak{S}}
\newcommand{\id}{{\sf{I}}} 
\newcommand{\esssup}{{\operatorname {esssup}}} 
\newcommand{\supp}{{\operatorname {supp}}} 
\newcommand{\tay}{{\operatorname {Tay}}}
\newcommand{\sign}{{\operatorname {sign}\,}}
\newcommand{\dist}{{\operatorname {dist}}} 
\newcommand{\diam}{{\operatorname {diam}}} 
\newcommand{\dime}{{\operatorname {dim}\,}} 
\newcommand{\const}{{\operatorname {cst.}\,}} 
\newcommand{\Dir}{{\operatorname {Dir}\,}}
\newcommand{\tr}{{\operatorname {Tr}}}
\newcommand{\Exp}{\mathscr{E}\kern-0.2mm{\operatorname{xp}}}
\newcommand{\Log}{\mathscr{L}\kern-0.2mm{\operatorname{og}}}
\newcommand{\Id}{{\operatorname {Id}}}
\newcommand{\heap}[2]{\genfrac{}{}{0pt}{}{#1}{#2}} 
\newcommand{\bc}{{\operatorname{bc}}}
\renewcommand{\emptyset} {\varnothing} 
\newcommand{\tO} {\tilde{\Omega}} 
\newcommand{\p} {\partial} 
\newcommand\embed{\hookrightarrow}
\newcommand\poly{{\mathrm{poly}}}

\DeclareMathOperator*{\argmin}{arg\,min}

\def\Bin{{\mathsf B}}
\def\Bet{\mbox{\boldmath$\beta$}}
\def\Gam{{\mathsf G}}
\newcommand{\Ex}{\mathsf{Exp}}
\def\Hyp{{\mathsf H}}
\def\Mult{{\mathsf M}}
\def\free{\mbox{free}}
\def\Norm{{\mathsf N}}
\def\Poi{{\mathsf P}}
\def\Uni{{\mathsf U}}
\def\Fi{{\mathsf F}}
\def\te{{\large\mathsf t}}
\def\Chi{{\mathbf\chi}^2}
\def\post{{\mathbf\pi}}
\def\prior{{\mathbf\a}}

\renewcommand{\Pr} {\mathsf{P}}
\newcommand{\Er} {\mathsf{E}}
\newcommand{\Vsf}{\mathsf{V}}
\newcommand{\Jsf}{\mathsf{J}}
\newcommand{\psf}{\mathsf{p}}
\newcommand{\Qsf}{\mathsf{Q}}

\newcommand{\Nsf}{\mathsf{N}}

\renewcommand{\baselinestretch}{1.06} 
\makeindex

\begin{document}

\title[\hfill \hfill]
{Large deviations analysis for random combinatorial partitions with counter terms }


\author{Stefan Adams and Matthew Dickson}
\address{Mathematics Institute, University of Warwick, Coventry CV4 7AL, United Kingdom}
\email{S.Adams@warwick.ac.uk,   dickson@math.lmu.de   }

\thanks{}
  
\date{}

 


\begin{abstract}
In this paper, we study various models for random combinatorial partitions using large deviation analysis for diverging scale of the reference process. Scaling limits of similar models have been studied recently \cite{FSa,FSb} going back to \cite{Ver96}. After studying the reference model, we provide a complete analysis of two mean field models, one of which is well-know \cite{BCMP05} and the other one is the cycle mean field model. Both models show critical behaviour despite their rate functions having unique minimiser. The main focus is then a model with negative counter term, the probabilistic version of the so-called \emph{Huang-Yang-Luttinger} (HYL) model \cite{BLP}. Criticality in this model is the existence of a  critical parameter for which two simultaneous minimiser exists. At criticality an order parameter is introduced as the double limits for the density of cycles with diverging length, and as such it extends recent work \cite{AD21}. 

\end{abstract} 

\maketitle

\section{Introduction}\label{Intro}
In this paper, we study various models for random combinatorial partitions using large deviation analysis for diverging scale of the reference process.
\subsection{The reference measure}
Random combinatorial partitions arise in many areas of mathematics as number theory, combinatorics, probability and statistical mechanics, as illustrated in \cite{Ver96} and further developed in recent work \cite{FSa,FSb}. 

The problem is about decomposing an integer $ N\in\N $ into a sum of positive integers, $ N=x_1+\cdots +x_m , m\in\N $. A partition of $N$ is then the equivalence class of sequences $ (x_1,\ldots, x_m) $ whose terms sum up to $N$ and where two sequences are equivalent if they differ by a permutation. Following Vershik \cite{Ver96}, we describe partitions by their \emph{occupation} sequences $ \lambda=(\lambda_k)_{k\in\N} $, where $ \lambda_k $ denotes the number of elements equal to $k$ in a sequence representative of the partition, thus $ \sum_{k\in\N} k\lambda_k =N $. The \emph{occupation numbers} $ \lambda_k $ represent also the number of cycles of length in a permutation of $N$ elements, in this way the sequence $\lambda=(\lambda_k)_{k\in\N} $ gives the cycle structure of permutations. We are concerned with random partitions models where we assign statistical weights for various models. All weights are so-called tilts of the reference measure. For the reference weights we denote $ \Ncal_k, k\in\N $, the  Poisson distributed occupation number with parameter $ \abs{\L_N}q_k^\ssup{\alpha} $, where $\abs{\L_N}=(2N)^d, d, N\in\N $, and
\begin{equation}\label{weightP}
q_k^\ssup{\alpha}=\frac{\e^{\beta k\alpha}}{(4\pi\beta)^{d/2}k^{1+d/2}}\,,\quad\mbox{ with }  \alpha\le 0,\beta>0, k\in\N\,.
\end{equation}
The reference measure is then the superposition of all Poisson processes and is itself a Poisson process with parameter
$$
\overline{q}:=\sum_{k\in\N} \; q_k^\ssup{\alpha}\,,
$$ see, e.g. \cite{R09}. We are concerned with diverging scales $ \abs{\L_N}\to\infty $ as $ N\to\infty $ and therefore  introduce the \emph{empirical cycle count} or the \emph{empirical occupation count} as
\begin{equation}\label{empirical}
\blambda_N:=\big(\Ncal_k/\abs{\L_N}\big)_{k\in\N}\,.
\end{equation}
For a sequence $ x=(x_k)_{k\in\N} \in\ell_1(\R_+) $ with $ \abs{\L_N} x_k\in\N_0 $, we denote $ \Qsf $ the probability distribution of the reference process, that is, the probability that the empirical cycle count is equal to $x$ is given by
\begin{equation}\label{referenceweight}
\Qsf\big(\blambda_N=(x_k)_{k\in\N}\big)=\e^{-\abs{\L_N}\overline{q}^\ssup{\alpha}}\prod_{k\in\N}\frac{\big(\abs{\L_N}q_k^\ssup{\alpha}\big)^{\abs{\L_N}q_k^\ssup{\alpha}}}{\big(\abs{\L_N}x\big)!}\,.
\end{equation}
The probability weights for the reference process are a special class of multiplicative weights. General multiplicative weights have been  introduced by Vershik \cite{Ver96} and  are studied and analysed  for scaling limits in \cite{FSa,FSb}. 
Our study is not concerned with scaling limits but with large deviation limits as $ N\to\infty $. All our models are given as probability measures in the sequence space $ \ell_1(\R_+) $, and we denote $ \nu_{N,\alpha}=\Qsf\circ\blambda_N^{-1} $ the reference measure on $ \ell_1(\R_+) $.  Before we introduce our tiled models in Section~\ref{sec-models}, we motivate the specific weights of the reference measure. These weight appear in calculations of  the trace of the Gibbs density operator for an ideal (non-interacting) gas of Bosons, a class of quantum particles obeying certain permutations statistics according to the representation of the permutation group of $N$ particle indices. As identical quantum particles cannot be distinguished, one needs to symmetrise their labels. Using the Feynman-Kac formula one obtains the trace as an expectation for $N$ Brownian bridges under symmetrised initial-terminal conditions, see \cite{A08,ACK} and \cite{AD06,AK08} for the symmetrisation and random permutations and partitions. In the so-called \emph{grand canonical ensemble} with random number $N$ of particles we recover our reference weights for a special choice of the scales $ \L_N=[-N,N]^d\subset\R^d $ and empty boundary conditions, see for instance \cite{ACK,AD21}. In Appendix~\ref{appIBG} we summarise the results for the reference measure and well-know results for the ideal Bose gas and its condensation phenomenon. 
In parts our work is related to \cite{BCMP05} and \cite{Lew86,BLP}. However, it should be noted that in these studies the random weights for the partitions originate from an energy (Fourier space) representation of the underlying physical models. In this energy setting, condensation is concerned with the zero energy mode in the systems whereas in our study we are concerned with infinitely long cycles, see \cite{S02,A08,AK08,ACK,AV17}. Recently there is some work on random Euclidean permutations \cite{EP19} which are different from the ones studied in this paper as their weights depend on the spatial distance as well. It can be promising to extend our models to include spatial dependence in the future.

\subsection{Models}\label{sec-models}
Our models are given as probability measures on $ \ell_1(\R_+) $  by various tilts of the reference measures. The tilts are so-called \emph{Hamiltonian functions} $H$ on $ \ell_1(\R_+) $ such that the new measure is given by the Radon-Nikodym density $\frac{\e^{-\beta \abs{\L_N} H}}{Z_N} $ with respect to the reference measure, where $ \beta > 0 $ and where $Z_N$ is the normalising constant also called \emph{partition function}.  Under the given tilts the new measures put higher probability weight to cycle counts with lower values of the Hamiltonian function. We define two so-called \emph{mean field} models and another  model with counter terms.  

\noindent For $ a\ge 0 $, define the \emph{cycle-mean field model} (CMF),
\begin{equation}
\begin{aligned}
H^\CMF(x)&=\frac{a}{2}\left(\sum_{k=1}^\infty\; x_k\right)^2\,,\quad x\in\ell_1(\R_+)\,,\\
\nu^{\CMF}_{N,\alpha}(\d x)&=\frac{\e^{-\beta\abs{\L_N} H^\CMF(x)}}{Z_N^\CMF(\beta,\alpha)}\nu_{N,\alpha}(\d x)\,,\quad Z_N^\CMF(\beta,\alpha)=\E_{\nu_{N,\alpha}}\left[\e^{-\beta\abs{\L_N} H^\CMF}\right]\,.
\end{aligned}
\end{equation}
The measure $ \nu_{N,\alpha}^\CMF $ gives higher weight to cycle counts with smaller values of the total number of cycles. Any given cycle count  is a partition of the number
\begin{equation}\label{particles}
\Nsf_{N}:=\sum_{k\in\N} k\Ncal_k\,,
\end{equation}
which we also call the '\emph{number of particles}'. 
The number $\Nsf_N $ is only lower semi continuous and not upper semi continuous (see \cite{ACK}). For $ a\ge 0 $ and any $\mu\in\R $, define the \emph{particle-mean-field model} (PMF), 
\begin{equation}\label{PMF}
\begin{aligned}
H^{\PMF}_\mu(x)& =  -\mu\sum_{k=1}^\infty k  x_k +    \frac{a}{2}\left(\sum_{k=1}^\infty kx_k\right)^2,\qquad x\in\ell_1(\R_+)\,, \\
\nu_{N,\mu,\alpha}^{\PMF} (\d x) &= \frac{\e^{-\abs{\L_N}\beta H^{\PMF}_\mu(x)}}{Z_N^\PMF(\beta,\alpha,\mu)}\nu_{N,\alpha}(\d x)\,,\quad Z_N^\PMF(\beta,\alpha,\mu)=\E_{\nu_{N,\alpha}}\left[\e^{-\abs{\L_N}\beta H^{\PMF}_\mu}\right]\,. 
\end{aligned}
\end{equation}
This models puts lower probability weight on cycle counts with large number of particles and has been studied in the literature, see \cite{BCMP05} for a nice summary. The name refers to the number of physical particles in the system. The measure, though well-known in the physics literature, is different in character from the partitions weights studied in \cite{Ver96} and \cite{FSa,FSb}. Namely, as the earlier weights are just product of the single weights, the squared term in the Hamiltonian $ H^\PMF $ creates product of weights of pairwise different weight numbers.

The major novelty of our large deviation analysis concerns an substantial extension of the so-called HYL-model (\emph{Huang-Yang-Luttinger model}) studied in \cite{BLP}. On one hand we replace the cycle weights originating from the energy representation in \cite{BLP} by our cycle weights stemming from spatial representation of the partition function, and on the other hand we obtain higher level large deviation principles allowing a detailed insight in the structure of the minimiser and possible phases and phase transitions.  But significantly, we can dispense a major technical assumption in \cite{BLP}, see details on this in \cite{AD21}.
For any $ a\ge b> 0 $ and any $ \alpha\le 0$, $\mu\in\R $, define the HYL-model by
\begin{equation}
\begin{aligned}
H^\HYL_\mu(x)&=-\mu\sum_{k=1}^\infty k x_k +\frac{a}{2}\left(\sum_{k=1}^\infty k x_k\right)^2-\frac{b}{2}\sum_{k=1}^\infty k^2 x_k^2\,,\qquad x\in\ell_1(\R_+)\,,\\
\nu^\HYL_{N,\alpha,\mu}(\d x)&=\frac{\e^{-\beta\abs{\L_N} H^\HYL_\mu(x)}}{Z^\HYL_{N,\alpha,\mu}}\,\nu_{N,\alpha}(\d x)\,,\quad Z_N^\HYL(\beta,\alpha,\mu)=\E_{\nu_{N,\alpha}}\left[\e^{-\abs{\L_N}\beta H^{\HYL}_\mu}\right]\,.
\end{aligned}
\end{equation}

\subsection{Organisation and summary of the paper}

The papers is structured into four chapters and an appendix. The first two sections present our results and the remaining ones collect our proofs. The appendix has three parts; Appendix~\ref{app-Bose} introduces the Bose functions as a class of poly-logarithmic functions. In Appendix~\ref{appIBG} we present a rigorous large deviation analysis of the reference process (ideal Bose gas) as courtesy for the reader and to present a new proof method using Baldi's theorem \cite{DZ09}. Appendix~\ref{app-Lambert} defines the Lambert $W$ function and collects some properties. This function is vital for the analysis of our rate functions for the CMF and the HYL model.
All our large deviation results are in Section~\ref{sec:LDPsection}. The main focus is on large deviation results for the empirical cycle counts, for the reference measure and the PMF model we complement this with large deviation principles for the  \emph{empirical density}
\begin{equation}
\brho_{N}:=\frac{1}{\abs{\L_N}}\Nsf_N\,.
\end{equation}
Note that $\brho_N=D\left(\blambda_{N}\right)$, where $D\colon \ell_1\left(\R_+\right)\to \R\cup\left\{+\infty\right\}$,
\begin{equation}\label{eqn:particledensity}
    D\left(x\right) := \sum^\infty_{k=1}kx_k\,,
\end{equation}
which is given in Section~\ref{sec:LDPTotalDensity}.
The large deviation principle for the CMF model in Theorem~\ref{THM-CMF} uses standard large deviation methods. The major obstacle for the large deviation principles for remaining models is that the number of 'particles' $ \Nsf_N $ or its density $ D$ is only lower semi continuous and not upper semi continuous. We prove the large deviation principle for the PMF model in Theorem~\ref{THM-PMF}  using a lower semi continuous regularisation of the Hamiltonian function following \cite{GZ93} in conjunction with a finite-dimensional approximation. The proof of the large deviation principle of the HYL model in Theorem~\ref{Thm:HYL} combines the methods for the PMF with two different representations of the Hamiltonian function adapted to the lower and the upper bound of the large deviation principle, all of which is in Section~\ref{sec-HYL}. The second main body of work is the variational analysis of the rate functions and models in Section~\ref{sec:Variational}. We analyse the zeroes of our rate functions and derive representations of the so-called limiting \emph{pressures}, that is, the limiting logarithmic moment generating functions of our models. We analyse the pressure functions as function of the parameter $ \alpha $ (CMF model) and $ \mu $ (PMF and HYL model) whose derivative give the density of particles. Critical behaviour is present when the derivative of the pressure is different from the expected density which is  given by the density of the zeroes of the rate function.

For the CMF model we find that the unique zeroes are given as functions of the Lambert $W$ function, and the zeroes and the analysis of the pressure  show that the CMF model has similar properties as the reference measure including the so-called criticality in terms of \emph{Bose-Einstein condensation} BEC defined for the ideal Bose gas in Appendix~\ref{appIBG}. The difference is only in terms of the critical density which is now a function of the parameter $a$ and the Lambert $W$ function.
 
The corresponding analysis of the PMF model in Section~\ref{sec-PMFa} shows also unique zeroes for the rate function but this time the phase transitions establishes as a change in the pressure density relation in Proposition~\ref{zeroPMF} and Proposition~\ref{P:fPMF}. This leads to the conjecture that the \emph{condensate density} is given by $ (\mu/a-\varrho(\alpha))_+ $, see for illustration Figure~\ref{fig:PMFdensity}, where the dashed red line represents the expected density and the horizontal line the derivative. Naturally, the analysis of the rate function for the HYL model is more complex, see Proposition~\ref{zeroHYL}. Here, the Lambert $W$ function plays a major role and uniqueness of the zeroes is only given for certain parameter regimes, see Theorem~\ref{THM-uniqueness}. The main result concerns a critical parameter $ \mu=\mu^* $ in Theorem~\ref{THM-nonuniqueness} when two simultaneous zeroes of the rate function exist. We believe that the zero with the lower density represents the system out of the condensate whereas the other one represents the system with condensation. The condensate itself is then the leftover probability mass and is conjecture to be the probability mass of the  '\emph{infinitely long cycles}', see \cite{S02,A08,ACK}. In order to shed some light on the condensation of cycles of diverging length we introduce in Section~\ref{sec-BEC} an order parameter as a double limit of the density of diverging cycle length. 

In future work we shall address scaling limits as in \cite{Ver96,FSa,FSb} where we hope that the detailed analysis of the rate function zeroes can help establishing different scaling limits. Another question concerns the concentration of measure around the two distant rate function zeroes of the HYL model at criticality.

\section{Large Deviation Principles}\label{sec:LDPsection}

\subsection{Large Deviations of the Empirical Cycle Count}\label{sec:LDPEmpiricalCycleCount}

We present our large deviation results which are all based on the large deviation principle for the reference measure  $ \nu_{N,\alpha}$ in Appendix~\ref{appIBG}. The large deviation principle and thermodynamic results  are given in Appendix~\ref{appIBG}. For the convenience of the reader and better understanding of our main results we present all details of the proofs of the large deviation principle and the exponential tightness using Baldi's Lemma in Appendix~\ref{appIBG}. For the following we recall the rate function for the ideal Bose gas, (see Theorem~\ref{THM-Ideal})
$$
I_\alpha(x) = \sum^\infty_{k=1}\frac{x_k}{\beta}\Big(\log\frac{x_k}{q_k^{\ssup{\alpha}}} - 1\Big) + \bar{q}^{\ssup{\alpha}}/\beta\,.
$$

From the construction of the CMF, PMF and HYL models via the ideal Bose gas model, it is natural to expect that their LDPs may be derived with an application of Varadhan's Lemma. Nevertheless our results are more sophisticated because the tilts are not continuous in the latter two cases.

\begin{theorem}[\textbf{Large deviations principle  for CMF  models}]\label{THM-CMF}
For any $ d\in\N,  a>  0 $ and $ \alpha\leq0$ the following holds. The sequence $ \big(\nu_{N,\alpha}^{\CMF}\big)_{N\ge 1}   $ satisfies an LDP on $ \ell_1(\R_+) $ with rate $ \beta\abs{\L_N} $ and rate function 
\begin{equation}
I^{\CMF}_\alpha(x) = H^{\CMF}(x)+ I_\alpha(x) - \inf_{y\in \ell_1\left(\R\right)}\{H^{\CMF}(y)+I_\alpha(y)\}\,.
\end{equation}
\end{theorem}

\begin{theorem}[\textbf{Large deviation principle for PMF models}]\label{THM-PMF}
For any $ d\in\N,  a>  0, \alpha\le 0 $, and $ \mu\in\R$  the following holds. The sequence $ \big(\nu_{N,\alpha,\mu}^{\PMF}\big)_{N\ge 1}   $ satisfies a LDP on $ \ell_1(\R_+) $ with rate $ \beta \abs{\L_N} $ and rate function 
\begin{equation}
I^\PMF_{\alpha,\mu}(x) =  I_\alpha(x)+ H^\PMF_{\mu,\lsc}(x)  - \inf_{y\in \ell_1\left(\R_+\right)}\{ I_\alpha(y)+H^{\PMF}_{\mu,\lsc}(y)\}\,,
\end{equation}
with
\begin{equation}\label{lscreg}
H^\PMF_{\mu,\lsc}(x)=H^\PMF_\mu(x)-\frac{1}{2a}\left(\mu-aD(x)\right)_+^2=\begin{cases} -\mu D(x)+\frac{a}{2}D(x)^2 &,D(x)\ge \frac{\mu}{a}\,,\\-\frac{\mu^2}{2a}& ,D(x)<\frac{\mu}{a}\,.\end{cases}
\end{equation}

\end{theorem}

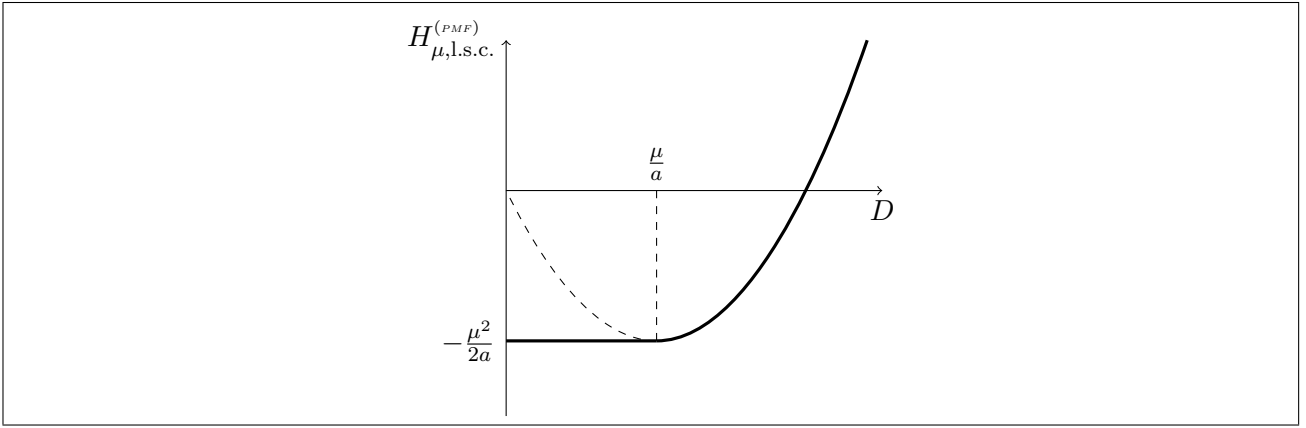
\begin{figure}
	\centering
	\begin{tikzpicture}[scale=2]
	\draw[->] (0,0) -- (2.5,0) node[below]{$D$};
	\draw[->] (0,-1.5) -- (0,1) node[left]{$H^{\PMF}_{\mu,\lsc}$};
	\draw[very thick] (1,-1) parabola (2.4,1);
	\draw[very thick] (0,-1) node[left]{$-\frac{\mu^2}{2a}$} -- (1,-1);
	\draw[dashed] (1,-1) parabola (0,0);
	\draw[dashed] (1,-1) -- (1,0) node[above]{$\frac{\mu}{a}$};
	\end{tikzpicture}
	\caption{Sketch of $H^\PMF_{\mu,\lsc}$ as a function of the total particle density $D$. \label{fig:PMF_lsc}}
\end{figure}

\begin{remark}
\begin{enumerate}[(a)]
\item In Figure~\ref{fig:PMF_lsc} we illustrate the lower semi continuous regularisation $H^\PMF_{\mu,\lsc} $ as a function of the density.

\item For $ \mu\le 0 $,  the rate function in Theorem~\ref{THM-PMF} reads
$$
I^{\PMF}_{\alpha,\mu}(x) = I_\alpha(x) + H^{\PMF}_{\mu}(x)- \inf_{y\in \ell_1\left(\R_+\right)}\{ I_\alpha(y)+H^{\PMF}_{\mu}(y)\}\,.
$$
\end{enumerate}
\hfill $ \diamond $
\end{remark}

\begin{theorem}[\textbf{Large deviations principle for HYL models}]\label{Thm:HYL}

For any $ d\in\N,  a\ge b\geq 0 , \alpha\le 0 $, and $ \mu\in\R$    the following holds. The sequence $ \big(\nu_{N,\alpha,\mu}^{\HYL}\big)_{N\ge 1}   $ satisfies a LDP on $ \ell_1(\R_+) $ with rate $ \beta \abs{\L_N} $ and rate function 
\begin{equation}
I^{\HYL}_{\alpha,\mu}(x) = I_\alpha(x)+  H^{\HYL}_{\mu,\lsc}(x) - \inf_{y\in \ell_1\left(\R_+\right)}\{I_\alpha(y)+H^{\HYL}_{\mu,\lsc}(y)\} \,,
\end{equation}
with
\begin{equation}
\begin{aligned}
H^{\HYL}_{\mu,\lsc}(x)&=H^\HYL_\mu(x)-\frac{1}{2(a-b)}\left(\mu-aD(x)\right)_+^2\\ &= -\frac{b}{2}\sum_{k=1}^\infty k^2 x_k^2 + \begin{cases} -\mu D(x)+\frac{a}{2}D(x)^2 &,D(x)\ge \frac{\mu}{a},\\ - \frac{b}{a-b}\left(-\mu D(x)+\frac{a}{2}D(x)^2\right) -\frac{\mu^2}{2(a-b)} &, D(x)<\frac{\mu}{a}.\end{cases}
\end{aligned}
\end{equation}

\end{theorem}

In both, Theorem~\ref{THM-PMF} and Theorem~\ref{Thm:HYL}, we have presented the lower semicontinuous regularisations of the interaction energy densities. It is this that allows us to overcome the lack of continuity in the original energy densities.

\subsection{Large Deviations of the Empirical Density}
\label{sec:LDPTotalDensity}

The large deviation principles for the empirical density require an independent proof as the contraction principle does not work directly. To see that,
recall that the 'particle` number $ \sum_{k\in\N}k\Ncal_k $   is only lower semicontinuous and not upper semicontinuous, a proof via the contraction principle is only feasible if one considers cut-off versions of the empirical density $ \brho_N^{\ssup{K}}=\frac{1}{\abs{\L_N}}\sum^K_{k=1}k\Ncal_k $ followed  by analysing the limit $ K\to\infty $ for the corresponding rate functions. We do not follow this approach here and briefly outline a direct approach as follows.

\begin{prop}\label{logmomentdensity}
Let $ \alpha\le 0 $, then for all $ t\in\R $, the logarithmic moment generating function is 
\begin{equation}
\begin{aligned}
\mathcal{L}(t)&:=\lim_{N\to\infty} \frac{1}{\beta \abs{\L_N}}\log\E_{\nu_{N,\alpha}^{\ssup{\bc}}}\big[\ex^{\beta\abs{\L_N} t\sum_{k=1}^\infty k\blambda_N^{\ssup{k}}}\big]=  
 \sum_{k=1}^\infty \frac{q_k^{\ssup{\alpha}}}{\beta}\big(\ex^{\beta tk}-1\big)=\begin{cases} +\infty &, \mbox{ if } t>\abs{\alpha}\,,\\
\in \R & ,\mbox{ if } \alpha+t\le 0\,.
\end{cases}
\end{aligned}
\end{equation}
\end{prop}
The following large deviation results uses the critical density for the ideal Bose gas, the thermodynamic limit of the pressure and the free energy  defined, respectively,  in Appendix~\ref{appIBG}. Denote $ Q_{N,\alpha}=\Qsf\circ \brho_N^{-1} $ the distribution of $ (\brho_N)_{N\ge 1} $ with chemical potential $ \alpha\le 0 $ and define the distribution $ Q_{N,\mu,\alpha}^\PMF $ via its Radon-Nikodym density
\begin{equation}
\frac{\d Q_{N,\mu,\alpha}^\PMF}{\d Q_{N,\alpha}}(x)=\frac{\exp\left(-\abs{\L_N}\beta\left(-\mu x+\frac{a}{2}x^2\right)\right)}{Z_N^\PMF(\beta,\mu,\alpha)}\,.
\end{equation}

\begin{theorem}\label{THM-densityLDP}
Let $ d\in\N $ and $  \beta>0 $.
\begin{enumerate}[(a)]
\item For any $ \alpha<0 $, the sequence $ \left(Q_{N,\alpha}\right)_{N\ge 1} $  satisfies a LDP on $ \R$ with rate $ \beta\abs{\L_N} $ and rate function
\begin{equation}
J_\alpha(x)= \begin{cases} p(\beta,\alpha)+f(\beta,x)-\alpha x &, \mbox{ if } x\in[0,\varrho_{\rm c}(d)] \mbox{ for } d\ge 3\wedge x\in[0,\infty)\mbox{ for } d=1,2\,,\\
+\infty &, \mbox{ if } x\notin[0,\varrho_{\rm c}(d)]\,.\end{cases}
\end{equation}

\item For any $ \alpha<0 $ and $ \mu\in\R $, the sequence $ \left(Q_{N,\mu,\alpha}^\PMF\right)_{N\ge 1} $ satisfies a LDP on $ \R $ with rate $ \beta\abs{\L_N} $ and rate function
\begin{equation}
J_{\mu,\alpha}^\PMF\left(x\right)= \begin{cases} J_\alpha\left(x\right) -\left(\mu+\alpha\right) x+\frac{a}{2}x^2 -\mathsf{N}  &,\mbox{ if } x\in[0,\varrho_{\rm c}(d)] \mbox{ for } d\ge 3\wedge x\in[0,\infty)\mbox{ for } d=1,2\,,\\
+\infty &, \mbox{ if } x\notin[0,\varrho_{\rm c}(d)]\,,\end{cases}
\end{equation}
where
$$
\mathsf{N} =  \inf_{y\in\R}\big\{J_\alpha\left(y\right) -\left(\mu+\alpha\right) y+\frac{a}{2} y^2\big\}\,.$$
\end{enumerate} 
\end{theorem}

\begin{remark}\label{RemdensityLDP}
The results in Theorem~\ref{THM-densityLDP} make the heuristic derivations in \cite{Lew86} rigorous and extend them to all $ \alpha<0 $ and $ \mu\in\R $. 
The free energy of the PMF model is $ f^\PMF(\beta,\varrho)=f(\beta,\varrho)+\frac{a}{2}\varrho^2 $, whereas the pressure
\begin{equation}\label{PMFpressure}
p^\PMF(\beta,\mu,\alpha)=\sup_{x\in\R}\left\{(\mu+\alpha)x-\frac{a}{2}x^2-f(\beta,x)\right\}=\sup_{x\in\R}\left\{ (\mu+\alpha)x-f^\PMF(\beta,x)\right\}\,.
\end{equation}
The  HYL model requires higher level  empirical functionals as the energy cannot be expressed as a functional of the empirical particle density. 
\hfill $ \diamond $
\end{remark}

\section{Variational analysis, pressure representations, and condensation}
\label{sec:Variational}
Our large deviation analysis in Section~\ref{sec:LDPsection} is complemented by a complete analysis for the rate functions and pressure representations in Section~\ref{sec-minimiser}. In Section~\ref{sec-BEC},  we finally study the onset of criticality know as the \emph{Bose-Einstein condensation} (BEC) and discuss the relevance of our results.

\subsection{Variational analysis and pressure representations}\label{sec-minimiser}

The results for our reference measure are collected in Appendix~\ref{appIBG} in Proposition~\ref{zeroIdeal}, Proposition~\ref{P:ideal1}, and Proposition~\ref{IdealFreeEnergy}. Both mean field models, the CMF and the PMF model, are closely related to the ideal Bose gas.  Using our large deviation principles in Section~\ref{sec:LDPsection} and the zeroes of the rate functions we obtain the thermodynamic limit of the pressure in our various models.

\subsubsection{CMF model}\label{sec-thermCMF}
We collect the results for the first mean-field model. All proofs of this section are in Section~\ref{sec-proofs-CMF}.

\medskip

\begin{prop}\label{zeroCMF}
\begin{enumerate}[(a)]	
\item The rate function $ I^{\CMF}_\alpha$ has a unique zero at $\xi^\CMF\in\ell_1\left(\R\right)$ given by
	\begin{equation}
	\xi^\CMF_k = \frac{W_0\left(a\beta\bar{q}^{\ssup{\alpha}}\right)}{a\beta\bar{q}^{\ssup{\alpha}}}q^{\ssup{\alpha}}_k, \quad k\in\N\,,
	\end{equation}
	where $W_0$ is the real branch of the Lambert W function for non-negative arguments.
\item  Let $ \beta>0 $, $\alpha\le 0 $ and $ a\ge 0 $,  then
\begin{equation}\label{p-cmf}\begin{aligned}
p^\CMF(\beta,\alpha) =\lim_{N\to\infty}\frac{1}{\beta\abs{\L_N}}\log Z_N^\CMF(\beta,\alpha)=\frac{1}{a\beta^2}W_0\left(a\beta\bar{q}^{\ssup{\alpha}}\right)\left(1+\frac{1}{2}W_0\left(a\beta\bar{q}^{\ssup{\alpha}}\right)\right)\,.	
\end{aligned}
\end{equation}
\end{enumerate}
\end{prop}

\medskip

\begin{remark}
Definition and properties of the Lambert function are given in Appendix~\ref{app-Lambert}.
\hfill $ \diamond$
\end{remark}

\begin{prop}\label{P:densityCMF}
	\begin{enumerate}[(a)] 
		\item For $ \beta> 0 $, we define $ p^\CMF(\beta,\alpha)=+\infty $ for $ \alpha > 0 $. Then $ p^\CMF(\beta,\cdot) $ is a closed convex function on $ \R $.
		
		\item For $\beta>0$, $\alpha<0$, the pressure $p^\CMF(\beta,\alpha)$ is smooth with respect to $\alpha$. In particular,
			\begin{equation*}
			\frac{\d p^\CMF}{\d \alpha} = D\left(\xi^\CMF\right) = \frac{W_0\left(a\beta\bar{q}^\ssup{\alpha}\right)}{a\beta\bar{q}^\ssup{\alpha}}D\left(q^\ssup{\alpha}\right)\,.
			\end{equation*}

\item In the thermodynamic limit $ N\to\infty $,
\begin{equation}\label{criticalCMF}
\varrho^{\CMF}_{\rm c}(d):=\lim_{\alpha\uparrow 0}\left(\frac{\d}{\d\alpha}p^{\CMF}_{\L_N}(\beta,\alpha)\right)=\begin{cases} +\infty &, d=1,2\,,\\[1.5ex] \frac{W_0(a\beta\bar{q}^{\ssup{0}})}{a\beta\bar{q}^{\ssup{0}}} \varrho_{\rm c}(d)&, d\ge 3\,,
\end{cases}
\end{equation} where $ \varrho_{\rm c}(d) $ is the critical density for the ideal Bose gas, see \eqref{critical} in Appendix~\ref{appIBG}.
	\end{enumerate}
\end{prop}

\medskip

\begin{prop}\label{P:LFT-CMF}
	For $\varrho>0$, the  free energy of the CMF model is defined as the Legendre-Fenchel transform of the pressure, 
	\begin{equation}
	f^\CMF\left(\beta,\varrho\right) := \sup_{\alpha\in\R}\left\{\alpha\varrho - p^\CMF\left(\beta,\alpha\right)\right\}
	= \begin{cases}
	\varrho\alpha - p^\CMF\left(\beta,\alpha\right) &, \varrho\leq \varrho^\CMF_{\rm c}(d),\\
	-p^\CMF\left(\beta,0\right) &, \varrho\geq \varrho^\CMF_{\rm c}(d),
	\end{cases}
	\end{equation}
	where $\alpha$ is a solution to
	\begin{equation*}
	\frac{1}{a\beta\bar{q}^{\ssup{\alpha}}}W_0\left(a\beta \bar{q}^{\ssup{\alpha}}\right)D(q^{\ssup{\alpha}}) = \varrho\,,
	\end{equation*}
	which exists and is unique for $\varrho\leq \varrho^\CMF_{\rm c}(d)$.
\end{prop}

\bigskip

\begin{remark}[\textbf{Conclusions CMF model}]   The CMF model shows similar results as the reference measure (ideal Bose gas) in Appendix~\ref{appIBG}, e.g.,  the free energy is constant in the density beyond its specific  critical density. The critical density of the CMF model is different from the ideal Bose gas one. Using properties of the Lambert function, see Appendix~\ref{app-Lambert}, we know that 
$$
\begin{aligned}
\lim_{c\downarrow 0}\frac{W_0(cx)}{cx}&=1,\\
\lim_{c\to\infty}\frac{W_0(cx)}{cx}&=0.
\end{aligned}
$$
Hence, as the coupling parameter $a \to 0 $ vanishes, we obtain the critical ideal Bose gas density, and as $ a\to\infty $ the critical density decreases indicating BEC for much lower particle densities. Here, we refer to the definition of BEC for the ideal Bose gas as outlined in Appendix~\ref{appIBG}.  It is shown by \cite{S02} and \cite{BCMP05}, that BEC corresponds to loss of probability weights on finite cycles. When the coupling parameter $a$ increases  the number of finite cycles is suppressed in the probability measure,  and therefore the system undergoes a transition to a regime where the particle density is realised in so-called infinite cycles. The CMF model has not been studied in the literature so far, it shows similar behaviour as the ideal Bose gas because the Hamiltonian adds only weight on large numbers of cycles present. 
\hfill $ \diamond$
\end{remark}

\subsubsection{PMF model}\label{sec-PMFa}

We collect our findings for the PMF model. We obtain all results in \cite{BCMP05} with a completely different method for all values of the chemical potential, in addition, we compute the condensate density in Theorem~\ref{thm:PMFcondensate} below.  We identify regimes where the expected particle density equals the density of the rate function zero or not.
All proofs of this section are in Section~\ref{proofs-PMF}.
\medskip

\bigskip

\begin{prop}\label{zeroPMF}
\begin{enumerate}[(a)]
\item The rate function $ I^\PMF_{\alpha,\mu}$ has a unique zero at $ \xi^\PMF\in\ell_1(\R_+) $ where
\begin{equation*}
\xi^\PMF_k=q_k^{\ssup{\alpha}} \exp\left(\beta k\left(\mu-a\delta^*\right)_-\right),\quad k\in\N,
\end{equation*}
and $ \delta^*=\delta^*(\beta,\alpha, \mu,a) $ is given implicitly as the unique solution to the equation
\begin{equation}\label{delta}
\delta^*=\sum_{k=1}^\infty k q_k^{\ssup{\alpha}}\exp\left(\beta k\left(\mu-a\delta^*\right)_-\right)=D(\xi^\PMF)\,,
\end{equation}
\begin{equation}
\delta^*=\begin{cases} 
\in(0,\varrho(\alpha+\mu-a\delta^*)) & \,,\mbox{ for }\mu\le 0\,,\\[1ex] 
\begin{cases}  \in(\mu/a,\varrho(\alpha)) &\,, \mbox{ for }\mu< a\varrho(\alpha)\,,\\ \varrho(\alpha) &\,,  \mbox{ for } \mu\ge a\varrho(\alpha) \,, \end{cases}  & \, ,\mbox{ for } \mu>0\,,  \end{cases}
\end{equation}

where 
\begin{equation}\label{rhoalpha}
\varrho(\alpha):=\sum_{k\in\N} kq_k^{\ssup{\alpha}}\,,\quad \alpha\le 0\,.
\end{equation}

\medskip

\item   Let $\beta>0, \alpha\le 0, \mu\in\R$, and $ a\ge 0$, then
\begin{equation}\label{PMFpressure2}
\begin{aligned}
p^\PMF(\beta,\alpha,\mu) =\begin{cases} p(\beta,\alpha)+\mu^2/2a &,\mbox{ for } \mu\ge a\varrho(\alpha), \delta^*=\varrho(\alpha)\,,\\ \frac{a}{2}\big(\delta^*\big)^2+p(\beta,\alpha+\mu-a\delta^*) &, \mbox{ for } \mu<a\varrho(\alpha), \delta^*\in(\mu/a,\varrho(\alpha))\,.\end{cases}
\end{aligned}
\end{equation}

\item 

For $\beta>0$ and $ \alpha\le 0 $, the pressure $p^{\PMF}\left(\beta,\alpha,\cdot\right)\in C^{1}\left(\R\right)$ and is convex. In particular,
$$
\begin{aligned}
\frac{\d p^{\PMF}}{\d\mu}=\begin{cases} \varrho(\alpha+\mu-a\delta^*)  & , \mbox{ for } \mu<a\varrho(\alpha), \delta^*\in (\mu/a,\varrho)\,,\\
\mu/a & ,\mbox{ for } \mu\ge a\varrho(\alpha), \delta^*=\varrho(\alpha)\,,\end{cases}
\end{aligned}
$$
where $ \delta^* $ is solution of \eqref{delta}.

\end{enumerate}

\end{prop}

\bigskip

The following conclusion shows differences of the density as function of the  dimension and $ \alpha\le 0 $. The critical density $ \varrho_{\rm c}(d) $ is derived  in \eqref{critical} for the reference measure, 
$$
\varrho_{\rm c}(d)= \begin{cases}  +\infty &, d=1,2,\\   \frac{1}{\left(4\pi\beta\right)^\frac{d}{2}}\zeta\left(\frac{d}{2}\right) &, d\ge 3\,, \end{cases}
$$
where $ \zeta $ is the Riemann zeta function defined in \eqref{zeta}.
\begin{cor}
We have that
$$
\varrho(\alpha)=\begin{cases} \in (0,\infty) &, \alpha <0, d\ge 1\,,\\
 \infty & , \alpha\equiv 0 \wedge d=1,2\,,\\ 
 \varrho_{\rm c}(d) \in (0,\infty) &,\alpha\equiv 0 \wedge d\ge 3\,,\end{cases}
$$
\end{cor}

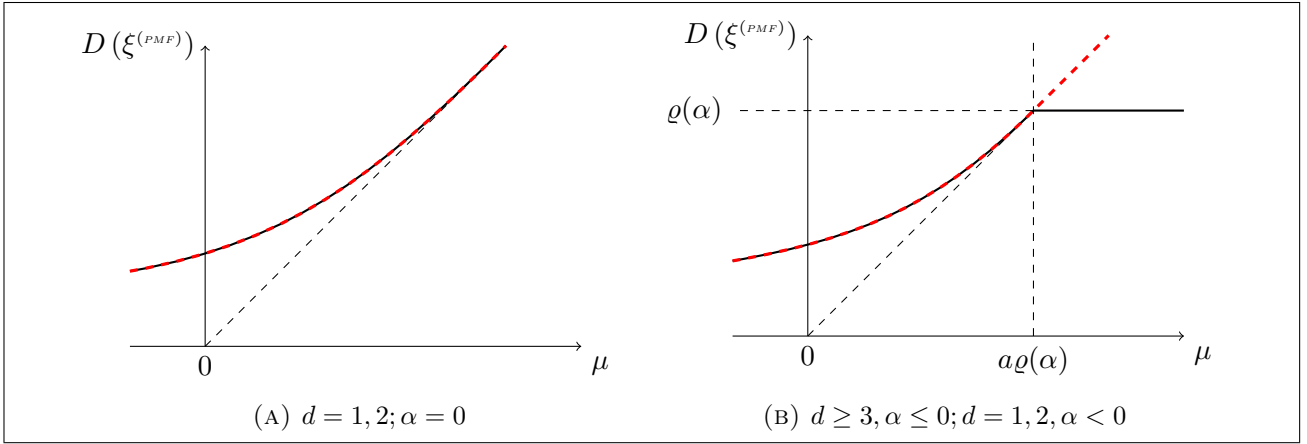
\begin{figure}
	\centering
	\begin{subfigure}[b]{0.45\textwidth}
		\begin{tikzpicture}[scale = 1]
		\draw[->] (-1,0) -- (5,0) node[below right]{$\mu$};
		\draw[->] (0,0) node[below]{$0$} -- (0,4) node[left]{$D\left(\xi^\PMF\right)$};
		\draw[dashed] (0,0) -- (4,4);
		\draw[thick] (-1,1) to [out=10,in=225] (4,4);
		\draw[very thick, red, dashed] (-1,1) to [out=10,in=225] (4,4);
		\end{tikzpicture}
		\caption{$d=1,2; \alpha=0 $}
		\label{fig:PMFdensity12}
	\end{subfigure}
	\begin{subfigure}[b]{0.45\textwidth}
		\begin{tikzpicture}[scale = 1]
		\draw[->] (-1,0) -- (5,0) node[below right]{$\mu$};
		\draw[->] (0,0) node[below]{$0$} -- (0,4) node[left]{$D\left(\xi^\PMF\right)$};
		\draw[dashed] (3,0) node[below]{$a\varrho(\alpha)$}--(3,4);
		\draw[dashed] (0,0) -- (3,3);
		\draw[dashed] (3,3) -- (-1,3) node[left]{$\varrho(\alpha)$};
		\draw[thick] (-1,1) to [out=10,in=225] (3,3) -- (5,3);
		\draw[very thick, red, dashed] (-1,1) to [out=10,in=225] (3,3) -- (4,4);
		\end{tikzpicture}
		\caption{$d\geq3, \alpha\le 0; d=1,2, \alpha<0$}
		\label{fig:PMFdensity>=3}
	\end{subfigure}
	\caption{Total particle density of the zero of $I^\PMF_{\alpha,\mu}$. The limiting expected particle density (including the condensate) only differs for $\mu>a\varrho(\alpha)$, where it follows the dashed plot.}
	\label{fig:PMFdensity}
\end{figure}

In Figure~\ref{fig:PMFdensity} we observe in case (B) that for large values of $ \mu >a\varrho(\alpha) $ the density  $ D(\xi^\PMF) $ of the zero of the rate functions differs from the expected density which is represented by the dashed line. This signals a so-called condensate density which we will investigate further below in Section~\ref{sec-BEC}. For $d=1,2 $ and $ \alpha=0 $, we do not have any critical behaviour as the density $ D(\xi^\PMF) $ is a growing function of $ \mu $, see left hand side (A) of Figure~\ref{fig:PMFdensity}.

\medskip

\bigskip

We conclude the analysis of the rate function deriving the free energy function from the pressure function. 

\begin{prop}\label{P:fPMF}
	let $\beta<0$. For $\varrho >  0 $, the  free energy of the PMF model is defined as the Legendre-Fenchel transform of the pressure,
	\begin{equation}
	f^\PMF\left(\beta,\varrho\right):=\sup_{\mu\in\R,\alpha\leq0}\left\{(\mu+\alpha)\varrho - p^\PMF\left(\beta,\mu,\alpha\right)\right\}= f\left(\beta,\varrho\right)+ \frac{a}{2}\varrho^2.
	\end{equation}
\end{prop}

\begin{remark}
The free energy of the PMF model shows that the density square term in the definition of the measure stabilises the distribution and contributes towards the free energy. 
\hfill $ \diamond$
\end{remark}

\begin{remark}[\textbf{Conclusion.}]  The so-called BEC phase transition, see Appendix~\ref{appIBG},  is established in various equivalent ways, in Theorem~\ref{thm:PMFcondensate} below it is shown that the excess particle density is carried by so-called loops of unbounded length. Alternatively, Proposition~\ref{zeroPMF} and Proposition~\ref{P:fPMF} establish the phase transition via the change of the pressure density relation. The advantage of our LDP approach is that the rate function has unique zero and not an approximating sequence of minimiser. This is due to the fact that we are using the lower semicontinuous regularisation of the energy proving the large deviation principle. A close inspection of Figure~\ref{fig:PMFdensity}
reveals this. For $ d\ge 3 $ and $ \alpha\le 0$ or $ d\ge 1 $ and $ \alpha<0 $, we know that $ a\varrho(\alpha) <\infty $, and thus the density of the zero  of the rate function is constant for all $ \mu\ge a\varrho(\alpha) $. In this region, the total particle density is the dashed line intersecting the point $ (a\varrho(\alpha),\varrho(\alpha)) $. The so-called condensate density is then $ (\frac{\mu}{a}-\varrho(\alpha))_+  $. This will be confirmed in Theorem~\ref{thm:PMFcondensate} below.
\hfill $ \diamond $
\end{remark}

\subsubsection{HYL model}

The zeroes and the analysis for the HYL model are more complex and involved. We collect different statements, the main ones are Theorem~\ref{THM-uniqueness} and Theorem~\ref{THM-nonuniqueness} for uniqueness  and non-uniqueness of the rate function zeroes respectively. All proofs of this section are in Section~\ref{sec-thermHYL}.

\begin{prop}\label{zeroHYL} 
\begin{enumerate}[(a)]
\item The zeroes $ \left\{\xi^\HYL\right\}\subset\ell_1(\R_+) $ of rate the function $ I^\HYL_{\alpha,\mu} $ satisfy the following expression,
	\begin{equation*}
	\xi_k^\HYL = -\frac{1}{b\beta k^2}W_{\chi^*_k}\Big(-b\beta k^2 q_k^\ssup{\alpha} \exp\big[\beta k \big(\mu - a \delta^*\big)
	\begin{Bmatrix}
	1&:a\delta^* \geq \mu \\
	-\frac{b}{a-b}&:a\delta^* \leq \mu 
	\end{Bmatrix}
	\big]\Big), \qquad k\in\N,
	\end{equation*}
	where $\left(\delta^*,\chi^*\right) \in \R_+\times \left\{0,-1\right\}^\N$ is a solution to
	\begin{equation}
	\label{eqn:HYLconsistency}
	\delta = g^\chi\left(\delta\right) := -\frac{1}{b\beta} \sum^{\infty}_{k=1} \frac{1}{k}W_{\chi_k}\left(-b\beta k^2 q_k^\ssup{\alpha} \exp\left[\beta k \left(\mu - a \delta\right)
	\begin{Bmatrix}
	1&:a\delta \geq \mu \\
	-\frac{b}{a-b}&:a\delta \leq \mu 
	\end{Bmatrix}
	\right]\right)\,.
	\end{equation}
\item

For $ \mu\in\R,\alpha\leq 0 $ and $ a > b \geq 0 $ we the pressure is given by the rate function,

\begin{equation}
\begin{aligned}
p^\HYL(\beta,\mu,\alpha)&= \lim_{N\to\infty}\frac{1}{\beta\abs{\L_N}}\log Z_N^\HYL(\beta,\mu,\alpha) = -\inf_{x\in\ell_1(\R_+)}\big\{I_\alpha(x) + H^\HYL_{\mu,\lsc}(x)\big\}\,.
\label{p-hyl}		
\end{aligned}
\end{equation}

\end{enumerate}

\end{prop}

Note that our results hold for all parameter $ a>b $ wheres the ones in \cite{BLP} apply only to $ a=2b$, and we can dispense a technical assumption necessary in \cite{BLP}. Furthermore we are able to derive a full large deviation principle, whereas \cite{BDLP90} was only able to find the pressure for their related model. We shall study the relation between the derivatives of the pressure and the expected densities. As this models has more than one potential minimiser (zero) of the rate function, we first single out uniqueness parameter regimes.

\bigskip

The next proposition shows that there exist two regimes for the parameter  $ \mu $ depending on all other parameters such that the rate function has a unique zero. The zeroes of the rate function equal minimiser of the functional
\begin{equation}
F_\mu(x):=I_\alpha(x)+H^\HYL(x)-\frac{1}{2(b-a)}\big(\mu-aD(x)\big)^2_+\,.
\end{equation}

\begin{theorem}\label{THM-uniqueness}
For all $d\ge 1, \alpha\le 0 $ and $ \beta>0 $, there exist $\mu^+(d,\beta,\alpha,a,b)\ge \mu_-(d,\beta,\alpha,a,b) $ such that for $ \mu<\mu_- $ and $ \mu>\mu^+ $ the function $F_\mu$ has a unique minimiser, corresponding to the parameter $ (\delta,0) $ where $ \delta $ is the unique solution to the equation
$$
\delta=g^{\chi\equiv 0}(\delta)\,,
$$ 
where $ g^\chi $ is defined in \eqref{eqn:HYLconsistency}.
\end{theorem}

Theorem~\ref{THM-uniqueness} establishes parameter regimes for which the rate function has a  unique zero and thus no phase transition or  critical behaviour is present. 
In these `unique zero' regimes, we have an expression for the derivative of the thermodynamic pressure.
\begin{prop}
\begin{enumerate}[(a)]

\item     For the regime described in Theorem~\ref{THM-uniqueness}, i.e., for $ \mu< \mu_- $ and $ \mu>\mu^+ $, the pressure $p^\HYL\left(\beta,\alpha,\mu\right)$ is smooth and convex in $\mu$. In particular,
	\begin{align*}
		\frac{\d p^\HYL}{\d \mu} = \begin{cases}
		    D\left(\xi^\HYL\right) &, \mu<\mu_-\\
		    \frac{b}{a-b}\left(\frac{\mu}{b} - D\left(\xi^\HYL\right)\right) &, \mu > \mu_+,
		\end{cases} 
	\end{align*}
	for such $\mu$. 
	
\item The zeroes $\xi^\HYL$ in the uniqueness  regime 	$ \mu< \mu_- $ and $ \mu>\mu^+ $ are $\xi^\HYL\in\ell_1\left(\R_+\right)$ where
	\begin{equation*}
	\xi^\HYL_k=-\frac{1}{b\beta k^2}W_0\left(-b\beta k^2 q_k^\ssup{\alpha}\exp\left[\beta k \left(\mu  - a \delta^*\right)\right]\right),\qquad k\in\N,
	\end{equation*}
	and $ \delta^*=\delta^*(\beta,\mu+\alpha,a,b) $ is given implicitly as the unique solution to the equation $\delta^*=g^0\left(\delta^*\right)$,
	where $ g^0=g^{\chi\equiv 0} $.

\end{enumerate}
	
\end{prop}

The remaining crucial question is whether within the other  parameter regimes one can identify parameter values with multiple zeros of the rate function. This in turn signals critical behaviour and is of fundamental interest.

Before we present our non-uniqueness results in the next theorem we collect some facts about the possible solutions for $ \chi\equiv 0 $. Recall that then
\begin{equation}
\begin{aligned}\label{mingeneral}
\xi_k^\HYL&=-\frac{1}{b\beta k^2} W_0(-b\beta k^2q_k^\ssup{\alpha}\exp(\Escr_k(\delta))\big) \,,\quad\mbox {where } \delta \mbox{ solves} \\
\delta&=g^0(\delta)=\sum_{k\in\N} k\xi^\HYL_k\,,
\end{aligned}
\end{equation}
and
$$ 
\Escr_k(\delta)= \beta k \big(\mu - a \delta^*\big)
	\begin{Bmatrix}
	1&:a\delta^* \geq \mu \\
	-\frac{b}{a-b}&:a\delta^* \leq \mu 
	\end{Bmatrix} \le 0\, \mbox{ and } \Escr_k(\delta)\uparrow 0 \mbox{ as } \delta\uparrow \mu/a \mbox{ and } \Escr_k(\delta)\to -\infty \mbox{ as } \delta\to \infty\,.
$$
The solution in \eqref{mingeneral} is only well-defined as long as
\begin{equation}\label{argW}
0\ge -b\beta k^2q_k^\ssup{\alpha}\exp(\Escr_k(\delta))\ge -\e^{-1}\,,\quad \mbox{ for all }k\in\N\,.
\end{equation}
In case \eqref{argW} fails for some $k $ we  have no solution for that value of $\delta$, and hence no stationary point and zero. The condition \eqref{argW} is satisfied for every $\delta$ as long as 
\begin{equation} 
\beta\ge \beta^*:=\Big(\frac{(\e b)^2}{(4\pi)^d}\Big)^{\frac{1}{d-2}} \quad\text{ or equivalently } \quad b\le \frac{1}{\e\beta q_1^\ssup{\alpha}}=\frac{(4\pi\beta)^{d/2}}{\e\beta}=:b^*\,.
\end{equation}
The cases $d=1,2 $ are noteworthy: in $d=2 $ the upper bound for $b$ does not depend on $ \beta $ whereas in $d=1 $ that upper bound is a decreasing function of $ \beta $. 

\bigskip

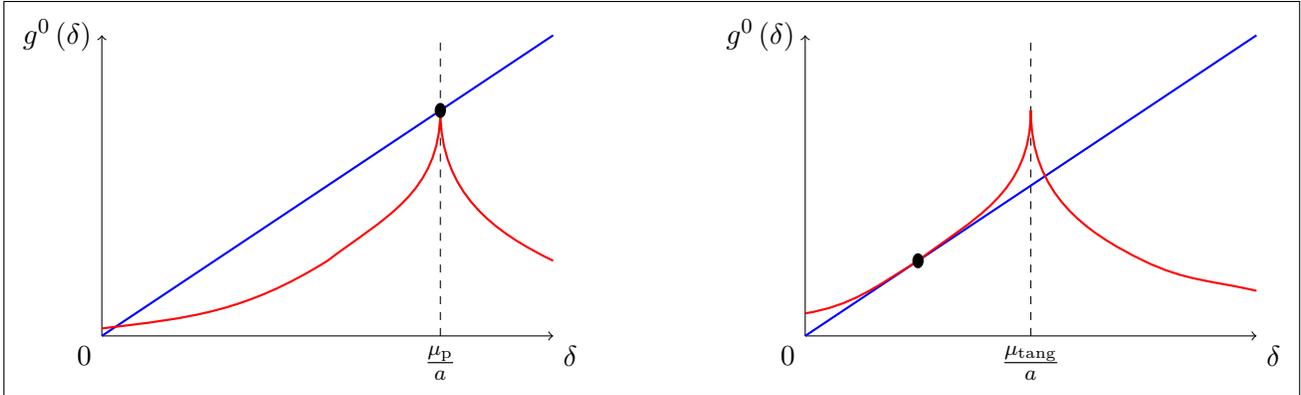
\begin{figure}
\centering
\begin{subfigure}[b]{0.45\textwidth}
	\centering
	\begin{tikzpicture}[yscale = 2,xscale=1.5]
	\draw[->] (0,0) node[below left]{$0$} -- (4,0) node[below right]{$\delta$};
	\draw[->] (0,0) -- (0,2) node[left]{$g^0\left(\delta\right)$};
	\draw[dashed] (3,0) node[below]{$\frac{\mu_{\rm p}}{a}$} -- (3,2);
	\draw[thick,blue] (0,0) -- (4,2);
	\draw[thick,red] (0,0.05) to [out=5,in=205] (2,0.5) to [out=30,in=270] (3,1.5) to [out=270,in=160] (4,0.5);
	\fill (3,1.5) circle (0.05);
	\end{tikzpicture}
\end{subfigure}
\hfill
\begin{subfigure}[b]{0.45\textwidth}
	\centering
	\begin{tikzpicture}[yscale = 2,xscale=1.5]
	\draw[->] (0,0) node[below left]{$0$} -- (4,0) node[below right]{$\delta$};
	\draw[->] (0,0) -- (0,2) node[left]{$g^0\left(\delta\right)$};
	\draw[dashed] (2,0) node[below]{$\frac{\mu_{\rm tang}}{a}$} -- (2,2);
	\draw[thick,blue] (0,0) -- (4,2);
	\draw[thick,red] (0,0.15) to [out=7,in=205] (1,0.5) to [out=30,in=270] (2,1.5) to [out=270,in=160] (3,0.5) to [out=340,in=170] (4,0.3);
	\fill (1,0.5) circle (0.05);
	\end{tikzpicture}
\end{subfigure}
	\caption{Sketch of $g^0\left(\delta\right)$ against the diagonal $g^0\left(\delta\right)=\delta$ for $d=3,4$, $\beta\geq \beta^*$ and $\mu=\mu_{\rm p} $ and $\mu= \mu_{\rm tang}<\mu_{\rm p}$. The chemical potential $\mu_{\rm p}$ is defined so that there is an intersection at the peak of $g^0\left(\delta\right)$, and $\mu_{\rm tang}$ is defined so that there is a tangential intersection. Note that $g^0\left(\delta\right)$ translates with $\mu$.}\label{fig:gHYL}
\end{figure}

In order to formulate the conditions for different  cases we introduce the parameters $ \mu_{\rm p} $ and $ \mu_{\rm tang} $. To define them, it will be convenient to define the function $\tilde{h}\colon(-\infty,0] \to \R_+$ as
\begin{equation}
    \tilde{h}(x) := -\frac{1}{b\beta} \sum^{\infty}_{k=1} \frac{1}{k}W_{\chi_k}\left(-b\beta k^2 q_k^\ssup{\alpha} \exp\left(\frac{ab}{a-b}\beta k x\right)\right).
\end{equation}
Note that $g^0(\delta) = \tilde{h}(\delta-\frac{\mu}{a})$ for $\delta \leq \frac{\mu}{a}$. We then define $\mu_{\rm p}$ as the chemical potential that produces a solution to \eqref{mingeneral} at the peak of $g^0(\delta)$ - represented pictorially in Figure~\ref{fig:gHYL}. This value is given by
\begin{equation}
\mu_{\rm p}:= a\tilde{h}(0).
\end{equation}
Then there is only an unique solution to $\delta = g^0(\delta)$ possible for any $ \mu>\mu_{\rm p} $, and this solution is for small values of $ \delta $. 
Note that $ \mu_{\rm p}=\mu_{\rm p}(b)  $ is increasing in $ b\le b^*$. For values of $ \mu<\mu_{\rm p} $, the peak wanders to the left and lies above the diagonal identity line. The blue line in Figure~\ref{fig:gHYL} may have up to three intersections with the graph $g^0 $ (red line). Now define
\begin{equation}
\mu_{\rm tang} := a\inf_{x\leq 0}\big\{\tilde{h}(x) - x\big\}\,.
\end{equation}
At $ \mu=\mu_{\rm tang} $ there is a tangential intersection whereas for $ \mu<\mu_{\rm tang} $ there is a unique intersection at some $ \delta^*>\mu/a $.  We see  from the construction that $ \mu_{\rm tang}\leq\mu_{\rm p} $. In particular, if $\lim_{x\uparrow 0}\tilde{h}'(x) > 1$ then the inequality is strict. This is always the case for $d=3,4$ which can be seen calculating the derivative. For the remaining cases $ d\ge 5 $ we need the assumption \eqref{dge5} below.

Figure~\ref{fig:gHYL} demonstrates how up to three possible minimisers can exist. There are three cases to distinguish. (i) The blue line intersects the red line only once. This happens for $ \mu<\mu_{\rm tang} $ at some $ \delta^*>\mu/a $ and for $ \mu>\mu_{\rm p} $ at some $ \delta^*<\mu/a $. This case is the unique regime in Theorem~\ref{THM-uniqueness} above. (ii) For $ \mu=\mu_{\rm tang} $ the blue line is tangent to $g^0 $ on the left hand side of the peak and intersects $g^0 $ to the right hand side of the peak. Likewise for $ \mu=\mu_{\rm p} $, the blue line intersects once to the left hand side of the peak and once to the right hand side of the peak.  (iii) For $ \mu\in(\mu_{\rm tang},\mu_{\rm p}) $, the blue line has two intersections with $g^0 $ to the left hand side of the peak and one with $g^0 $ to the right hand side of the peak. 

\bigskip

In the context of the following theorem, for $ d\ge 5 $ the following condition on the $\chi=0 $ solutions is required.
\begin{equation}\label{dge5}
\lim_{\delta\uparrow\frac{\mu}{a}} \big(g^0\big)^\prime(\delta) \equiv \lim_{x\uparrow0}h'(x) >1 \,\quad\mbox{ for } \mu\le\mu_{\rm p}\,.
\end{equation}
If this condition does not hold, then the behaviour is qualitatively different. In particular, the simultaneous minimiser behaviour in the following theorem does not occur.

\begin{theorem}\label{THM-nonuniqueness}
Let $d\ge 3, \beta>0, a>0 $ and $ 0<b<\min\big\{a,(\beta q_1^\ssup{\alpha})^{-1}\e^{-\beta\mu_{\rm p}/a}\big\} $. If $ d\ge 5$, suppose that \eqref{dge5} also holds.
Then there exist a unique $ \mu^*\in(\mu_{\rm tang},\mu_{\rm p}) $, such that $F_\mu $ has a unique minimiser for $ \mu<\mu^* $ and for $\mu>\mu^* $, and two simultaneous minimiser at $ \mu=\mu^* $. In particular, $\inf\{F_\mu\} $ is differentiable on $ (-\infty,\mu^*)\cup(\mu^*,\infty) $, and not on any open neighbourhood of $ \mu^* $.
The densities of the two simultaneous minimisers $ \xi^1(\mu^*) $ and $ \xi^2(\mu^*) $ satisfy $ D(\xi^1(\mu^*))>\mu^*/a > D(\xi^2(\mu^*)) $. 
\end{theorem}

\begin{remark}
This theorem proves our conjectures above. The critical value $ \mu^* $ for the two simultaneous minimiser is not know explicitly but we can see that for dimensions $ d\ge 3 $, the value decreases with increasing the parameter $ \beta $. If we increase the parameter $b$ of the negative counter terms in the measure (Hamiltonian) we obtain smaller values in the critical parameter as well. This behaviour of the system is been expected and our result establishes a proof for the whole parameter regime $ b\in[0,b^*) $.

\hfill $ \diamond $
\end{remark}

\bigskip

\bigskip

\subsection{Criticality and Bose-Einstein condensation}\label{sec-BEC}
The results in Section~\ref{sec-minimiser} signal three different types of criticality in our analysis. First the so-called BEC according to the reference measure and the CMF model, see Appendix~\ref{appIBG} and Section~\ref{sec-thermCMF}: The reference measure and the CMF model  show some critical behaviour, namely, the density and chemical potential relation breaks down in dimensions $ d\ge 3 $. The excess density above the critical value of the density is identified as the BEC condensate density. The free energy is constant for all densities above the critical one, showing that the condensate density does not contribute to the free energy. In \cite{BCMP05}, using the energy (Fourier) representation, it is shown that the excess density for the reference measure (for dimensions $ d\ge 3 $) equals the expected density of particles in the zero-energy mode. We believe that the same holds for our CMF model.

For the PMF model in Section~\ref{sec-PMFa} critical behaviour is the change in the derivative of the pressure-density relation, see Figure~\ref{fig:PMFdensity}, where the dashed line in (B) shows the density and signals that the excess density is given by $ (\frac{\mu}{a}-\varrho(\alpha))_+ $.  The third type of criticality is  established only for the HYL model and is given as the existence of two simultaneous  minimiser (zeroes) of the LDP rate function for a critical parameter $ \mu=\mu^* $, see Theorem~\ref{THM-nonuniqueness}. As one of the two densities is strictly smaller than the other one, the difference shall represent the so-called condensate density. However, we do not know if this excess density corresponds to the density of the so-called infinitely long loops. 

To investigate this further, we shall find an order parameter. In \cite{Gir60}, Girardeau suggests an order parameter for the (so-called \emph{generalised}) Bose-Einstein condensation phase transition which has been further studied by \cite{Lew86,BLP}. Having in mind the density of particles with zero single particle energy in the thermodynamic limit, they first take the finite volume expected density of particles with energy below some cut-off. Then the thermodynamic limit is taken before finally the cut-off goes to zero. In contrast to these momentum-space perspectives, we do not keep track of the particles' energy. Instead, we partition our gas by loop type, and expect the condensate to occupy loops of diverging length. Therefore we want to evaluate the `condensate density' given by
	\begin{equation*}
	\Delta\left(\beta,\alpha\right) := \lim_{K\to\infty}\lim_{N\to\infty} \E_{\nu_{N,\alpha}}\left[D-D_K\right],\qquad D_K\left(x\right):=\sum^K_{k=1}kx_k\,.
	\end{equation*}

\bigskip

We analyse first the reference measure (ideal gas) and the CMF model in the following theorem.

\begin{theorem}\label{thm:IDEALcondensate}
The following statements hold for the reference measure and the CMF model.
For $\beta>0$, $\alpha<0$ $a\ge 0$ , we have $\Delta\left(\beta,\alpha\right) =\Delta^\CMF(\beta,\alpha) = 0$, and for $\beta>0$, $\alpha=0$, $a\ge 0 $
		\begin{equation*}
		\Delta\left(\beta,0\right) =\Delta^\CMF(\beta,0)
		=\begin{cases} +\infty & ,\mbox{ for } d=1,2\,,\\ 0 & ,\mbox{ for }  d\ge 3\,.\end{cases}		\end{equation*}
	\end{theorem}

Mathematically, the critical behaviour is seen in the thermodynamic limit of distribution function of the empirical particle density (called Kac distribution), the limit of the Kac distribution  exists only for densities $ \varrho<\varrho_{\rm c} $ and is the degenerate distribution 
$$
\mathbb{K}(x)=\begin{cases} 0 &, x <\varrho, \\ 1  & , x\ge \varrho.\end{cases}
$$
This critical behaviour is also seen in the rate function for the density large deviation principle in Theorem~\ref{THM-densityLDP}. 	We believe that the unexpected results in Theorem~\ref{thm:IDEALcondensate} are due to this degeneracy.
	
\bigskip

The next theorem confirms the conjecture that$ (\mu/a-\varrho(\alpha))_+ $ represents the condensate density. 

\begin{theorem}\label{thm:PMFcondensate}
	For all $\beta>0$, $\mu\in\R$ and $ \alpha\leq 0 $,
	\begin{equation*}
	\begin{aligned}
		\Delta^\PMF\left(\beta,\alpha,\mu\right) &= \Big(\frac{\mu}{a} - \varrho(\alpha)\Big)_+ = \Big(\frac{\partial}{\partial \mu}\left(H^\PMF - H^\PMF_{\mu,\lsc}\right)\Big)\left(\xi^\PMF\right),
		\end{aligned}
	\end{equation*}
	where $\xi^\PMF$ is the unique minimiser (zero) of the rate function $I^\PMF_{\alpha,\mu}$. Furthermore,
	$$
	\lim_{\alpha\uparrow 0} \Delta^\PMF(\beta,\alpha,\mu)=\big(\frac{\mu}{a}-\varrho_{\rm c}(d)\big)_+=\begin{cases} 0 &,\mbox{ for } d=1,2\,,\\ > &,\mbox{ for } d\ge 3, \mu>a\varrho_{\rm c}(d)\,.\end{cases}
	$$
\end{theorem}

\bigskip

The following result shows that the condensate density has a limit for certain regimes of thermodynamic parameter $ \beta $ and $ \mu $ and energy parameter $ a $ and $ b$. Explicit expressions for the minimiser and condensate density in the critical regime are not available.
	
	\begin{theorem}\label{THM:HYL-Cond}
	For $\beta>0$, $\mu\in\R$, $\alpha\leq0$, where the derivative is defined,
	\begin{equation}\label{HYLcondenstate1}
	\Delta^\HYL(\beta,\mu,\alpha) = -\lim_{K\to\infty}\frac{\d}{\d s}\Big(\inf_{\ell_1\left(\R_+\right)}\left\{I^\HYL_{\alpha,\mu+s} + s D_K\right\} - p^\HYL\left(\beta,\alpha,\mu+s,\right)\Big)\Big|_{s=0}.
	\end{equation}
	In particular, if the infimum is achieved by $\xi\left(s\right)\in C^1\left(\left(-\epsilon,\epsilon\right):\ell_1\left(\R_+\right)\right)$ for some $\epsilon>0$, then
	\begin{equation}\label{HYLcondenstate2}
	\Delta^\HYL\left(\beta,\alpha,\mu\right) = \frac{a}{a-b}\Big(\frac{\mu}{a} - D\left(\xi\left(0\right)\right)\Big)_+ = \Big(\frac{\partial}{\partial \mu}\left(H^\HYL_{\mu} - H^\HYL_{\mu,\lsc}\right)\Big)\left(\xi\left(0\right)\right).
	\end{equation}
\end{theorem}

\begin{remark}[\textbf{Conclusion HYL model}]  The BEC phase transition for the HYL model  is established as follows. Proposition~\ref{HYLboringcase} establishes a subcritical regime for which the pressure is smooth and its derivative gives the particle density with no condensation.    
Depending on the density $ D(\xi(0)) $, for large enough $ \mu $ the particle density in loops of unbounded length is not vanishing, see Theorem~\ref{THM:HYL-Cond}. In Figure~\ref{fig:HYL3sols} we can identify the regime $ \mu \ge \mu_{\rm p} $ when the density of the zero of the rate function is decreasing with increasing $ \mu $ such that the excess density is carried by loops of unbounded length. According to our Theorem~\ref{THM-nonuniqueness} at criticality $ \mu=\mu^* $ we have two zeroes $ \xi^1(\mu^*) $ and $ \xi^2(\mu^*) $ of the rate function with $ D(\xi^1(\mu^*))>\mu^*/a> D(\xi^2(\mu^*)) $. Using \eqref{HYLcondenstate2} in Theorem~\ref{THM:HYL-Cond}, we can see that $ \Delta^\HYL(\beta,\alpha,\mu^*)=\frac{a}{a-b}\big(\frac{\mu^*}{a}-D(\xi^2(\mu^*))\big)_+>0 $ but we are unable to establish that the difference of the densities of the two simultaneous minimiser equals $\Delta^\HYL(\beta,\alpha,\mu^*)$.

Finally, if we choose $ a=2b $ in \eqref{HYLcondenstate2}, we can recover the results in \cite{Lew86} and \cite{BLP}. 
It shows in fact, that for increasing values of the coupling parameter $ a $, the condensate density decreases. On the other hand, if the parameter for the counter energy term, $ b $, is approaching $ a $, the condensate density increases. This is due to the fact that with large counter terms the system distributes the physical particles in as few as possible different cycles lengths. To accommodate the particle density, the only way is to put them in infinitely long cycles.  Our analysis actually shows that the BEC phase transition for HYL is more complex and requires further detailed study. For example, at criticality $ \mu=\mu^* $ one would like to study concentration of the underlying measure around the two distant rate function zeroes. 

\hfill $ \diamond $
\end{remark}

\section{Proof of the Large deviations principles}\label{sec-LDPproofs}
This section contains the proofs for all large deviation principles.  We use the large deviation principle of the reference measure in Appendix~\ref{appIBG}  as a stepping stone towards arriving at LDPs for our models. Section~\ref{sec-CMF} and Section~\ref{sec-PMF} give the proofs for  both mean-field  models, i.e, CMF model  and PMF model respectively, and Section~\ref{sec-HYL} the proof for the HYL model.

\subsection{Proof of Theorem~\ref{THM-CMF} - Cycle Mean Field LDP}\label{sec-CMF}
We are going to apply Varadahn's Lemma in \cite[Theorem~4.3.1]{DZ09}.  To show continuity of $H^{\CMF}$, let us show sequential continuity. This implies continuity since $\ell_1(\R)$ is a metric space. Let $x^{(n)}\rightarrow x$ be a convergent sequence in $\ell_1(\R)$, so $\lim_{n\rightarrow\infty}\sum_{j\in\N}\abs{x^{(n)}_k-x_k}=0$.
 Let $S(x):=\sum_{k\in\N}x_k$. Then
	\begin{equation*}
	\lim_{n\rightarrow\infty}\abs{S(x^{(n)})-S(x)}= \lim_{n\rightarrow\infty}\big|\sum_{k\in\N}x^{(n)}_k - \sum_{k\in\N}x_k\big| \leq \lim_{n\rightarrow\infty}\sum_{k\in\N}\abs{x^{(n)}_k-x_k} =0.
	\end{equation*}
	Hence $S$ is continuous. We can then write the Hamiltonian as the composition of continuous functions $H^{\CMF} = T\circ S$, where $T\colon \R\rightarrow\R$, $x\rightarrow \frac{a}{2}x^2$.
We can now simply apply Varadhan's Lemma. The lower bound
\begin{equation}\label{CMFlower}
\liminf_{N\to\infty}\frac{1}{\beta \abs{\L_N}}\log \E_{\nu_{N,\alpha}}\big[\ex^{-\beta\abs{\L_N} H^\CMF}\big]\ge \sup_{x\in\ell_1(\R)}\big\{- H^\CMF(x)-I_\alpha(x)\big\}
\end{equation}
follows easily with \cite[Lemma~4.3.4]{DZ09}.  For the corresponding upper bound we simply note that the tail-condition in \cite[Theorem~4.3.1]{DZ09} holds due to $ H^\CMF(x)\ge 0 $ for all $ x\in\ell_1(\R_+) $. Therefore, with \cite[Lemma~4.3.6]{DZ09} we obtain the corresponding upper bound
\begin{equation}\label{CMFupper}
\limsup_{N\to\infty}\frac{1}{\beta \abs{\L_N}}\log \E_{\nu_{N,\alpha}}\big[\ex^{-\beta\abs{\L_N} H^\CMF}\big]\le \sup_{x\in\ell_1(\R)}\big\{- H^\CMF(x)-I_\alpha(x)\big\}.
\end{equation}
We conclude with the statement in Theorem~\ref{THM-CMF} by combining the lower bound \eqref{CMFlower} and the upper bound \eqref{CMFupper}.\qed

\subsection{Proof of Theorem~\ref{THM-PMF} - Particle Mean Field LDP}\label{sec-PMF}
 To prove the large deviation principle for $ \nu_{N,\mu,\alpha}^\PMF $ one would simply use Varadhan's Lemma. However, the first term of the Hamiltonian $ H^\PMF_\mu $ for $ \mu>0 $ is not lower semicontinuous whereas the second term is only lower semicontinuous.  Using \cite{GZ93}, one would arrive at lower and upper bounds for $ \E_{\nu_{N,\alpha}}[\ex^{-\abs{\L_N}\beta H^\PMF_\mu}] $ using the upper and the lower semicontinuous regularisation of $ H^\PMF_\mu $, respectively.    Unfortunately, the upper semicontinuous regularisation of the Hamiltonian equals infinity, and thus it does not provide a lower bound for the large deviation principle. Our strategy is therefore twofold. For the large deviation upper bound we use the lower semicontinuous regularisation in conjunction with the corresponding bound in Varadhan's Lemma.  We obtain the corresponding large deviation lower bound by conditioning that the empirical cycle count is supported on a finite-dimensional subspace. On this event we can replace our measure by the corresponding measure with finite dimensional mark space. On this subspace the Hamiltonian is in fact continuous and thus application of Varadhan's Lemma provides a lower bound. To remove the cutoff parameter we will construct finite-dimensional sequences approximating the infimum of the corresponding lower bound.
We start with a couple of observations.

\begin{lemma}
For all $ \mu>0  $, the lower semicontinuous regularisation of $ H^\PMF_\mu $ is given as
\begin{equation}\label{lsc}
H^{\PMF}_{\mu,\lsc}(x)=H^\PMF_\mu(x)-\frac{1}{2a}\left(\mu-aD(x)\right)_+^2=\begin{cases} -\mu D(x)+\frac{a}{2}D(x)^2 &,D(x)\ge \frac{\mu}{a},\\-\frac{\mu^2}{2a}&, D(x)<\frac{\mu}{a}, \end{cases}\quad x\in\ell_1(\R_+), 
\end{equation}
whereas for all $ \mu\le 0 $, $ H^\PMF_{\mu,\lsc}\equiv H^\PMF_\mu $.
\end{lemma}

\begin{proofsect}{Proof}
Suppose $h$ is lower semicontinuous and $h\left(x\right)\leq H^\PMF_\mu\left(x\right)$ for all $x\in\ell_1\left(\R_+\right)$. Let $x^\ssup{n}\to x$ in $\ell_1\left(\R_+\right)$. Then $h\left(x\right)\leq\liminf_{n\to\infty}h\left(x^\ssup{n}\right) \leq \liminf_{n\to\infty}H^\PMF_\mu\left(x^\ssup{n}\right)$. Let $x^{\ssup{n}}_k=x_k+\frac{1}{na}\left(\mu-aD(x)\right)_+\1\{n=k\}$, so this inequality implies that $h\left(x\right)\leq H^{\PMF}_{\mu,\lsc}\left(x\right)$.

Note that $H^{\PMF}_{\mu,\lsc}\left(x\right) = g\circ D\left(x\right)$, where $g\colon D\mapsto -\mu D + \frac{a}{2}D^2 - \frac{1}{2a}\left(\mu - a D\right)^2_+$ for $0\leq D<+\infty$ and $+\infty\mapsto+\infty$. Since $g$ is continuous and non-decreasing (see Figure~\ref{fig:PMF_lsc}), and $D\left(x\right)$ is lower semicontinuous, the composition is also lower semicontinuous.
\qed
\end{proofsect}

\begin{prop}[\textbf{Upper  bound PMF model}]
For all $ \mu\in\R, \alpha\le 0 $, and $ a>0 $,
\begin{equation}\label{ubound}
\limsup_{N\to\infty}\frac{1}{\beta\abs{\L_N}}\log \E_{\nu_{N,\alpha,\mu}^\PMF}\big[\ex^{-\beta\abs{\L_N} H^\PMF_\mu}\big]\le -\inf_{x\in\ell_1(\R_+)}\big\{ I_\alpha(x)+ H^\PMF_{\mu,\lsc}(x)\big\}\,.
\end{equation}
\end{prop}

\begin{proofsect}{Proof}
The statement follows easily with the upper bound estimate in Varadhan's Lemma given in \cite[Lemma~4.3.6]{DZ09} using  the inequality $ H^\PMF_\mu(x)\ge H^\PMF_{\mu,\lsc}(x)  $ for all $ x\in\ell_1(\R_+) $, the lower semicontinuity of $ H^\PMF_{\mu,\lsc} $,  and the fact that $ H^\PMF_{\mu,\lsc}(x)\ge -\frac{\mu^2}{2a} $. The later estimate provides the tail-condition necessary to apply \cite[Lemma~4.3.6]{DZ09}.
\qed
\end{proofsect}

\begin{prop}[\textbf{Lower  bound PMF model}]\label{PMFlower}
For all $ \mu\in\R, \alpha\le 0 $, and $ a>0 $,
\begin{equation}\label{lbound}
\liminf_{N\to\infty}\frac{1}{\beta \abs{\L_N}}\log \E_{\nu_{N,\alpha,\mu}^\PMF}\big[\ex^{-\beta\abs{\L_N} H^\PMF_\mu}\big]\ge -\inf_{x\in\ell_1(\R_+)}\big\{ I_\alpha(x)+H^\PMF_{\mu,\lsc}(x)\big\}\,.
\end{equation}
\end{prop}

\begin{proofsect}{Proof}
The strategy for proving the lower bound is to first introduce a cut-off parameter as done in \cite{ACK}, that is, we change the measure to obtain a finite-dimensional problem which gives continuity of the Hamiltonian and thus a large deviation lower bound for the finite-dimensional space. The final step is then to remove the cut-off parameter. As our Hamiltonian is not positive,  removing of the cut-off is not as straightforward as in \cite{ACK}.  We thus need to construct a sequence for the finite dimensional spaces which allows for energy estimates and at the same time gives convergence towards the lower semicontinuous regularisation. The last step is crucial as  a  lower bound can be  obtained via the upper semicontinuous regularisation which in our case is identical to infinity.

\medskip

\textbf{Step 1: Restriction of the mark space.}  We will approximate the infinite-dimensional space $ \ell_1(\R_+) $ by finite-dimensional spaces. Pick some $ K\in\ N $, and the corresponding measure on $ \R_+^K $ which is isomorphic to $ \pi_K(\ell_1(\R_+)) $ with $ \pi_K\colon\ell_1(\R_+)\to\R_+^K,x\mapsto (x_1,\ldots,x_K) $, is denoted 
$$ 
\nu_{N,\alpha}^{\ssup{K}}=\nu_{N,\alpha}\circ\pi_K^{-1}\,.
$$   We obtain a large deviation principle for the cut-off version in the following.

\begin{lemma}
	\label{Thm:LDPidealRJ}
	For given $K\in\N $ and $\alpha\leq0$, the sequence $(\nu_{N,\alpha}^{\ssup{K}})_{N\in\N}$ satisfies a LDP on $\R_+^K$ with rate $\beta \abs{\Lambda_N}$ and rate function
	\begin{equation}\label{RFK}
	I^{\ssup{K}}_\alpha(x) = \sum^K_{k=1}\Big[\frac{x_k}{\beta}\Big(\log\frac{x_k}{q_k^{\ssup{\alpha}}} - 1\Big)+ \frac{q_k^{\ssup{\alpha}}}{\beta}\Big]\,.
	\end{equation}
	
\end{lemma}

\begin{proofsect}{Proof}
	Since the projection $\pi_K$ is continuous, we can apply the contraction principle to obtain a variational form of the rate function
	\begin{equation*}
	I^{\ssup{K}}_\alpha(x) = \inf_{\widetilde{x}\in\ell_1(\R_+)\colon\pi_K(\widetilde{x})=x}  I_\alpha(\widetilde{x})\,,
	\end{equation*}
	where $ I_\alpha $ is the rate function for $\nu_{N,\alpha}$ in Proposition~\ref{THM-Ideal}. We conclude using that projection is continuous. 
%
\qed
\end{proofsect}

\medskip

\textbf{Step 2: Lower bound.}   
 We obtain a lower bound 
 $$
 \E_{\nu_{N,\alpha}}\big[\ex^{-\beta\abs{\L_N} H^\PMF_\mu}\big]\ge \E_{\nu_{N,\alpha}}\big[\ex^{-\beta\abs{\L_N} H^\PMF_\mu}\1\{\blambda^\ssup{N}\in\R_+^K\}\big]\,, 
 $$
where we identified $ \R_+^K $ with the corresponding subspace in $ \ell_1(\R_+) $. On that event we can replace $ H^\PMF_\mu $ by 
$$
 H_{\mu,K}^\PMF(x)=-\mu\sum_{k=1}^K k x_k +\frac{a}{2}\Big(\sum_{k=1}^K k x_k\Big)^2\,,
 $$ 
and $ \E_{\nu_{N,\alpha}} $ by $ \E_{\nu_{N,\alpha}^{\ssup{K}}} $ - up to a factor of $\nu_{N,\alpha}\left(\blambda^\ssup{N}\in\R_+^K\right) = \exp\left(-\abs{\L_N}\sum^\infty_{k=K+1}q_k^{\ssup{\alpha}}\right)$. The finite-dimensional approximation $ H^\PMF_{\mu,K} $ is in fact continuous, and thus we obtain a large deviation lower bound  using Lemma~\ref{Thm:LDPidealRJ} and Varadhan's Lemma, see \cite[Lemma~4.3.4]{DZ09},
\begin{equation}
\begin{aligned}
\liminf_{N\to\infty}\frac{1}{\beta\abs{\L_N}}\log \E_{\nu_{N,\alpha}}\big[\ex^{-\beta\abs{\L_N} H^\PMF_\mu}\big] \ge -\inf_{x\in\R_+^K}  \big\{ I^{\ssup{K}}_\alpha(x)+ H^\PMF_{\mu,K}(x) \big\} - \sum^\infty_{k=K+1}q_k^{\ssup{\alpha}}\,.
\end{aligned}
\end{equation}

\medskip

\textbf{Step 3: Removing the cut-off parameter.} We are left to remove the cut-off by taking $ K\to\infty $ and to prove that the $K\to\infty $ limit of $ H^\PMF_{\mu,K} $ is replaced by the lower semicontinuous regularisation of $ H^\PMF_\mu $. The sum vanishes in the $K\to\infty$ limit because it converges.
\begin{lemma}
\begin{equation}\label{removing}
\limsup_{K\to\infty} \inf_{x\in\R_+^K}\big\{ I^{\ssup{K}}_\alpha(x)+H^\PMF_{\mu,K}(x)\big\} \le \inf_{x\in\ell_1(\R_+)}\big\{ I_\alpha(x)+ H^\PMF_{\mu,\lsc}(x)\big\}\,.
\end{equation}
\end{lemma}
\begin{proofsect}{Proof}
Fix $ \widetilde{x}\in \ell_1(\R_+) $ satisfying $ I_\alpha(\widetilde{x})+ H^\PMF_{\mu}(\widetilde{x}) <\infty $. For $ K\in\N$, consider $ x^K=\pi_K(\widetilde{x}) $. By \eqref{RFK}, we have $ I^{\ssup{K}}_\alpha(x^K)\le I_\alpha(\widetilde{x}) $. We shall replace  $ x^K $ by $ \widehat{x}^K $, defined as  
$$
\widehat{x}^K_k=\begin{cases}  x^K_k, &, k=1,\ldots, K-1,\\
x^K_K +\frac{1}{K}\left(\frac{\mu}{a}-D(\widetilde{x})\right)_+, &, k=K.\end{cases}
$$
Clearly, $ \norm{x^K-\widehat{x}^K}_{\ell_1} \to 0 $ as $ K\to\infty $, and $ x^K\to\widetilde{x} $ as $K\to\infty $. Furthermore, 
$$
I^{\ssup{K}}_\alpha(\widehat{x}^K)\le I^{\ssup{K}}_\alpha(x^K)+\abs{I^{\ssup{K}}_\alpha(\widehat{x}^K)-I^{\ssup{K}}_\alpha(x^K)}\le I_\alpha(\widetilde{x}) + \bigO(\frac{1}{K}\log K)\,.
$$ We turn to the energy term which needs extra care as the Hamiltonian is not positive. Observe that
$$
D(\widehat{x}^K) =\begin{cases}  \frac{\mu}{a} +D(x^K)-D(\widetilde{x}) &, D(\widehat{x}) <\frac{\mu}{a}\,,\\ D(x^K) &, D(\widehat{x})\ge \frac{\mu}{a}\,. \end{cases} 
$$
Assume that $ D(\widetilde{x}) \ge \frac{\mu}{a} $, so 
$$H^\PMF_{\mu,K}(\widehat{x}^K)=-\mu\sum_{k=1}^K k x^K_k +\frac{a}{2}\Big(\sum_{k=1}^Kk x^K_k\Big)^2\le  H^\PMF_{\mu,\lsc}(\widetilde{x})\,.$$ In the other case, $ D(\widetilde{x})<\frac{\mu}{a} $, for every  $ \eps>0 $ choose $ K$ sufficiently large such that $ \abs{D(x^K)-D(\widetilde{x})}<\eps $, and estimate
$$
\begin{aligned}
H^\PMF_{\mu,K}(\widehat{x}^K)&=-\mu\Big(\frac{\mu}{a}+D(x^K)-D(\widetilde{x})\Big)+\frac{a}{2}\Big(\frac{\mu}{a}+D(x^K)-D(\widetilde{x})\Big)^2\\
& \le H^\PMF_{\mu,\lsc}(\widetilde{x}) +2\mu\eps+\eps^2\frac{a}{2}
\end{aligned}
$$
to conclude with  \eqref{removing}.
\qed
\end{proofsect}

We finally combine Proposition~\ref{ubound} and Proposition~\ref{PMFlower}  to finish the proof for Theorem~\ref{THM-PMF}.

\qed
\end{proofsect}

\subsection{Proof of Theorem~\ref{Thm:HYL}} \label{sec-HYL}
This section proves Theorem~\ref{Thm:HYL} using  techniques which are similar to the ones in the  proof in Section~\ref{sec-PMF}. However, there are significant differences  to address due to the fact that the Hamiltonian $ H^\HYL_\mu $ has positive and negative contributions. We rewrite the Hamiltonian in two equivalent ways for any $ a\ge b>0 $ and $ \mu\in\R $,
\begin{align}
H^\HYL_\mu(x)&=-\mu\sum_{k=1}^\infty k x_k +\frac{(a-b)}{2}\Big(\sum_{k=1}^\infty k x_k\Big)^2+\frac{b}{2}\sum_{\heap{j,k=1}{j\not= k}}^\infty jkx_jx_k\label{version1}\\
&=-\mu\sum_{k=1}^\infty k x_k +\frac{a}{2}\Big(\sum_{k=1}^\infty k x_k\Big)^2-\frac{b}{2}\sum_{k=1}^\infty k^2x_k^2.\label{version2}
\end{align}
Note that that the right hand side in \eqref{version1} is the sum of a PMF Hamiltonian with interaction strength $ (a-b) $ and a lower semicontinuous and non-negative term. On the other hand, \eqref{version2}  expresses $ H^\HYL $ as the sum of a PMF Hamiltonian, and an upper semicontinuous and non-positive term. Let us introduce the following notations
\begin{align*}
H_{\mu,a}^\PMF(x)=-\mu\sum_{k=1}^\infty k x_k +\frac{a}{2}\Big(\sum_{k=1}^\infty k x_k\Big)^2, \quad a>0,\\
H_+(x)=\frac{b}{2}\sum_{\heap{j,k=1}{j\not= k}}^\infty jkx_jx_k,\qquad
H_-(x)= -\frac{b}{2}\sum_{k=1}^\infty k^2x_k^2.
\end{align*}

Thus $ H^\HYL_\mu = H_{\mu,a}^\PMF + H_-=H^\PMF_{\mu,a-b}+H_+$.
\begin{lemma}
For $ b>0 $, $H_- $ is upper semicontinuous and $H_+ $ is lower semicontinuous on $ \ell_1(\R_+) $. 
\end{lemma} 
\begin{proofsect}{Proof}
We shall show that $\sum_{\heap{j,k=1}{j\not= k}}^\infty jkx_jx_k $ and $  \sum_{k=1}^\infty k^2x_k^2 $ are both lower semicontinuous on $ \ell_1(\R_+) $. Suppose $ x^{\ssup{n}}\to x $ in $ \ell_1(\R_+) $. Clearly, $\abs{x^{\ssup{n}}_k-x_k}\le \norm{x^{\ssup{n}}-x}_{\ell_1}\quad\mbox{ for all }k\in\N$. Furthermore, due to the $ \ell_1 $-convergence and $ \ell_1\subset\ell_\infty $,  the term $ \abs{x^{\ssup{n}}_k-x_k} $ is bounded in both $k$ and $ n$. Hence
$$
\abs{(x^{\ssup{n}}_k)^2-x_k^2}=\abs{x^{\ssup{n}}_k-x_k}\abs{x^{\ssup{n}}_k+x_k}\to 0\; \mbox{ as } n\to\infty.
$$
Applying Fatou's Lemma here proves that $ H_- $ is upper semicontinuous. Similarly, for all $ (j,k)\in\N^2 $,
$$
\abs{x_j^{\ssup{n}}x_k^{\ssup{k}}-x_jx_k}\le \abs{x^{\ssup{n}}_j}\abs{x^{\ssup{n}}_k-x_k}+\abs{x_k}\abs{x^{\ssup{n}}_j-x_j}\to 0 \; \mbox{ as } n\to\infty.
$$
This convergence in conjunction with Fatou's Lemma shows that $ H_+$ is lower semicontinuous.
\qed
\end{proofsect}

\begin{lemma}\label{L:regHYL}
For $ a>b $, the $ \ell_1(\R_+) $ lower semicontinuous regularisation of $ H^\HYL_\mu $ is given by
\begin{equation}\label{lscHYL}
H^\HYL_{\mu,\lsc}(x)=H^\HYL_\mu(x) -\frac{\left(\mu-aD(x)\right)_+^2}{2(a-b)},\quad x\in\ell_1(\R_+).
\end{equation}
\end{lemma}

\begin{proofsect}{Proof}
Denote the right hand side of \eqref{lscHYL} by $ h$. Clearly, $ h(x)\le H^\HYL_\mu(x) $ and $ H^\PMF_{\mu,(a-b),\lsc} (x) \le h(x)=H^\PMF_{\mu,(a-b),\lsc}(x)+H_+(x)\le  H^\HYL_\mu(x) $ for all $ x\in\ell_1(\R_+) $. We need to show that $ h $ is the greatest lower semicontinuous function less or equal to $ H^\HYL_\mu $.

Suppose that $ x\in\ell_1(\R_+) $ with $ D(x)=\infty $. Then since $ H^\PMF_{\mu,(a-b),\lsc}(x) =\infty $, we have $ h(x)=\infty $. Suppose now that $ x\in\ell_1(\R_+) $ with $ D(x) <\infty $. For any sequence $ (x_n)_{n\in\N} $ with $ x_n\to x $ as  $ n\to\infty $ there exists $ (\eps^{\ssup{n}})_{n\in\N}\subset\ell_1(\R) $ such that $ x_n=x+\eps^{\ssup{n}} $, $ \eps^{\ssup{n}}\to 0 $ as $ n\to\infty $,  and $ x+\eps^{\ssup{n}}\in\ell_1(\R_+) $. Furthermore, $ \liminf_{n\to\infty}\left(D(x^{\ssup{n}})-D(x)\right)\ge 0 $ and thus $ \liminf_{n\to\infty}D(\eps^{\ssup{n}})\ge 0$. We show that $ h $ is lower semicontinuous by proving that 
\begin{multline}
\liminf_{n\to\infty} \left(h(x^{\ssup{n}})-h(x)\right)=\liminf_{n\to\infty}\Big\{-\mu\left(D(x^{\ssup{n}})-D(x)\right)+\frac{a}{2}\left(D(x^{\ssup{n}})^2-D(x)^2\right)\\
-\frac{b}{2}\Big(\sum_{k=1}^\infty k^2\left((x_k^{\ssup{n}})^2-x_k^2\right)\Big)+\frac{1}{2(a-b)}\Big(\left(\mu-aD(x)\right)_+^2 -\left(\mu-aD(x^{\ssup{n}})\right)_+^2\Big)\Big\}\ge 0 \label{hlsc}.
\end{multline}     
If $D\left(\varepsilon^\ssup{n}\right)\to+\infty$, then $\liminf_{n\to\infty}h\left(x^\ssup{n}\right) = +\infty$, so we suppose that $D\left(\varepsilon^\ssup{n}\right)$ is finite and bounded.

We write $ \eps_k^{\ssup{n}}=\eps^{+\ssup{n}}_k -\eps^{-\ssup{n}}_k $ with $ \eps^{+\ssup{n}}_k,\eps^{-\ssup{n}}_k\ge 0 $. Clearly, $ \eps^{-\ssup{n}}_k\le x_k $ for all $ k\in\N $.
We shall show both
\begin{align}
    \label{limsupzero}
\limsup_{n\to\infty}\sum_{k=1}^\infty k^2 x_k\eps^{\ssup{n}}_k &= 0,\\
\label{upperbmixed}
\limsup_{n\to\infty} \sum_{k=1}^\infty k^2(x_k+\eps^{\ssup{n}}_k)^2 &\le \sum_{k=1}^\infty k^2x_k^2 +\limsup_{n\to\infty} \left(D(\eps^{\ssup{n}})\right)^2.
\end{align}
As $ D(x) <\infty $, for all $ C>0 $ there exits $K_C\in\N$ such that $kx_k<C$ for all $k>K_C$. Therefore we have $\sum^\infty_{k=K_C+1}k^2x_k\varepsilon^\ssup{n}_k < C\sum^\infty_{k=K_C+1}k\varepsilon^\ssup{n}_k \leq CD\left(\varepsilon^\ssup{n}\right)$. Then since $\lim_{n\to\infty}\sum^{K_C}_{k=1}k^2x_k\varepsilon^\ssup{n}_k = 0$, and $D\left(\varepsilon^\ssup{n}\right)$ is bounded, we can choose $C$ arbitrarily small to get \eqref{limsupzero}. To obtain \eqref{upperbmixed} we just expand
\begin{equation}
\sum_{k=1}^\infty k^2(x_k+\eps^{\ssup{n}}_k)^2\le \sum_{k=1}^\infty k^2 x_k^2+2\sum_{k=1}^\infty k^2(x_k-\eps^{-\ssup{n}}_k)\eps^{+\ssup{n}}_k + \left(D(\eps^{+\ssup{n}})\right)^2.
\end{equation}
The middle term vanishes due to \eqref{limsupzero}. To show that 
\begin{equation}
    \limsup_{n\to\infty}\left(D(\eps^{+\ssup{n}})\right)^2\le \limsup_{n\to\infty}\left(D(\eps^{\ssup{n}})\right)^2,
\end{equation} 
note that $ D(x)<\infty $ implies that $D(\eps^{\ssup{n}})-D(\eps^{+\ssup{n}})=D(\eps^{-\ssup{n}})\le D(x)<\infty$. Hence, for any $ \delta $ there exists $ K\in\N $ such that $\sum_{k=K+1}^\infty k\eps^{-\ssup{n}}_k\le\sum_{k=K+1}^\infty k x_k <\frac{\delta}{2}$. On  the other hand, for this $ \delta $ there exists a $ n(K)\in\N $ such that
\begin{equation}
    \sum_{k=1}^K k\eps^{-\ssup{n}}_k <\frac{\delta}{2},\quad \mbox{ for all }  n\ge n(K),
\end{equation}
thus showing \eqref{upperbmixed}. We continue with
\begin{multline}\label{hlsc2}
\mbox{ r.h.s. of \eqref{hlsc} }  \ge \liminf_{n\to\infty}\Big\{  \frac{1}{2(a-b)}\left(\left(\mu-aD(x)\right)_+^2 -\left(\mu-aD(x^{\ssup{n}})\right)_+^2\right) \\-(\mu-aD(x))D(\eps^{\ssup{n}})+\frac{a-b}{2}\left(D(\eps^{\ssup{n}})\right)^2\Big\}.
\end{multline}
Recall that $ \liminf_{n\to\infty}D(\eps^{\ssup{n}})\ge 0 $, and thus we know that $ (\mu-aD(x))< 0 $ implies that eventually $ (\mu-aD(x)-aD(\eps^{\ssup{n}}))<0 $ and $ (\mu-aD(x)-aD(\eps^{\ssup{n}}) +bD(\eps^{\ssup{n}}))< 0 $.
Suppose that $ \mu/a<D(x) $. Then
\begin{equation}
    \mbox{r.h.s. of \eqref{hlsc2} }=\liminf_{n\to\infty}\Big\{ -\left(\mu-aD(x)\right)D(\eps^{\ssup{n}})+\frac{a-b}{2}D(\eps^{\ssup{n}})^2\Big\}\ge 0.
\end{equation}
Suppose $ \mu/a\ge D(x) $ and $ \mu-aD(x)-aD(\eps^{\ssup{n}}) \le 0 $. Then 
\begin{equation}
    \mbox{ r.h.s. of \eqref{hlsc2}  }\ge \frac{1}{2(a-b)}\liminf_{n\to\infty}\left\{\left(\mu-aD(x)-aD(\eps^{\ssup{n}})+bD(\eps^{\ssup{n}})\right)^2\right\}\ge 0,
\end{equation}
and likewise for $  \mu/a\ge D(x)   $ and $   \mu-aD(x)-aD(\eps^{\ssup{n}}) > 0  $,
\begin{multline}
\mbox{r.h.s. of \eqref{hlsc2} } \ge \frac{1}{2(a-b)} \liminf_{n\to\infty}\Big\{\left(\mu-aD(x)-aD(\eps^{\ssup{n}})+bD(\eps^{\ssup{n}})\right)^2 \\ 
- \left(\mu-aD(x)-aD(\eps^{\ssup{n}})\right)^2  \Big\}\ge 0.
\end{multline}

We have established \eqref{hlsc} and thus the lower semicontinuity of $ h $. Suppose $f$ is lower semicontinuous and $f\left(x\right)\leq H^\HYL_\mu\left(x\right)$ for all $x\in\ell_1\left(\R_+\right)$. Let $x^\ssup{n}\to x$ in $\ell_1\left(\R_+\right)$. Then $f\left(x\right)\leq\liminf_{n\to\infty}f\left(x^\ssup{n}\right) \leq \liminf_{n\to\infty}H^\HYL_\mu\left(x^\ssup{n}\right)$. Let $ x^{\ssup{n}}=x+ \eps^{\ssup{n}} $ with
\begin{equation}\label{seqeps}
\eps^{\ssup{n}}_k=\1\{k=n\}\frac{\left(\mu-aD(x)\right)_+}{n(a-b)},
\end{equation} so this inequality implies that $f\left(x\right)\leq h\left(x\right)$.
\qed

\end{proofsect}

\begin{prop}[\textbf{Upper  bound HYL model}]\label{prop-uboundHYL}
For all $ \mu\in\R, \alpha\le 0 $, and $ a > b\geq0 $,
\begin{equation}\label{uboundHYL}
\limsup_{N\to\infty}\frac{1}{\beta\abs{\L_N}}\log \E_{ \nu_{N,\alpha}}\Big[\ex^{-\beta\abs{\L_N} H^\HYL_\mu}\Big]\le -\inf_{x\in\ell_1(\R_+)}\big\{ I_\alpha(x)+H^\HYL_{\mu,\lsc}(x)\big\}.
\end{equation}
\end{prop}

\begin{proofsect}{Proof}
The statement follows easily with the upper bound estimate in Varadhan's Lemma given in \cite[Lemma~4.3.6]{DZ09} using  the inequality $ H^\HYL_{\mu}(x)\ge H^\HYL _{\mu,\lsc}(x) \ge H^\PMF_{\mu,(a-b),\lsc}(x)   $ for all $ x\in\ell_1(\R_+) $, the lower semicontinuity of $ H^\HYL_{\mu,\lsc} $,  and the fact that $ H^\HYL_{\mu,\lsc}(x)\ge -\frac{\left(\mu\right)^2}{2(a-b)} $. The latter estimate provides the tail-condition necessary to apply \cite[Lemma~4.3.6]{DZ09}.
\qed
\end{proofsect}

For the lower bound we are using the lower bound \eqref{lbound} for the PMF model and $ H^\HYL_\mu=H^\PMF_{\mu,a} + H_- $ with $ H_- $ being upper semicontinuous.

\begin{prop}[\textbf{Lower  bound HYLmodel}]\label{prop-lboundHYL}
For all $ \mu\in\R, \alpha\le 0 $, and $ a > b\geq0 $,
\begin{equation}\label{lboundHYL}
\liminf_{N\to\infty}\frac{1}{\beta\abs{\L_N}}\log \E_{\nu_{N,\alpha}}\Big[\ex^{-\beta\abs{\L_N} H^\HYL_\mu}\Big]\ge -\inf_{x\in\ell_1(\R_+)}\big\{ I_\alpha(x)+H^\HYL_{\mu,\lsc}(x)\big\}.
\end{equation}
\end{prop}

\begin{proofsect}{Proof}
Using
$$
 \E_{\nu_{N,\alpha}}\big[\ex^{-\beta\abs{\L_N} H^\HYL_\mu}\big]=\E_{\nu_{N,\alpha,\mu}^\PMF}\big[\ex^{-\beta\abs{\L_N} H_-}\big] \E_{\nu_{N,\alpha}}\big[\ex^{-\beta\abs{\L_N} H^\PMF_{\mu,a}}\big]
$$ 
in conjunction with the LDP in Theorem~\ref{THM-PMF} and in particular the lower bound  \eqref{lbound}  we arrive at
\begin{align*}
\liminf_{N\to\infty}&\frac{1}{\beta\abs{\L_N}}\log\Big( \E_{\nu_{N,\alpha,\mu}^\PMF}\big[\ex^{-\beta \abs{\L_N}H_-}\big] \E_{\nu_{N,\alpha}}\big[\ex^{-\beta\abs{\L_N} H^\PMF_{\mu,a}}\big]\Big)\\&\ge -\inf_{x\in\ell_1(\R)}\big\{I^\PMF_{\alpha,\mu}(x)+H_-(x)\big\} -\inf_{x\in\ell_1(\R)}\big\{I_\alpha(x) + H^\PMF_{\mu,a,\lsc}\big\} \\
& =-\inf_{x\in\ell_1(\R)}\big\{ I_\alpha(x)+H^\PMF_{\mu,a,\lsc}(x)+\beta H_-(x)\big\}= -\inf_{x\in\ell_1(\R)}\big\{ I_\alpha(x)+ H^\HYL_{\mu,\lsc}(x)\big\},
\end{align*}
where the last equality follows from Lemma~\ref{regularisation} below.

\qed
\end{proofsect}

\begin{lemma}\label{regularisation}
$$
\inf_{x\in\ell_1(\R_+)}\big\{ I_\alpha(x)+H^\PMF_{\mu,a,\lsc}(x)+H_-(x)\big\} =\inf_{x\in\ell_1(\R_+)}\big\{I_\alpha(x)+H^\HYL_{\mu,\lsc}(x)\big\}\,. 
$$
\end{lemma}
\begin{proofsect}{Proof}
The infimum of any function on an open set is equal to the infimum of its lower semicontinuous regularisation over the same set. Note that the $\ell_1$ topology restricted to $\ell_1\left(\R_+\right)$ has $\ell_1\left(\R_+\right)$ as open. We thus need to show that
\begin{equation}\label{reg}
\big(H^\PMF_{\mu,a,\lsc}+H_-\big)_{\lsc}(x)=H^\HYL_{\mu,\lsc}(x) \quad \mbox{ for all } x\in\ell_1(\R_+).
\end{equation}
Note that $H^\HYL_\mu(x)=H^\PMF_{\mu,a}(x) +H_-(x) \ge H^\PMF_{\mu,a,\lsc}(x)+H_-(x) \ge \big(H^\PMF_{\mu,a,\lsc}+H_-\big)_{\lsc}(x)$. To show \eqref{reg} we use the proof of Lemma~\ref{L:regHYL} and choose the sequence according to \eqref{seqeps}. We thus obtain $H^\HYL_{\mu,\lsc}(x)\ge \big(H^\PMF_{\mu,a,\lsc}+H_-\big)_{\lsc}(x)$.
\qed
\end{proofsect}

We finally combine Proposition~\ref{prop-uboundHYL} and Proposition~\ref{prop-lboundHYL}  to finish the proof for Theorem~\ref{Thm:HYL}.\qed

\subsection{Proof of Theorem~\ref{THM-densityLDP}}

To prove the first part of Theorem~\ref{THM-densityLDP} requires the limiting logarithmic moment generating function in Proposition~\ref{logmomentdensity} and  an application of the G\"artner-Ellis theorem.

\begin{proofsect}{Proof of Proposition~\ref{logmomentdensity}}  For $ t\in \R $ we get, using the independence of the Poisson point processes,
\begin{equation*}
{\tt E}\big[\ex^{\beta t\sum_{k+1}^\infty k\Ncal_k}\big]=\prod_{k=1}^\infty {\tt E}\big[\ex^{\beta t\sum_{k=1}^\infty k\Ncal_k}\big]=\prod_{k=1}^\infty \sum_{m=0}^\infty \ex^{\beta tkm}\P(\Ncal_k=m)=\prod_{k=1}^\infty \ex^{\abs{\L_N}q_k^{\ssup{\alpha}}(\ex^{\beta tk}-1)},
\end{equation*}
and thus
$$
\mathcal{L}(t)=\lim_{N\to\infty}\frac{1}{\beta \abs{\L_N}}\log {\tt E}\big[\ex^{\beta t\sum_{k=1}^\infty k\Ncal_k}\big]=\sum_{k=1}^\infty \frac{q_k^{\ssup{\alpha}}}{\beta}\big(\ex^{\beta tk}-1\big).
$$
\qed
\end{proofsect}

\begin{proofsect}{Proof of Theorem~\ref{THM-densityLDP}}
(a) This is a straightforward application of the G\"artner-Ellis theorem, see \cite{DZ09}. To ensure that $ 0 $ is in the domain of the logarithmic moment generating function we need to have $ \alpha<0 $. Then the rate function is giving as the Legendre-Fenchel transform
$$
\begin{aligned}
J_\alpha(x)&=\sup_{t\in\R}\big\{t x-\sum_{k=1}^\infty \frac{q_k^{\ssup{\alpha}}}{\beta}\big(\ex^{\beta tk}-1\big)\big\}&=\big(p\left(\beta,\alpha\right)+\sup_{t\in\R}\big\{(t+\alpha)x-p(\beta,\alpha+t)\big\} -\alpha x\big)\\
&=\big(p\left(\beta,\alpha\right)+f\left(\beta,x\right)-\alpha x\big)\,,
\end{aligned}
$$
where we used that for $ x\le \varrho_{\rm c} $,
$$
f(\beta,x)=\sup_{\alpha\in\R}\{\alpha x-p(\beta,\alpha)\}.
$$
Clearly, $ J_\alpha(x)=\infty $ when $ x<0 $ as the empirical density only takes positive values. Suppose that $ x>\varrho_{\rm c} $ for dimensions $ d\ge 3 $ ($ \varrho_{\rm c}=\infty $ for $ d=1,2 $). Then
\begin{equation}
\sup_{t\in\R}\left\{  \left(t+\alpha\right)\varrho_{\rm c} +\left(t+\alpha\right)\left(x-\varrho_{\rm c}\right)-p\left(\beta,\alpha+t\right)  \right\}\ge f\left(\beta,\varrho_{\rm c}\right)+\sup_{t\in\R}\left\{\left(t+\alpha\right)\left(x-\varrho_{\rm c}\right)\right\}=+\infty\,,
\end{equation}
and thus $ J_\alpha(x) \equiv +\infty $ for $ x\notin [0,\varrho_{\rm c}] $.

\medskip

\noindent (b) This is straightforward application of Varadhan's Lemma in \cite[Theorem~4.3.1]{DZ09} using that $ h(x):=x\mapsto -\mu x+\frac{a}{2}x^2 $ is continuous and the fact that the tail-condition is satisfied
$$
\lim_{M\to\infty}\limsup_{N\to\infty}\frac{1}{\beta \abs{\L_N}}\log {\tt E}\big[\ex^{-\beta \abs{\L_N}h(\brho_n)}\1\{h(\brho_n)\ge M\}\big]=-\infty.
$$
\qed
\end{proofsect}

\section{Variational analysis, pressure representation and condensation - Proofs}

This section collect all variational analysis proofs. We use frequently the following technical lemma for calculating derivatives of the rate functions.

\begin{lemma}
	\label{diffPressureTech}
	Let $I\subset\R$ be an open interval, $F:\ell_1 \times I \to \R$, and $\xi\in C^1\left(I;\ell_1\right)$. Also define
$$
	\mathcal{G}:I\to\R;\quad  s\mapsto F\left(\xi\left(s\right),s\right).
	$$
	Then if $F\left(x,s\right)$ is G{\^a}teaux differentiable in its first argument at $\xi\left(s\right)$ with $\left.\frac{\partial F}{\partial x_k}\right|_{\xi}=0$ $\forall k\in\N$, then
	\begin{equation*}
	\frac{\d \mathcal{G}}{\d s} = \left.\frac{\partial F}{\partial s}\right|_{\xi}.
	\end{equation*}
\end{lemma}

\begin{proofsect}{Proof}
An application of the chain rule gives $\frac{\d \mathcal{G}}{\d s} = \left.\frac{\partial F}{\partial s}\right|_{\xi} + \sum^\infty_{j=1}\frac{\d \xi_k}{\d s}\left.\frac{\partial F}{\partial x_k}\right|_{x=\xi}$. Since the partial derivatives of $F$ with respect to $x_k$ vanish at $\xi$, we only keep the first term.
\qed
\end{proofsect}
\smallskip

\subsection{Proofs for the CMF model}\label{sec-proofs-CMF} 
\begin{proofsect}{Proof of Proposition~\ref{zeroCMF}}
Recall that $I^{\CMF}_\alpha$ is lower semicontinuous and has compact level-sets. Also note that $I_\alpha$ is strictly convex where it is finite, and that $H^{\CMF}$ is convex. Therefore $I^\CMF_\alpha$ is strictly convex where it is finite (a non-empty set) and uniqueness of the minimiser follows.
	
	To calculate the minimiser, we search for stationary points. Since $I^\CMF$ is strictly convex where it is finite, if we find a stationary point then it is the global minimiser. By considering the coordinate derivatives, we know that the minimiser must satisfy all the following equations
	\begin{equation*}
	\frac{1}{\beta}\log \frac{x_k}{q^{\ssup{\alpha}}_k} + a \sum^\infty_{k=1}x_k  = 0, \quad k\in\N.
	\end{equation*}
	To make this more manageable, we introduce the dummy variable $\Gamma \in \R_+$ and corresponding equation $\Gamma = \sum^\infty_{k=1}x_k$. We shall  solve
	\begin{align}
	\label{eqn:CMFstat1}
	\log \frac{x_k}{q^{\ssup{\alpha}}_k} + a\beta \Gamma &= 0, \quad k\in\N,\\
	\label{eqn:CMFstat2}
	\Gamma - \sum^\infty_{k=1}x_k &= 0.
	\end{align}
	Given $\Gamma$, \eqref{eqn:CMFstat1} is uniquely solved by $ x_k = q^{\ssup{\alpha}}_k \exp\left( - a \beta \Gamma\right), k\in\N $, and therefore \eqref{eqn:CMFstat2} becomes
	\begin{equation*}
	\Gamma = \exp\left( -a\beta \Gamma\right)\bar{q}^{\ssup{\alpha}}.
	\end{equation*}
	This has the unique solution $ \Gamma = \frac{1}{a\beta}W_0\left(a\beta \bar{q}^{\ssup{\alpha}}\right) $,
	and so \eqref{eqn:CMFstat1} and \eqref{eqn:CMFstat2}  are uniquely jointly solved by $ x=\xi $ given by
	\begin{equation*}
	\xi_k = \frac{W_0\left(a\beta\bar{q}^{\ssup{\alpha}}\right)}{a\beta\bar{q}^{\ssup{\alpha}}}q^{\ssup{\alpha}}_k, \quad k\in\N.
	\end{equation*}
	\qed
\end{proofsect}


\medskip

\begin{proofsect}{Proof of Proposition~\ref{P:densityCMF}}
(a)  The continuity of $p^\CMF$ for $\alpha\leq 0$ follows from \eqref{p-cmf} and the continuity of $ a\beta \bar{q}^\ssup{\alpha}$ and $W_0$. Convexity follows from considering the derivatives of $W_0 $ with respect to $\alpha$ for $\alpha<0$, see Appendix~\ref{app-Lambert} for derivatives of the Lambert function.\\
(b) Smoothness follows from $W_0$ and the Bose functions being differentiable on the appropriate regions. The form of the first derivative can be found by either directly differentiating \eqref{p-cmf}, or by using Lemma~\ref{diffPressureTech} with the zero found in Proposition~\ref{zeroCMF}.\\
(c) We obtain \eqref{criticalCMF} from \eqref{critical} and the continuity of $ W_0 $ and $ \bar{q}^{\ssup{\alpha}} $. 
\qed
\end{proofsect}

\begin{proofsect}{Proof of Proposition~\ref{P:LFT-CMF}}
	This is proven in the same way as Proposition~\ref{IdealFreeEnergy}.
	\qed
\end{proofsect}

\bigskip

\bigskip

\subsection{Proofs for the PMF model}\label{proofs-PMF}

We collect our proofs for the PMF model.

\begin{proofsect}{Proof of Proposition~\ref{zeroPMF}} 
(a) To obtain the unique zero of the rate function we shall find the unique minimiser of the un-normalised rate function $ F(x):=I_\alpha(x)+ H^\PMF_{\mu,\lsc}(x) $. For the existence of a minimiser, recall that $F$ is lower semicontinuous and has compact level-sets. Also note that $ I_\alpha $ is strictly convex where it is finite, and $ H^\PMF_{\mu,\lsc} $ is also convex in the linear function $ D(x) $.  Therefore $F$ is strictly convex where it is finite ( a non-empty set) and uniqueness of the minimiser follows. To calculate the minimiser, we search for stationary points. Since $F$ is strictly convex where it is finite, if we find a stationary point then it is the global minimiser.  By considering again as in the proof of Proposition~\ref{zeroIdeal} the coordinate derivatives, we know that the minimiser must satisfy all the following equations
$$
\frac{1}{\beta}\log\frac{x_k}{q_k^{\ssup{\alpha}}}+ k\left(aD(x)-\mu\right)_+=0, \quad k\in\N\,.
$$
To make this more manageable, we introduce the dummy variable $ \delta\in\R_+ $ and corresponding equation $ \delta=D(x) $.
\begin{align}
\frac{1}{\beta}\log\frac{x_k}{q_k^{\ssup{\alpha}}}+k\left(a \delta-\mu\right)_+&=0, \quad k\in\N\,,\label{1stcond}\\
\delta-D(x)&=0.\label{seccond}
\end{align}
Given the value $ \delta $, \eqref{1stcond} is uniquely solved by $x_k=q_k^{\ssup{\alpha}}\exp\left(\beta k\left(\mu-a\delta\right)_-\right), k\in\N\,,
$ and therefore \eqref{seccond} becomes
\begin{equation}\label{solutionDelta}
\delta=\sum_{k=1}^\infty k q_k^{\ssup{\alpha}} \exp\left(\beta k\left(\mu-a\delta\right)_-\right)\,.
\end{equation}
Denote the right hand side in \eqref{solutionDelta} by $h(\delta) $,  and note that $ h(\delta)\to 0 $ as $ \delta\to\infty $. Furthermore,
$$
\lim_{\delta\to 0} h(\delta)=\begin{cases} \sum_{k\in\N} kq_k^{\ssup{\alpha+\mu}}=:\varrho(\alpha+\mu)\in(0,\infty)  & ,\mu<0,\\
\varrho(\alpha) & , \mu\ge 0\,,\end{cases}
$$
with
$$
\varrho(\alpha)=\begin{cases} \in (0,\infty) &, \alpha <0, d\ge 1\,,\\
 \infty & , \alpha\equiv 0 \wedge d=1,2\,,\\ 
 \varrho_{\rm c}(d) \in (0,\infty) &,\alpha\equiv 0 \wedge d\ge 3\,,\end{cases}
$$
see Figure~\ref{fig:PMFconsistency}. In all cases there exists a unique solution  which we denote $ \delta^* $. For $ \mu\le 0 $ the solution is $ \delta^*\in(0,\varrho(\alpha+\mu-a\delta^*)) $, and  for $ \mu>0 $ we have two additional cases, that is,
$$
\delta^*=\begin{cases}  \in(0,\varrho(\alpha+\mu-a\delta^*))  & , \mbox{ for } \mu \le 0\,,\\ \begin{cases} \in(\mu/a,\varrho(\alpha)) &\,, \mbox{ for } \mu< a\varrho(\alpha)\,,\\
\varrho(\alpha) &\,,  \mbox{ for } \mu\ge a \varrho(\alpha) \,, \end{cases} & , \mbox{ for } \mu>0\,. \end{cases}
$$

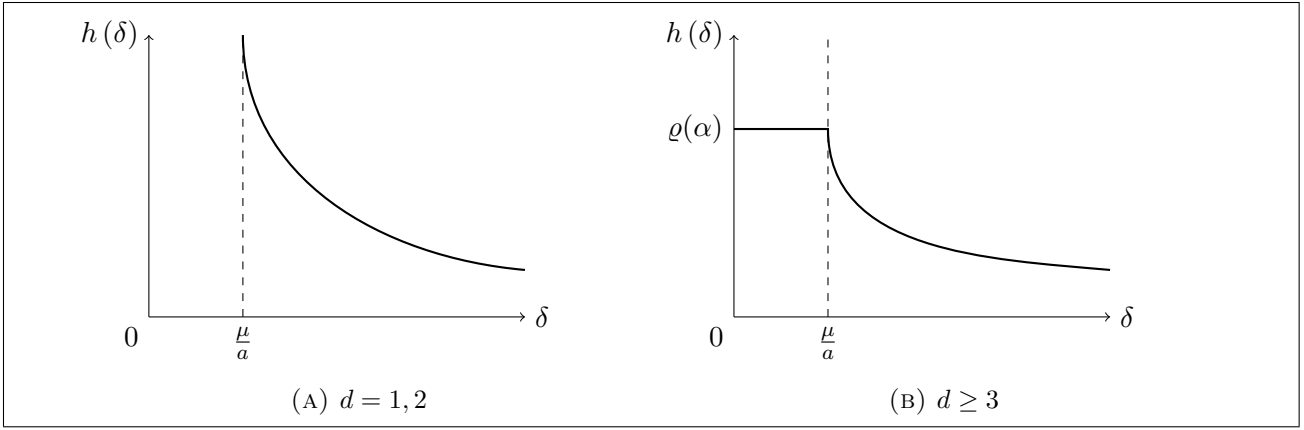
\begin{figure}
	\centering
	\begin{subfigure}[b]{0.45\textwidth}
		\begin{tikzpicture}[scale = 1.25]
		\draw[->] (0,0) node[below left]{$0$} -- (0,3) node[left]{$h\left(\delta\right)$};
		\draw[->] (0,0) -- (4,0) node[right]{$\delta$};
		\draw[dashed] (1,0) node[below]{$\frac{\mu}{a}$} -- (1,3);
		\draw[thick] (1,3) to [out=270,in=175] (4,0.5);
		\end{tikzpicture}
		\caption{$d=1,2$}
		\label{fig:PMFconsistency12}
	\end{subfigure}
	\begin{subfigure}[b]{0.45\textwidth}
		\begin{tikzpicture}[scale =1.25]
		\draw[->] (0,0) node[below left]{$0$} -- (0,3)node[left]{$h\left(\delta\right)$};
		\draw[->] (0,0) -- (4,0) node[right]{$\delta$};
		\draw[dashed] (1,0) node[below]{$\frac{\mu}{a}$} -- (1,3);
		\draw[thick] (0,2) node[left]{$\varrho(\alpha)$} -- (1,2) to [out=270,in=175] (4,0.5);
		\end{tikzpicture}
		\caption{$d\geq 3$}
		\label{fig:PMFconsistency>=3}
	\end{subfigure}
	\caption{Sketch of $h\left(\delta\right)$. This shows $\mu>0$, but the sketch translates with $\mu$.}
	\label{fig:PMFconsistency}
\end{figure}

\noindent (b) For the thermodynamic pressure we use the large deviation rate function and obtain, using the unique zero in (a),
$$
\begin{aligned}
p^\PMF(\beta,\alpha,\mu)&=-\inf_{y\in\ell_1(\R_+)}\big\{I_\alpha(y)+H^\PMF_{\mu,\lsc}(y)\big\}\\
&= p^\PMF(\beta,\alpha,\mu) =\begin{cases} p(\beta,\alpha)+\mu^2/2a &,\mbox{ for } \mu\ge a\varrho(\alpha), \delta^*=\varrho(\alpha)\,,\\ \frac{a}{2}\big(\delta^*\big)^2+p(\beta,\alpha+\mu-a\delta^*) &, \mbox{ for } \mu<a\varrho(\alpha), \delta^*\in(\mu/a,\varrho(\alpha))\,.\end{cases}
\end{aligned}
$$

\noindent (c)  This follows from either directly differentiating \eqref{PMFpressure2}, or by using Lemma~\ref{diffPressureTech} with the zero found in (a). Convexity also follows from (a) and (b), noting that $D\left(\xi^\PMF\right)$ is continuous and increasing in $\mu$. \qed
\end{proofsect}

\begin{proofsect}{Proof of Proposition~\ref{P:fPMF}}
For $ \mu \ge a \varrho(\alpha) $ we have
$$
\frac{\d p^{\PMF}}{\d\alpha}=\varrho(\alpha) \,\mbox{ and }\, \frac{\d p^\PMF}{\d\mu}=\frac{\mu}{a}=\delta^*=\varrho(\alpha)\,.
$$
Thus $ \mu\varrho(\alpha)-\mu^2/2a=\frac{a}{2}\varrho(\alpha)^2 $, and then the supremum over $ \alpha \le 0 $ gives the free energy for the ideal Bose gas and therefore the statement for this case. For the remaining case $ \mu<a\varrho(\alpha) $, we note that
$$
\frac{\d p^\PMF}{\d\alpha}=\frac{\d p^\PMF}{\d\mu}=\varrho(\alpha+\mu-a\delta^*)  
$$
with $ \delta^*\in (\mu/a,\varrho(\alpha)) $. Thus
$$
\begin{aligned}
\sup_{\alpha\le 0,\mu\in\R}&\big\{(\alpha+\mu)\varrho(\alpha+\mu-a\delta^*)-a/2(\delta^*)^2-p(\beta,\alpha,\mu-a\delta^*)\big\}=\\
\sup_{\alpha\le 0,\mu\in\R}&\big\{(\alpha+\mu-a\delta^*)\varrho(\alpha+\mu-a\delta^*)+a(\delta^*)^2-\frac{a}{2}(\delta^*)^2 -p(\beta,\alpha,\mu-a\delta^*) \big\}=f(\beta,\varrho)+a/2 \varrho^2\,.\end{aligned}
$$
\qed
\end{proofsect}


\bigskip

\subsection{Proofs variational analysis of the  HYL model}\label{sec-thermHYL}

\medskip

We present the variational analysis proofs for the HYL model.

\begin{proofsect}{Proof of Proposition~\ref{zeroHYL}}
\noindent (a) Suppose that the partial derivative $\frac{\partial I^\HYL_{\alpha,\mu}}{\partial x_k}$ is defined and non-zero at $x\in \mathrm{Int}\left\{\ell_1\left(\R_+\right)\right\}$. Then $I^\HYL_{\alpha,\mu}$ does not achieve its infimum at $x$. Since the boundary $\partial\ell_1\left(\R_+\right) = \left\{x\in \ell_1\left(\R_+\right) : \exists\, k \text{ s.t. }x_k=0\right\}$, and $\frac{\partial I^\HYL_{\alpha,\mu}}{\partial x_k} = -\infty$ here, the infimum is not achieved here.
	
For $x\in\mathrm{Int}\left\{\ell_1\left(\R_+\right)\right\}$ we have
	\begin{equation*}
	\frac{\partial I^\HYL_{\alpha,\mu}}{\partial x_k} (x)= \beta^{-1}\log \frac{x_k}{q_k^\ssup{\alpha}} - bk^2 x_k -  k \left(\mu - a D\left(x\right)\right)
	\begin{Bmatrix}
	1&,aD\left(x\right) \geq \mu\\
	-\frac{b}{a-b}&,aD\left(x\right) \leq \mu
	\end{Bmatrix}, \qquad k\in\N,
	\end{equation*}
	which is defined everywhere in $\mathrm{Int}\left\{\ell_1\left(\R_+\right)\right\}$. Hence a solution $\xi$ must solve $\frac{\partial I^\HYL_{\alpha,\mu}}{\partial x_k} (\xi)= 0$ for all $k\in\N $. To make this more manageable, we introduce the dummy variable $\delta\in\R_+$ and corresponding equation $\delta = D\left(x\right)$.  Our problem is then to solve
	\begin{align}
	\label{eqn:HYLstationary1}
	\log \frac{x_k}{q_k^\ssup{\alpha}} - b\beta k^2 x_k - \beta k \left(\mu - a \delta\right)
	\begin{Bmatrix}
	1&,a\delta \geq \mu\\
	-\frac{b}{a-b}&,a\delta \leq \mu
	\end{Bmatrix} &= 0, \qquad k\in\N,\\[1.5ex]
	\label{eqn:HYLstationary2}
	\delta - D\left(x\right) &= 0.
	\end{align}
	Unfortunately - unlike in the corresponding PMF case - even when we are given $\delta$ we are not guaranteed to have a solution for \eqref{eqn:HYLstationary1}, or that such a solution would be unique. If we fix $\delta$, then the $k^{th}$ equation of \eqref{eqn:HYLstationary1} either has no solution or is solved by 
	\begin{equation*}
	x_k = -\frac{1}{b\beta k^2}W_{\chi_k}\Big(-b\beta k^2 q_k^\ssup{\alpha} \exp\Big[\beta k \left(\mu - a \delta\right)
	\begin{Bmatrix}
	1&,a\delta \geq \mu\\
	-\frac{b}{a-b}&,a\delta \leq \mu
	\end{Bmatrix}
	\Big]\Big)
	\end{equation*}
	for all $\chi_k\in\left\{0,-1\right\}$, where $W_0$ and $W_{-1}$ are the two real branches of the Lambert W function. The `no solution' case corresponds precisely to $W_0$ and $W_{-1}$ not being defined for this input.  Substituting these $x_k$ back into \eqref{eqn:HYLstationary2} gives the condition \eqref{eqn:HYLconsistency} as required.
	
	\noindent (b) This follows immediately from the large deviation principle in Theorem~\ref{Thm:HYL}.
\qed
\end{proofsect}

\begin{proofsect}{Proof of Theorem~\ref{THM-uniqueness}}
We know that the global minimiser $ \xi$ of $F_\mu$ will equal
\begin{equation}\label{min1}
\xi_k=\xi_k(\delta,\chi)=-\frac{1}{b\beta k^2} W_{\chi_k}   \Big(-b\beta k^2 q_k^\ssup{\alpha} \exp\Big[\beta k \left(\mu - a \delta\right)\,,\quad k\in\N,
	\begin{Bmatrix}
	1&,a\delta \geq \mu\\
	-\frac{b}{a-b}&,a\delta \leq \mu
	\end{Bmatrix}
	\Big]\Big)
\end{equation}
for a choice of $ (\delta,\chi)\in\R_+\times\{0,-1\}^\N $ that satisfies
\begin{equation}\label{min2}
\delta=g^\chi(\delta):=\sum_{k\in\N} k\xi_k(\delta,\chi)\,.
\end{equation}
We first observe that the argument of the Lambert function in \eqref{min1} is decreasing in $k$ and in $ \delta\ge \mu/a $. Therefore, to ensure that the function $ g^\chi $ is finite and well-defined (i.e., series converges) the only admissible sequences $\chi $ as such that only finitely many indices are $-1 $ with all remaining ones pointing to the zero branch of the Lambert function. Furthermore, note that $ W_-(x) \ge \log (-x) $  for $ x\in [-\e^{-1},0) $. So if $ \chi_k=-1 $, then for $ \delta>\mu/a $ we have 
$$ 
g^\chi(\delta)>\frac{(a\delta-\mu)}{b} k-\frac{1}{b\beta}\log\big(b\beta k^2q_k^\ssup{\alpha}\big)\,. 
$$
Therefore there exists $ \mu_1\in\R $ such that for $ \mu<\mu_1 $ the equation $ \delta=g^\chi(\delta) $ has no solution for any of these $\chi $ with a finite number of negative branches. Furthermore, if $ \chi\equiv 0 $, then $g^\chi $ is strictly decreasing and continuous for sufficiently small values of $ \mu $. Note that $g^\chi $ in this case is decreasing for any $ \delta>\mu/a $ with maximum at $ \mu/a $. Hence for sufficiently small values of $ \mu <\mu_-=\mu_-(d,\beta,\alpha,a,b) $ there is a unique solution for $ \delta=g^\chi(\delta) $ \eqref{min2}. 

For the large positive $\mu $ case, we need to reduce the admissible number of values of $(\delta,\chi) $ we are investigating. We do this separately for the two regimes of $\delta $, $  \delta <\mu/a $ and $ \delta >\mu/a $. First suppose $ \delta<\mu/a $ and there exists $k\in\N $ such that we have a negative branch $ \chi_k=-1 $. Then we compute the second derivative
\begin{equation}
\frac{\partial^2 F_\mu}{\partial x_k^2}\Big|_{x=\xi(\delta,\chi)}=\frac{1}{\xi(\delta,\chi)}-b\beta k^2\big(1+\frac{a}{a-b}\big)<-\frac{ab}{(a-b)} \beta k^2\,,
\end{equation}
where the inequality holds because $ \chi_k=-1 $ implies that $ \xi_k>\frac{1}{\b \beta k^2} $. Since this $1$-dimensional Hessian is negative, this $ \xi(\delta,\chi) $ cannot be a local minimum and therefore not a global minimum. Now suppose that $ \delta>\mu/a $ and there exist $ k_1,k_2\in\N $ with $ k_1\not= k_2 $ such that $ \chi_{k_1}=\chi_{k_2}=-1 $. We then consider the two-dimensional Hessian in the $k_1 $ and the $k_2 $ coordinate directions. We find
$$
\det\Big(\frac{\partial^2 F_\mu}{\partial x_i\partial x_j}\Big|_{x=\xi(\delta,\chi)}\Big)_{i,j\in\{k_1,k_2\}}=\big(\frac{1}{\xi_{k_1}}+(a-b)\beta k_1^2\big)\big(\frac{1}{\xi_{k_2}}+(a-b)\beta k^2\big)-a^2\beta^2k_1^2 k_2^2< 0\,,
$$
where the inequality follows from $ \xi_k>\frac{1}{b\beta k^2} $ if $ \chi_k=-1 $ like before. Since this two-dimensional Hessian has a negative eigenvalue, this $ \xi(\delta,\chi) $ cannot be a local minimum and therefore not a global minimum. A similar analysis shows this for a finite number of distinct negative branches. In summary, we are left with solution of $ \delta=g^\chi(\delta) $ \eqref{min2} for which $ \chi\equiv 0 $ or there exists at most a single negative branch $k$ with $ \chi_k=-1 $ and $ \delta\ge \mu/a $. In order to compare  the remaining different candidates for minimiser, we calculate
$$
F_\mu(\xi)=-\sum_{k\in\N} \frac{\xi_k}{\beta}+\frac{b}{2}\sum_{k\in\N} k^2 \xi_k^2-\frac{a}{2}\delta^2-\frac{~(\mu-a\delta)_+^2}{2(a-b)}\,.
$$
Here we used \eqref{min2} and the property of the Lambert $W$ function that $W_{\chi_k}(x)\exp(W_{\chi_k}(x))=x $. We then approximate the remaining candidates for large values of $ \mu $. For $ \mu\gg 1$, $g^0(\delta)=\delta $ is solved only by a $ \delta\ll 1 $. We refer to this $ \chi\equiv 0 $ solution as $ \xi^\ssup{0} $. Then
$$
F_\mu(\xi^\ssup{0})=-\frac{\mu^2}{2(a-b)}+o(1)\,,
$$ as $ \mu\to\infty $. For the remaining competitors for minimiser, we let $ \xi^{\ssup{K}} $ refer to any solution corresponding to $ \chi_k=-\1\{k=K\}, k\in\N $. We are solely concerned with $ \delta\ge \mu/a $ for these cases. Hence $ \xi^\ssup{K}_K=\frac{1}{K}(\delta+O(1)) $, and $ \xi_k^\ssup{K}<\frac{1}{b\beta k^2} $ otherwise. This means
$$
F_\mu(\xi^\ssup{K})=-\frac{\delta}{K}-\frac{(a-b)}{2} \delta^2 +O(1)\,,
$$ for $ \mu\to\infty $.

If we investigate the difference $ F_\mu(\xi^\ssup{K})-F_\mu(\xi^\ssup{0}) $, we find that $ \xi^\ssup{K} $ is preferred if
\begin{equation}
\delta\gg \frac{1}{\beta K(a-b)}\big(-1+\sqrt{1+\beta^2\mu^2K^2}\big)<\frac{\mu}{(a-b)}-\frac{1}{\beta K(a-b)}\,.
\end{equation}
Conversely, $\xi^\ssup{0} $ is preferred if 
$$ 
\delta\ll \frac{\mu}{(a-b)}-\frac{1}{\beta K(a-b)} \,.
$$ In particular, this means that we need $ a\delta\gg \mu $ if $ \xi^\ssup{K} $ is to stand a chance. Now for $ a\delta\gg \mu $, and $ \chi_k=-\1\{k=K\} $ we have
$$
g^\chi(\delta)=\frac{a\delta-\mu}{a}+\frac{1}{b\beta K}\log(a\delta-\mu)+O(1)\,.
$$
For this calculation we use that
$$
W_{-1}(x)=\log(-x)-\log(-\log(-x))+o(1)\quad\mbox{ as } x\uparrow 0\,.
$$
Using this approximation with $ \delta-g^\chi(\delta) $ tells us that
$$
\delta=\frac{\mu}{(a-b)}-\frac{1}{\beta K(a-b)}\log(a\delta-\mu)+O(1)\,,
$$ and therefore $ \xi^\ssup{0} $ is always preferred for large values of $ \mu \ge \mu^+(d,\beta,\alpha,a,b)$. 
\qed
\end{proofsect}

\begin{proofsect}{Proof of Theorem~\ref{THM-nonuniqueness}}
We split the proof into three parts. In part I we discuss general solutions for selections of $ \chi $. In Step II we show that for the parameter regime given in the theorem one can rule out all $ \chi\not=0 $ solutions. In the final step III we prove the existence of the critical $ \mu^* $ with two simultaneous minimiser (zeros) for the $\chi=0 $ solution. 

\noindent\textbf{Step I:}  Note that for $d\geq 3$, the arguments for the Lambert W functions are strictly increasing in the summation index $k$, approaching $0$. This means that since the difference $W_0\left(x\right) - W_{-1}\left(x\right)\geq 0$ is strictly increasing in $x$ and equals $0$ if and only if $x=-\ex^{-1}$, we only need to consider finitely many $\chi$ for a given $\mu$ (all of which are eventually $0$). Now since any non-convexity in $g^\chi$ can only arrive via the finitely many $\chi_k=-1$ terms, solutions to $\delta = g^\chi\left(\delta\right)$ are locally finite in $\R$. To complement this, note that $\lim_{\delta\to+\infty}g^0\left(\delta\right)=0$ whilst for $\chi\ne 0$ we have $g^\chi\left(\delta\right) \gg \delta$. Hence we only need to consider a finite range of $\delta$, and therefore for a given $\mu$ there are only finitely many solutions for $\delta$.
	
	Because $g^\chi$ is continuous for each $\chi$, we can collect solutions uniquely and maximally into continuous paths $\xi^j\left(\mu\right)$ defined on closed (possibly infinite) intervals $I^j$ with non-empty interior. We allow families to overlap at endpoints of these intervals. Because we are only considering $\mu\leq \mu_{\rm p}$ and there are only finitely many solutions for each $\mu$, we will only have finitely many families being relevant to our discussion.
	
	For each of these families, we will denote 
	\begin{equation*}
	D^j\left(\mu\right) := D\left(\xi^j\left(\mu\right)\right), \qquad
	P^j\left(\mu\right) := -\left(I_\alpha + H^\HYL_{\mu,\lsc}\right)\left(\xi^j\left(\mu\right)\right),
	\end{equation*}
	defined on the interval $I^j$. From Proposition~\ref{zeroHYL}  we know that
	\begin{equation*}
	p^\HYL\left(\beta,\alpha,\mu\right) = -\inf_{y\in\ell_1(\R_+)}\left\{I_\alpha(y)+ H^\HYL_{\mu,\lsc}(y)\right\} = \max_{j}P^j\left(\mu\right).
	\end{equation*}
	Therefore for each $\mu$, there exists a $J$ such that $p^\HYL\left(\beta,\mu\right)=P^J\left(\mu\right)$.
	
	From the continuity of $g^\chi$ we know that all $D^j$ are continuous on their $I^j$. Then Lemma~\ref{diffPressureTech} tells us that each $P^j$ is differentiable on the interior $\mathrm{Int}(I^j)$, with derivative
	\begin{equation*}
	\frac{\d P^j}{\d \mu} = D^j + \frac{\left(\mu - a D^j\right)_+}{a-b}.
	\end{equation*}
	Continuity of this derivative follows from the continuity of $D^j$.

\noindent \textbf{Step II:} We show that for the following values of $ b $ and $ \mu $, i.e.,
$$
b <\min\big\{a,\frac{1}{\beta q_1^\ssup{\alpha}}\e^{-\mu_{\rm p}/a}\big\} \quad \mbox{ and } \; \mu<\mu_{\rm p}\,,
$$
we do not need to consider the $ \chi\not= 0 $ solutions of $ \delta=g^\chi(\delta) $. Note that if $ \chi\not=0 $, then $ g^\chi(\delta)>g^{\chi^\ssup{1}}(\delta) $, where $ \chi_k^\ssup{1}=-\1\{k=1\} $. By using $ W_{-1}(x)\le \log(-x) $, we find that if $ \chi\not= 0 $ we have
$$
g^\chi(\delta)>\frac{a}{b}\delta-\frac{\mu}{b}-\frac{1}{\beta}\log(\beta q_1^\ssup{\alpha}) +\frac{1}{\beta}\log (1/b)\,.
$$
Therefore we have no solutions if 
$$
b<\frac{1}{\beta q_1^\ssup{\alpha}}\exp\big(\beta\big(\frac{(a-b)}{b} \delta-\frac{\mu}{a}\big)\big)\,.
$$
The right hand side of the last inequality is bounded by $ \frac{1}{\beta q_1^\ssup{\alpha}}\e^{-\beta\mu/a} $, which is in turn bounded below by $ \frac{1}{\beta q_1^{\ssup{\alpha}}}\e^{-\beta \mu_{\rm p}/a} $ for the range of $ \mu $ we are considering. Since - by our assumptions  - $b$ is smaller than this bound, we can ignore $ \chi $ other than $ \chi=0 $.

	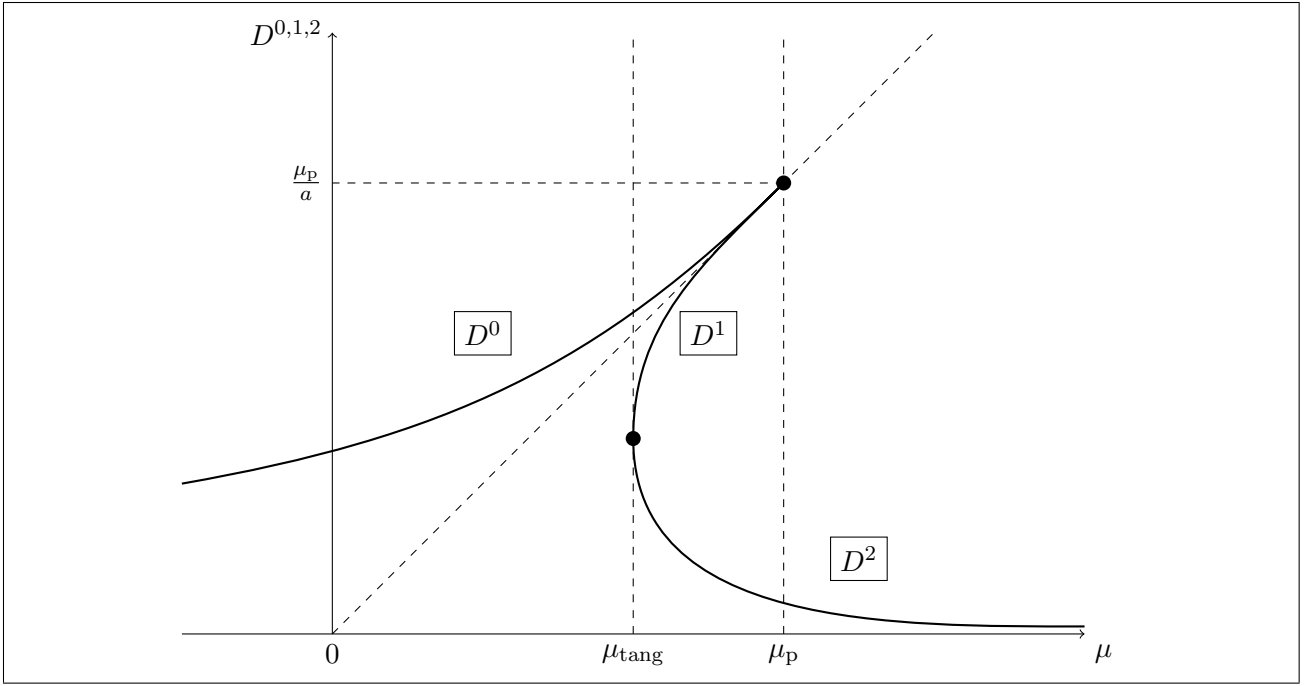
\begin{figure}
		\centering
		\begin{tikzpicture}[scale = 2]
		\draw[->] (-1,0) -- (5,0) node[below right]{$\mu$};
		\draw[->] (0,0) node[below]{$0$} -- (0,4) node[left]{$D^{0,1,2}$};
		\draw[dashed] (3,0) node[below]{$\mu_{\rm p}$}--(3,4);
		\draw[dashed] (2,0) node[below]{$\mu_{\rm tang}$}--(2,4);
		\draw[dashed] (0,0) -- (4,4);
		\draw[thick] (-1,1) to [out=10,in=225] (3,3) to [out=225,in=90] (2,1.3) to [out=270,in=180] (5,0.05);
		\fill (2,1.3) circle (0.05);
		\fill (3,3) circle (0.05);
		\draw[dashed] (3,3) -- (0,3) node[left]{$\frac{\mu_{\rm p}}{a}$};
		\node[draw] at (1,2) {$D^0$};
		\node[draw] at (2.5,2) {$D^1$};
		\node[draw] at (3.5,0.5) {$D^2$};
		\end{tikzpicture}
		\caption{Sketch of the total particle density of the three $\chi=0$ solutions for $d=3,4$.}
		\label{fig:HYL3sols}
	\end{figure}

\noindent\textbf{Step III:}   Let us now consider the $\chi=0$ solutions. Since $g^0$ is convex when restricted to $\delta\leq \frac{\mu}{a}$, $\frac{\d g^\chi}{\d \delta}\to+\infty$ as $\delta\uparrow\frac{\mu}{a}$, and $g^0$ is decreasing for $\delta \geq \frac{\mu}{a}$, there exists $\mu_{\rm tang}<\mu_{\rm p}$ such that this branch has multiple solutions if and only if $\mu\in\left[\mu_{\rm tang},\mu_{\rm p}\right]$. Let us label these $\xi^0$, $\xi^1$, and $\xi^2$ such that $D^0 > D^1\geq D^2$. Note that $I^0=\left(-\infty,\mu_{\rm p}\right]$, $I^1=\left[\mu_{\rm tang},\mu_{\rm p}\right]$, and $I^2=\left[\mu_{\rm tang},+\infty\right)$. For a visualisation of these solutions, see Figure~\ref{fig:HYL3sols}.

	Since $\xi^0\left(\mu_{\rm p}\right) = \xi^1\left(\mu_{\rm p}\right)$ and $\xi^1\left(\mu_{\rm tang}\right) = \xi^2\left(\mu_{\rm tang}\right)$, we have $P^0\left(\mu_{\rm p}\right) = P^1\left(\mu_{\rm p}\right)$ and $P^1\left(\mu_{\rm tang}\right) = P^2\left(\mu_{\rm tang}\right)$. Because $D^0\geq \frac{\mu}{a}$ and $D^{1,2}\leq \frac{\mu}{a}$, we have
	\begin{equation*}
		\frac{\d P^0}{\d \mu} = D^0, \qquad \frac{\d P^{1}}{\d \mu} = \frac{b}{a-b}\left(\frac{\mu}{b} - D^{1}\right) < \frac{b}{a-b}\left(\frac{\mu}{b} - D^{2}\right) = \frac{\d P^{2}}{\d \mu},
	\end{equation*}
	on $\left(\mu_{\rm tang},\mu_{\rm p}\right)$. Together these mean that $P^2\left(\mu_{\rm p}\right) > P^1\left(\mu_{\rm p}\right) = P^0\left(\mu_{\rm p}\right)$. Now extending our attention to all $\xi^j$ defined on some part of $\left(-\infty,\mu_{\rm p}\right]$, we define
	\begin{equation*}
		M:= \left\{\mu\leq\mu_{\rm p}: \exists\, j \text{ such that } P^j\left(\mu\right) > P^0\left(\mu\right)\right\}\,.
	\end{equation*}
	We have just shown that $M\ne\emptyset$, so $ \mu^*:=\inf M \in (-\infty,\mu_{\rm p}] $. Since $P^2$ and $P^0$ are continuous, $\mu^*<\mu_{\rm p}$. If $\mu^* \in M$, then $\max_{j}P^j\left(\mu\right)$ is discontinuous at $\mu=\mu^*$. In this case we are done. If otherwise $ \mu^* \notin M$, then since the $P^j$ are each continuous, $\exists\, J$ and $\epsilon>0$ such that $P^J\left(\mu\right) > P^0\left(\mu\right)$ for $\mu\in\left(\mu^*,\mu^* + \epsilon\right)$. Now we have to show that the derivatives of $P^0$ and $P^J$ necessarily have different limits as we take $\mu\to\mu^*$ from their respective sides. First note that
	\begin{equation*}
	\lim_{\mu\uparrow\mu^*}\frac{\d P^0}{\d\mu} = D^0\left(\mu^*\right)\,.
	\end{equation*}
Note that $g^\chi\left(\delta\right)\geq g^0\left(\delta\right)$ with equality only if $\chi=0$ or if we have both $\beta = \beta^*$ and $\delta = \frac{\mu}{a}$. Therefore $D^J\ne\frac{\mu}{a}$ and
	\begin{align*}
	D^J < \frac{\mu}{a} &\implies D^J\in\left[D^2,D^1\right]\,,\\
	D^J > \frac{\mu}{a} &\implies D^J > D^0\,.
	\end{align*}
This last inequality is strict because equality could only occur at $\mu = \mu_{\rm p}$, but $\mu^*<\mu_{\rm p}$. If $D^J>\frac{\mu}{a}$, then
	\begin{equation*}
	\lim_{\mu\downarrow\mu^*}\frac{\d P^J}{\d\mu} = D^J\left(\mu^*\right) > D^0\left(\mu^*\right) = \lim_{\mu\uparrow\mu^*}\frac{\d P^0}{\d\mu}
	\end{equation*}
	and we are done. From the symmetry of $g^0$ about $\delta = \frac{\mu}{a}$ and from $g^0$ being decreasing for $\delta \geq \frac{\mu}{a}$,  we have
	\begin{equation}
	\label{RelDensity}
	D^0 - \frac{\mu}{a} < \frac{b}{a-b}\left(\frac{\mu}{a}-D^1\right)
	\end{equation}
	for $\mu \in\left[\mu_{\rm tang},\mu_{\rm p}\right)$, see Figure~\ref{fig:HYL-ineq}. To see that, note that from the definition of $g^0(d)$, there is a symmetry to the function: the function to the right of $\mu/a$ is a scaled and reflected version of the function to the left. In Figure~\ref{fig:HYL-ineq}, this means that the distance from the leftmost intersection (the $D^1 $ intersection) to $\mu/a$ gives you the distance from $\mu/a$ to the imagined rightmost point. Since $g^0$ is decreasing to the right of $\mu/a$, this imagined point is to the right of the $D^2$ intersection, giving the inequality. This implies that if $D^J<\frac{\mu}{a}$,
	\begin{equation*}
	\frac{\d P^J}{\d \mu}  = \frac{b}{a-b}\left(\frac{\mu}{b} - D^J\right) \geq \frac{b}{a-b}\left(\frac{\mu}{b} - D^1\right) > D^0 = \frac{\d P^0}{\d \mu}.
	\end{equation*}
	Taking the limit to $ \mu^*$ gives our result.

	\begin{figure}
		\centering
		\begin{tikzpicture}
			\draw[dashed] (2,-0.5) -- (2,4);
			\draw[thick] (2,4) to [out=270,in=40] (0,0);
			\draw[thick] (2,4) to [out=270,in=170] (8,0);
			\draw[thick] (0,0.1) -- (8,4);
			\draw[dashed] (0,0.25) -- (8,0.25);
			\draw[<->] (0.3,-0.5) -- (2,-0.5);
			\draw[<->] (2,-0.5) -- (6.75,-0.5);
			\draw[<->] (2,4) -- (3.15,4);
			\node[draw] at (1,-1) {$\frac{\mu}{a}-D^1$};
			\node[draw] at (2.5,4.5) {$D^0-\frac{\mu}{a}$};
			\node[draw] at (3.5,-1) {$\frac{b}{a-b}\left(\frac{\mu}{a}-D^1\right)$};
			\fill (0.3,0.25) circle (0.1);
			\fill (6.75,0.25) circle (0.1);
			\fill (3.15,1.65) circle (0.1);
			\draw[dashed] (0.3,-0.5) -- (0.3,0.25);
			\draw[dashed] (6.75,-0.5) -- (6.75,0.25);
			\draw[dashed] (3.15,1.65) -- (3.15,4);
		\end{tikzpicture}
		\caption{\label{fig:HYL-ineq} Sketch of inequality \eqref{RelDensity}.}
	\end{figure}

\qed
\end{proofsect}

%
%

\bigskip

\begin{proofsect}{Proof of Theorem~\ref{thm:IDEALcondensate}}
		Let us begin with the $\alpha<0$ case. For $s\leq-\alpha$, fixed $N$, and fixed $K$, define
		\begin{equation*}
		g^\ssup{K}_N\left(s\right) := \frac{1}{\beta\left|\Lambda_N\right|}\log\E_{\nu_{N,\alpha}}\left[\exp\left(\left|\Lambda_N\right|s\beta\left(D-D_K\right)\right)\right].
		\end{equation*}
		Then
		\begin{equation*}
		\frac{\d g^\ssup{K}_N}{\d s}\Big|_{s=0} = \mathbb{E}_{\nu_{N,\alpha}}\big[D-D_K\big] \quad \text{and}\quad \Delta\left(\beta,\alpha\right) = \lim_{K\to\infty}\lim_{N\to\infty}\frac{\d g^\ssup{K}_N}{\d s}\Big|_{s=0}.
		\end{equation*}
		Since $g^\ssup{K}_N$ are all convex in $s$, we will use Griffith's Lemma to get the point-wise limit of the derivative from the derivative of the point-wise limit.
%
		To calculate the point-wise limit of $g^\ssup{K}_N$, we first rewrite $g^\ssup{K}_N$ as
		\begin{equation*}
		g^\ssup{K}_N\left(s\right) = \frac{1}{\beta\left|\Lambda_N\right|}\log\E_{\nu_{N,\alpha+s}}\left[\exp\left(-\left|\Lambda_N\right|s\beta D_K\right)\right] + \frac{1}{\beta\left|\Lambda_N\right|}\log \frac{Z_{N}\left(\beta,\alpha+s\right)}{Z_{N}\left(\beta,\alpha\right)}.
		\end{equation*}
		We use Varadhan's Lemma with  the tilt $\Phi = -s\beta D_K$. This $\Phi$ is continuous, but we need to pay attention to the boundedness conditions. We need to show that
		\begin{equation}\label{IDEALVaradhanConstraint}
		\lim_{M\to\infty}\limsup_{N\to\infty}\frac{1}{\beta\abs{\L_N}}\E_{\nu_{N,\alpha+s}}\big[\ex^{\left|\L_N\right|\Phi}\mathds{1}\left\{\Phi\geq M\right\}\big] = -\infty.
		\end{equation}
		For $0\leq s \leq -\alpha$, we have $\Phi\leq 0$ almost surely, so \eqref{IDEALVaradhanConstraint} holds trivially. For $s<0$ we have to work a little harder. Since $\Phi$ is continuous, the set $\left\{\Phi=m\right\}$ is closed (and measurable). Hence the upper bound of the LDP gives
		\begin{equation*}
		\limsup_{N\rightarrow\infty}\frac{1}{\beta\left|\Lambda_N\right|}\log \nu_{N,\alpha+s}\left(\Phi=m\right) \leq -\inf_{y\colon \Phi=m}\{I_{\alpha+s}(y)\}  \leq \bar{q}^\ssup{\alpha+s} -\bar{q}^\ssup{\alpha}+ \left(\alpha+s\right)\big(\frac{m}{\abs{s}\beta}\big).
		\end{equation*}
		This means that for sufficiently large $N$ there exists a $m$ and $N$ independent constant $C>\bar{q}^\ssup{\alpha+s}-\bar{q}^\ssup{\alpha}$, such that
		\begin{equation*}
		\ex^{m\left|\Lambda_N\right|}\nu_{N,\alpha+s}\left(\Phi=m\right) \leq \exp\big(\abs{\Lambda_N}\big[C + \frac{\alpha}{\abs{s}}m\big]\big)\,.
		\end{equation*}
		Since $\alpha<0$, we have sufficiently fast decay in $m$ to prove that \eqref{IDEALVaradhanConstraint} holds even for $s<0$, and Varadhan gives us
		\begin{equation*}
		\lim_{N\rightarrow\infty}g^\ssup{K}_N\left(s\right) = -\inf_{y\in\ell_1}\left\{ I_{\alpha+s}(y) + sD_K(y)\right\} + p\left(\beta,\alpha+s\right) -p\left(\beta,\alpha\right),\qquad \forall s\leq-\alpha.
		\end{equation*}
		
		In the style of Lemma \ref{zeroIdeal}, we can find that this infimum is achieved at $\xi=\xi(s)\in\ell_1\left(\R_+\right)$, where
		\begin{equation*}
		\xi_k=\begin{cases}
		q_k^{\ssup{\alpha}} &,k=1,\ldots,K\,,\\
		q_k^{\ssup{\alpha+s}} &, k>K\,.
		\end{cases}
		\end{equation*}
		Hence
		\begin{equation*}
			\lim_{N\rightarrow\infty}g^\ssup{K}_N\left(s\right) = \frac{1}{\beta}\sum^\infty_{k=K+1}\left(q_k^{\ssup{\alpha+s}}-q_k^{\ssup{\alpha}}\right),\qquad \frac{\d}{\d s}\Big(\lim_{N\to\infty}g^\ssup{K}_N\left(s\right)\Big)\Big|_{s=0} = \sum^\infty_{k=K+1} kq_k^{\ssup{\alpha}}\,.
		\end{equation*}
		Finally the sum vanishes as $K\to\infty$.
		
		For the case $\alpha=0$ with $d=1,2$, we take a more direct approach. It is clear from our construction that the point-wise limit $\lim_{N\to\infty}\E_{\nu_{N,\alpha}}\left[kx_k\right] = kq_k^{\ssup{\alpha}}$. Then for all $M\in\N$,
		\begin{align*}
		\liminf_{N\to\infty}\E_{\nu_{N,\alpha}}\Big[\sum^\infty_{k=K+1} kx_k\Big] \geq \lim_{N\to\infty}\E_{\nu_{N,\alpha}}\Big[\sum^M_{k=K+1}kx_k\Big] = \sum^M_{k=K+1}kq_k^{\ssup{\alpha}}\,.
		\end{align*}
		Since this lower bound diverges as $M\to\infty$ if $\alpha=0$ and $d=1,2$, we have our result for this case.
		
		For $\alpha=0$ and $d\geq3$ we use direct methods similar to the $\alpha>0$ case to get the required results.
For the CMF-model we follow the same steps as above.  For $s\leq-\alpha$, fixed $N$, and fixed $K$, define
	\begin{equation*}
	g^\ssup{K}_N\left(s\right) := \frac{1}{\beta\left|\Lambda_N\right|}\log\E_{\nu_{N,\alpha}^\CMF}\left[\exp\left(\left|\Lambda_N\right|s\beta\left(D-D_K\right)\right)\right].
	\end{equation*}
	Then by rearranging terms and applying Varadhan, we find that for $s\leq -\alpha$,
	\begin{equation*}
	\lim_{N\rightarrow\infty}g^\ssup{K}_N\left(s\right) = -\inf_{\ell_1(\R_+)}\big\{\frac{1}{\beta}I^\CMF_{\alpha+s} + sD_K\big\} + p^\CMF\left(\beta,\alpha+s\right) -p^\CMF\left(\beta,\alpha\right),\qquad \forall s\leq-\alpha\,.
	\end{equation*}
	
	In the style of Lemma \ref{zeroCMF}, we can find that this infimum is achieved at $\xi\left(s\right)\in\ell_1\left(\R_+\right)$, where
	\begin{equation*}
	\xi_k=\begin{cases}
	\frac{W_0\left(a\beta \bar{q}^\ssup{\alpha}\right)}{a\beta \bar{q}^\ssup{\alpha}}q_k^{\ssup{\alpha}} &,  k\leq K,\\
	\frac{W_0\left(a\beta \bar{q}^\ssup{\alpha+s}\right)}{a\beta \bar{q}^\ssup{\alpha+s}}q_k^{\ssup{\alpha+s}} &,  k>K.
	\end{cases}
	\end{equation*}
	Substituting this into the infimum, and then taking the derivative gives us
	\begin{equation*}
	\frac{\d}{\d s}\Big(\lim_{N\to\infty}g^\ssup{K}_N\left(s\right)\Big)\Big|_{s=0} = \frac{W_0\left(a\beta \bar{q}^\ssup{\alpha}\right)}{a\beta \bar{q}^\ssup{\alpha}}\sum^\infty_{k=K+1}kq^\ssup{\alpha}_k.
	\end{equation*}
	Finally the sum vanishes as $K\to\infty$, and Griffith's Lemma gives us the result.
	\qed

\end{proofsect}

\bigskip

\begin{proofsect}{Proof of Theorem~\ref{thm:PMFcondensate}}
	Our proof begins similarly to that of Theorem~\ref{thm:IDEALcondensate}. For $\mu,s\in\R$, fixed $N$, and fixed $K$, define
	\begin{equation*}
		g^\ssup{K}_N\left(s\right) = \frac{1}{\beta\left|\Lambda_N\right|}\log\E_{\nu_{N,\alpha,\mu}^\PMF}\left[\exp\left(\left|\Lambda_N\right|s\beta\left(D-D_K\right)\right)\right].
	\end{equation*}
	We once again rearrange terms to get
	\begin{equation*}
	g^\ssup{K}_N\left(s\right) = \frac{1}{\beta\left|\Lambda_N\right|}\log\E_{\nu_{N,\alpha,\mu+s}^\PMF}\left[\exp\left(-\left|\Lambda_N\right|s\beta D_K\right)\right] + \frac{1}{\beta\left|\Lambda_N\right|}\log \frac{Z^\PMF_{N}\left(\beta,\alpha,\mu+s\right)}{Z^\PMF_{N}\left(\beta,\alpha,\mu\right)}.
	\end{equation*}
	Then we want to use Varadhan's Lemma with our LDP for the PMF measure and the tilt $\phi = -s\beta D_K$. This $\phi$ is continuous, but we need to pay attention to the boundedness conditions. We will show
	\begin{equation}\label{VaradhanConstraint}
		\lim_{M\to\infty}\limsup_{N\to\infty}\frac{1}{\left|\L_N\right|}\E_{\nu_{N,\alpha,\mu+s}^\PMF}\Big[\ex^{\left|\L_N\right|\phi}\mathds{1}\left\{\phi\geq M\right\}\Big] = -\infty.
	\end{equation}
	For $s\geq0$, we have $\phi\leq 0$ almost surely, so \eqref{VaradhanConstraint} holds trivially. For $s<0$ we have to work a little harder. Our LDP for the PMF model gives us a bound on the probability of this set:
	\begin{align*}
		\limsup_{N\rightarrow\infty}\frac{1}{\left|\Lambda_N\right|}\log \nu^\PMF_{N,\alpha,\mu+s}\left(\phi=m\right) &\leq -\inf_{\phi=m}I^\PMF_{\alpha,\mu+s}\\
		\frac{m}{\left|s\right|\beta}\geq \frac{\mu+s}{a}\implies \qquad& \leq \inf_{\ell_1}\left\{I_\alpha +  H^\PMF_{\mu+s,\lsc}\right\} + \beta(\mu+s)\left(\frac{m}{\abs{s}\beta}\right) - \frac{a\beta}{2}\left(\frac{m}{\abs{s}\beta}\right)^2.
	\end{align*}
	This means that given $m \geq \left|s\right|\beta \frac{\mu+s}{a}$, then for sufficiently large $N$ there exists a $m$ and $N$ independent constant $C>\inf_{\ell_1}\left\{I_0 + \beta H^\PMF_{\mu+s,\lsc}\right\}$, such that
	\begin{equation*}
		\ex^{m\abs{\Lambda_N}}\nu^\PMF_{N,\alpha,\mu+s}\left(\phi=m\right) \leq \exp\Big(\abs{\Lambda_N}\left[C + \frac{\mu}{\abs{s}m} - \frac{a}{2\beta\abs{s}^2}m^2\right]\Big)\,.
	\end{equation*}
	The very fast decay with $m$ proves that \eqref{VaradhanConstraint} holds even for $s<0$, and Varadhan gives us
	\begin{equation*}
		\lim_{N\rightarrow\infty}g^\ssup{K}_N(s) = -\inf_{\ell_1}\left\{I^\PMF_{\alpha,\mu+s} + sD_K\right\} + p^\PMF\left(\beta,\alpha,\mu+s\right) -p^\PMF\left(\beta,\alpha,\mu\right),\qquad \forall s\in\R.
	\end{equation*}
	
	In the style of Lemma \ref{zeroPMF}, we can find that this infimum is achieved at $\xi\left(s\right)\in\ell_1\left(\R_+\right)$, where
	\begin{equation*}
		\xi_k=q_k^\ssup{\alpha} \exp\left(\beta k\left[\left(\mu + s-a\delta^*\right)_--s\mathds{1}\left\{k\leq K\right\}\right]\right),\quad k\in\N,
	\end{equation*}
	and $ \delta^*\left(s\right)$ is given implicitly as follows for the different cases. For $ \mu+s\le 0 $ and $ \mu\le 0 $, we have that
\begin{equation}
\delta^*=\sum_{k=1}^K kq_k^\ssup{\alpha} \exp(\beta k(\mu-a\delta^*))+\sum_{k>K} kq_k^\ssup{\alpha} \exp(\beta k(\mu+s-a\delta^*))=:\varrho^K(\alpha+\mu+s)\,,
\end{equation}
that is, $ \delta^*< \varrho^K(\alpha+\mu+s) $,  and for $ \mu+s>0 $ we have two cases as follows.
\begin{equation}
\delta^*=\begin{cases} \sum_{k=1}^K kq_k^\ssup{\alpha}\exp(\beta k(\mu-a\delta^*)) +\sum_{k>K} kq_k^\ssup{\alpha}\exp(\beta k(\mu+s-a\delta^*)) & , \mbox{ for } \delta^*\ge\frac{\mu+s}{a}\,,\\ \sum_{k=1}^K kq_k^\ssup{\alpha} \exp(-\beta k s)+\sum_{k>K} k q_k^\ssup{\alpha}=:\varrho^K(\alpha,s) & ,\mbox{ for } \delta^*\le \frac{\mu+s}{a}\,.\end{cases} 
\end{equation}
Thus
\begin{equation}
\delta^*=\begin{cases} \in (0,\varrho^K(\alpha+\mu+s) & , \mbox{ for } \mu+s\le 0,\mu\le 0 \,,\\
\begin{cases} \in (\frac{\mu+s}{a},\varrho^K(\alpha,s)) &,\mbox{ for } \mu+s< a\varrho^K(\alpha,s)\,,\\ \varrho^K(\alpha,s) &,\mbox{ for } \mu+s \ge a \varrho^K(\alpha,s)\,.\end{cases} &,\mbox{ for } \mu+s>0\,. \end{cases} 
\end{equation}
	
%
	If we denote
	\begin{equation}
		\delta^*_K=\begin{cases}
			\sum^K_{k=1}kq_k \exp\left(\beta k\left(\mu-a\delta^*\right)\right) &, \mbox{ for } \delta^* \geq \frac{\mu+s}{a},\\
			\sum^K_{k=1}kq_k \exp\left(-s\beta k\right)  &, \mbox{ for } \delta^* \leq \frac{\mu+s}{a},
		\end{cases}
	\end{equation}
	then Lemma \ref{diffPressureTech} tells us that
	\begin{align*}
		\frac{\d}{\d s}\Big(\lim_{N\rightarrow\infty}g^{\ssup{K}}_N\left(s\right)\Big) &= \left(\delta^* - \delta^*_K\right)\left(s\right) + \big(\frac{\mu+s}{a}-\varrho^K(\alpha,s)\big)_+\\
		\frac{\d}{\d s}\Big(\lim_{N\rightarrow\infty}g^\ssup{K}_N\left(s\right)\Big)\Big|_{s=0} &= \begin{cases}
			\left(\delta^* - \delta^*_K\right)\left(0\right) &, \mbox{ for } \mu\leq a\varrho^K(\alpha,0)\,,\\
			\frac{\mu}{a} - \delta^*_K\left(0\right) &, \mbox{ for } \mu \geq a\varrho^K(\alpha,0)\,,
		\end{cases}\\
		\lim_{K\to\infty}\frac{\d}{\d s}\Big(\lim_{N\rightarrow\infty}g^\ssup{K}_N\left(s\right)\Big)\Big|_{s=0} &= \left(\frac{\mu}{a}-\varrho(\alpha)\right)_+.
	\end{align*}
	\qed
\end{proofsect}

\bigskip

\begin{proofsect}{Proof of Theorem~\ref{THM:HYL-Cond}}
	The proof of \eqref{HYLcondenstate1} follows very similarly to the corresponding stage of the proof of Theorem \ref{thm:PMFcondensate}. Note that the HYL version of \eqref{VaradhanConstraint} follows because $H^\HYL_{\mu,\lsc}\geq H^\PMF_{\mu,a-b,\lsc}$ almost surely. The proof of \eqref{HYLcondenstate2} uses Lemma \ref{diffPressureTech}.
	\qed
\end{proofsect}

\section*{Appendices}
\begin{appendices}

\section{Bose function}\label{app-Bose}
The Bose functions are poly-logarithmic functions defined by 
\begin{equation}\label{defBosefunctionapa}
g(n,\alpha):=\Li_n(\ex^{-\alpha})=\frac{1}{\Gamma(n)}\int_0^\infty\frac{t^{n-1}}{{\rm e}^{t+\alpha}-1}\,\d t=\sum_{k=1}^\infty k^{-n}{\rm e}^{-\alpha k}\quad\mbox{ for all } n \mbox{ and } \alpha>0,
\end{equation}
and also for $ \alpha=0 $ and $ n>1 $. In the latter case,
\begin{equation}\label{zeta}
g(n,0)=\sum_{k=1}^\infty k^{-n}=\zeta(n),
\end{equation}
which is the zeta function of Riemann. The behaviour of the Bose functions about $ \alpha=0 $ is given by
\begin{equation}\label{boseexp}
g(n,\alpha)=\begin{cases}
\Gamma(1-n)\alpha^{n-1}+\sum_{k=0}^\infty\zeta(n-k)\frac{(-\alpha)^k}{k!} &, n\not= 1,2,3,\ldots, \\[1.5ex]
\frac{(-\alpha)^{n-1}}{(n-1)!}\left[-\log\alpha +\sum_{m=1}^{n-1}\frac{1}{m}\right]+\sum_{\heap{k=0}{k\not= n-1}}\zeta(n-k)\frac{(-\alpha)^k}{k!} &, n\in\N.                   
\end{cases}
\end{equation}
%

\section{The ideal Bose gas}\label{appIBG}


We review the  large deviation principle for the ideal Bose, that is, large deviation principle for the empirical cycles counts. However, we  present this well-known result by using rigorous large deviation proofs based on Baldi's  Theorem and exponential tightness. Denote $ \nu_{N,\alpha} $ the distribution of the empirical cycle count  with chemical potential  $ \alpha\le 0 $. 

\begin{theorem}[\textbf{Ideal Bose gas}]
\label{THM-Ideal}
	For $d\in\N,  \beta>0  $, and $\alpha\leq0$,  the sequence $ (\nu_{N,\alpha})_{N\in\N}$ satisfies a large deviation principle (LDP) on $\ell_1(\R_+)$ with rate $\beta \abs{\Lambda_N} $ and rate function
	\begin{equation*}
		I_\alpha\left(x\right) = \sum^\infty_{k=1}\frac{x_k}{\beta}\Big(\log\frac{x_k}{q_k^{\ssup{\alpha}}} - 1\Big) + \bar{q}^{\ssup{\alpha}}/\beta\,.
	\end{equation*}
\end{theorem}

Since each entry in the $\ell_1$-valued empirical cycle count is independent and related to a Poisson random variable under the ideal Bose gas model, Theorem~\ref{THM-Ideal} is proven by applying Baldi's Lemma.

\begin{remark}
	The condition that $\alpha\leq 0$ arises from the $\bar{q}^{\ssup{\alpha}}$ term. Clearly, $\bar{q}^{\ssup{\alpha}}$ is finite if and only if $\alpha\leq0$. \hfill $ \diamond $
\end{remark}

\bigskip

\begin{prop}[\textbf{Pressure}]\label{zeroIdeal}
\begin{enumerate}[(a)]
\item The rate function for the ideal Bose gas model, $ I_\alpha $, has a unique zero $\xi\in\ell_1\left(\R_+\right)$ given by
\begin{equation}
\xi_k=q_k^{\ssup{\alpha}},\quad k\in\N.
\end{equation}
\item  

Let $ \beta>0 $ and $ \alpha\le 0 $, then the thermodynamic limit of the pressure
\begin{equation}
\begin{aligned}
	p(\beta,\alpha) &=\lim_{N\to\infty}\frac{1}{\beta\abs{\L_N}}\log \ex^{\abs{\L_N}\bar{q}^{\ssup{\rm bc},\alpha}}=   \inf_{y\in\ell_1(\R_+)}\Big\{ \sum_{k=1}^\infty \frac{y_k}{\beta}\Big(\log\frac{y_k}{q_k^{\ssup{\alpha}}}-1\Big)\Big\}=\frac{\bar{q}^{\ssup{\alpha}}}{\beta}\\&\label{p}
\end{aligned}
\end{equation}
exists.
\end{enumerate}

\end{prop}

\bigskip

We denote 
$$ p_{\L_N} =\frac{1}{\beta\abs{\L_N}} Z_{\L_N}(\beta,\alpha) 
$$ the average finite-volume pressure.

\medskip

\begin{prop}[\textbf{Thermodynamics}]\label{P:ideal1}
\begin{enumerate}[(a)]
\item For $ \beta> 0$, $\alpha>0 $, we define $ p(\beta,\alpha)=+\infty $. Then the pressure $ p(\beta,\cdot) $ is a closed convex function on $ \R $.

\item For $\beta>0$, $\alpha<0$, the ideal gas pressure $p(\beta,\alpha)$ is smooth with respect to $\alpha$. In particular, 
\begin{equation*}
\frac{\d p}{\d \alpha} = D\left(q^{\ssup{\alpha}}\right).
\end{equation*}

\item For $ \beta > 0, \alpha< 0 $, and any $ N\in\N $,
\begin{equation}
\frac{1}{\abs{\L_N}}{\tt E}\big[N_{\L_N}^{\ssup{\ell}}\big]=\frac{\d}{\d\alpha} p_{\L_N}(\beta,\alpha).
\end{equation}
The function $ \alpha\mapsto \frac{\d}{\d\alpha} p_{\L_N}(\beta,\alpha) $ is increasing on $ (-\infty,0) $. It follows that we can give $ \frac{1}{\abs{\L_N}}{\tt E}[N_{\L_N}^{\ssup{\ell}}]$ any pre-assigned value $ \varrho\in(0,\infty) $ by choosing $ \alpha = \alpha_N(\varrho)\in(-\infty,0) $.

\item In the thermodynamic limit $ N\to\infty $,
\begin{equation}\label{critical}
\varrho_{\rm c}:=\lim_{\alpha\uparrow 0}\Big(\frac{\d}{\d\alpha} p_{\L_N}(\beta,\alpha)\Big)=\begin{cases}  +\infty &, d=1,2,\\   \frac{1}{\left(4\pi\beta\right)^\frac{d}{2}}\zeta\left(\frac{d}{2}\right) &, d\ge 3. \end{cases}
\end{equation}
Let $ \alpha_N(\varrho) $ denote the unique root of 
\begin{equation}
\frac{\d}{\d\alpha} \big(p_{\L_N}(\beta,\alpha)\big)=\varrho
\end{equation}
then $ \alpha(\varrho)=\lim_{N\to\infty} \alpha_N(\varrho) $ exists and is equal to the unique root of
\begin{equation}
\frac{\d}{\d\alpha} \big(p(\beta,\alpha)\big)=\varrho \quad\mbox{ if } \varrho<\varrho_{\rm c},
\end{equation}
and it is equal to zero otherwise.
\end{enumerate}

\end{prop}

\bigskip

\begin{prop}[\textbf{Free energy}]
	\label{IdealFreeEnergy}
	For $\varrho>0$, we define the ideal Bose gas free energy as the Legendre-Fenchel transform of the pressure,
	\begin{align}
		f\left(\beta,\varrho\right):=\sup_{s\in\R}\left\{s\varrho-p\left(\beta,s\right)\right\}
		= \begin{cases}
		\frac{-1}{\beta\left(4\pi\beta\right)^{\frac{d}{2}}}g(1+\frac{d}{2},-\beta\gamma) + \varrho\gamma &,  \varrho\leq\varrho_{\rm c}\; ,\\
		\frac{-1}{\beta\left(4\pi\beta\right)^\frac{d}{2}}\zeta\left(1+\frac{d}{2}\right) &,\varrho\geq\varrho_{\rm c}\;, 
		\end{cases}
	\end{align}
	where $\gamma\le 0 $ is a solution to
	\begin{equation*}
	\frac{1}{\left(4\pi\beta\right)^\frac{d}{2}}g\big(\frac{d}{2},-\beta\gamma\big) = \varrho,
	\end{equation*}
	which exists and is unique for $\varrho\leq\varrho_{\rm c}$.
\end{prop}

It is easy to see that $ \varrho\mapsto f(\beta,\varrho) $ is a decreasing convex function; it is given by
$$
f(\beta,\varrho)=\alpha(\varrho)\varrho-p(\beta,\alpha(\varrho)) \quad\mbox{ for } \varrho<\varrho_{\rm c}.
$$
The linear segment in the graph of $ f$ where $ f $ is constant and equal to $ -p(\beta,0) $ for $ \varrho\ge \varrho_{\rm c} $, signals a first-order phase-transition at $ \alpha=0 $. This phase-transition is called \emph{Bose-Einstein condensation (BEC)}.

\bigskip

\noindent \textbf{Proofs for the reference measure (ideal Bose gas)}:


\begin{proofsect}{Proof of Theorem~\ref{THM-Ideal}}
We recall the following theorem for the convenience of the reader.
\begin{lemma}[\textbf{Baldi's Theorem}]\label{Baldi}
	Suppose $(\nu_N)_{N\in\N} $ is an exponentially tight sequence of measures on $\ell_1(\R)$. Let $\L\colon\ell_\infty(\R)\to [0,\infty] $ be the limiting cumulant generating function, and suppose that it exists and is finite for every $t\in\ell_\infty(\R)$. If $\L$ is G{\^a}teaux differentiable, and lower semicontinuous on $\ell_\infty(\R)$, then $(\nu_N)_{N\in\N} $ satisfies an LDP with rate function
	\begin{equation}
		\L^*(x) = \sup_{t\in\ell_\infty(\R)}\big\{\langle t,x\rangle - \L(t)\big\},\qquad x\in\ell_1(\R).
	\end{equation}
\end{lemma}

We shall now set about establishing that the hypotheses of Baldi's Theorem are satisfied. We adapt a beautiful proof in a recent study of Bosonic loop measures on graphs given in \cite{Dan15}.

\begin{lemma}[\textbf{Exponential tightness}]\label{exptight}
For every  $ \alpha\le 0 $, 	$\big(\nu_{N,\mu}\big)_{N\in\N} $ is an exponentially tight sequence of measures.
\end{lemma}

\begin{proofsect}{Proof}
	Suppose there exists an $x=x\left(\gamma\right)\in\ell_1\left(\R\right)$ such that for all $k\geq1$,
	\begin{equation*}
		\limsup_{N\rightarrow\infty}\frac{1}{\left|\Lambda_N\right|}\log\nu_{N,\alpha}\left(\blambda^{\ssup{k}}_{N}\geq x_k\right) <-2^{-k}\gamma,
	\end{equation*}
	where $\blambda_N=\big(\blambda^{\ssup{k}}_{N}\big)_{k\in\N} $ is an $\ell_1\left(\R\right)$-valued random variable with law $\nu_{N,\alpha}$. Also, define the set
	\begin{equation*}
		K = \left\{y\in\ell_1\left(\R\right)\colon \left|y_k\right| \leq \left|x_k\right| \forall k \geq 1\right\},\qquad x\in\ell_1(\R).
	\end{equation*}
To show compactness of $K$ it is  easy to see that $K$ is bounded and closed. Boundedness follows from $ \norm{y}_{\ell_1(\R)}\le \norm{x}_{\ell_1(\R)} $ for all $ y\in K $. Suppose that $K$ is not closed, that is, there exists a sequence $ y^{\ssup{n}}\in K $ with limit $ y^{\ssup{n}}\to y\notin K $ as $ n\to\infty $. Suppose that $ \abs{y_k}>\abs{x_k} $ for $ y\notin K $ and some $k\in\N $. Choose $ \eps=\frac{1}{2}(\abs{y_k}-\abs{x_k}) $, then, for $ n $ sufficiently large, 
$$
\abs{y_k^{\ssup{n}}-y_k}\le \sum_{j\in\N}\abs{y^{\ssup{n}}_j-y_j}< \frac{1}{2}(\abs{y_k}-\abs{x_k}),
$$
which implies that 
$$
\abs{y^{\ssup{n}}_k}>\abs{y_k}-\eps=\frac{1}{2}(\abs{y_k}+\abs{x_k})>\abs{x_k},
$$ contradicting $ y^{\ssup{n}}\in K $. Hence, $ K$ is closed. It remains to show that $K$ is totally bounded.  From that, we shall find a finite cover of $\eps$-open balls for $K$. Pick $ \eps> 0 $, and choose $ N\in\N$ such that $ \sum_{k>N} \abs{x_k}<\eps/2 $, and define the so-called cut-off sequences $ \widetilde{K}=\{y\in K\colon y_k=0, k> N\} $. Clearly, $ \widetilde{K} $ is isomorphic to the totally bounded set 
$$
[-\abs{x_1},\abs{x_1}]\times\cdots[-\abs{x_N},\abs{x_N}] \subset \R^N,
$$
and thus it is itself totally bounded. There exist $ w^{\ssup{1}},\ldots,w^{\ssup{M}}\in \widetilde{K} $ such that 
$$ 
\widetilde{K}\subset \bigcup_{i=1}^M B(w^{\ssup{i}},\frac{\eps}{2}).
$$
For any $y\in $ denote $ \widetilde{y}\in\widetilde{K} $ the sequences which agrees with $y$ on the first $ N$ terms, and choose $ w^{\ssup{i}} $ such that $ \widetilde{y}\in B(w^{\ssup{i}},\frac{\eps}{2}) $. Then,
$$
\norm{y -w^{\ssup{i}}}_{\ell_1(\R)}=\sum_{k=1}^N\abs{\widetilde{y}_k-w^{\ssup{i}}_k}+\sum_{k>N}\abs{y_k}<\frac{\eps}{2}+\frac{\eps}{2}.
$$
Thus, $ K\subset\bigcup_{i=1}^M B(w^{\ssup{i}},\frac{\eps}{2}) $, and we conclude with the compactness of $K$.

	Now, since the $\blambda^{\ssup{k}}_N$ are independent, we have
	\begin{equation*}
		\limsup_{N\rightarrow\infty}\frac{1}{\left|\Lambda_N\right|}\log\nu_{N,\alpha}\left(K^{\rm c}\right) = \limsup_{N\rightarrow\infty}\frac{1}{\left|\Lambda_N\right|}\sum_{k\in\N}\log\nu_{N,\alpha}\left(\blambda^{\ssup{k}}_N>x_k\right)<-\gamma,
	\end{equation*}
	and conclude with the statement in the lemma.
	All that remains now is to find such a sequence $x$. We consider each $x_k$ in turn. For all constants $c\geq0$, and $\tau>0$, we have the Chernoff bound
	\begin{align*}
		\nu_{N,\alpha}\left(\blambda^{\ssup{k}}_N> c\right) &= \nu_{N,\alpha}\left(\ex^{\frac{\tau}{\abs{\L_N}} \Ncal_k} > \ex^{\tau c}\right)\\
		&\leq \ex^{-\tau c}{\tt E}\left[\ex^{\frac{\tau}{\abs{\L_N}} \Ncal_k}\right]\\
		&= \ex^{-\tau c}\exp\left(\left|\Lambda_N\right|q^{\ssup{\alpha}}_k\left(\ex^{\frac{\tau}{\left|\Lambda_N\right|}}-1\right)\right).
	\end{align*}
	
	Differentiating this bound with respect to $\tau$ gives us that the minimum occurs at $\tau^* = \abs{\Lambda_N}\log\frac{c}{q^{\ssup{\alpha}}_k}$. If $c>0$, then $\tau^*>0$ for sufficiently large $N$. This means that we can optimise this form of bound as
	\begin{equation*}
		\nu_{N,\alpha}(\blambda^{\ssup{k}}_N>c) \leq \Big(\frac{c}{q^{\ssup{\alpha}}_k}\Big)^{-\left|\Lambda_N\right|c}\exp\left(\left|\Lambda_N\right|\left(c-q^{\ssup{\alpha}}_k\right)\right).
	\end{equation*}
	Taking $N\rightarrow\infty$ then gives us
	\begin{equation*}
			\limsup_{N\rightarrow\infty}\frac{1}{\abs{\Lambda_N}}\log\nu_{N,\alpha}(\blambda^{\ssup{k}}_N>c) \leq c - q^{\ssup{\alpha}}_k -c\log\frac{c}{q^{\ssup{\alpha}}_k}.
	\end{equation*}
	
	Now note that on $c>0$, the maps
	\begin{equation*}
		c\mapsto c - q^{\ssup{\alpha}}_k- c\log\frac{c}{q^{\ssup{\alpha}}_k} + 2^{-k}\gamma, \qquad k\in\N,
	\end{equation*}
	are differentiable, strictly decreasing, and have at most a unique zero $c^*_k$. If there does not exist such a zero, then the map is strictly negative, and it will suffice in what follows to set $c^*_k=0$. Since our maps are strictly negative for $c>c^*_k$, we only need to find a sequence $x$ such that $x_k>c^*_k$ for all $k$. Now we only need to find such an $x\in\ell_1\left(\R\right)$.
	
	Consider $x_k = c^*_k + 2^{-k}$. Therefore $x\in\ell_1\left(\R\right)$ if and only if $c^*\in\ell_1\left(\R\right)$. If we defined $c^*_k$ as a zero, then $c^*_k$ solves
	\begin{equation*}
		c^*_k\Big(1-\log\Big(\frac{c^*_k}{q^{\ssup{\alpha}}_k}\Big)\Big) = q^{\ssup{\alpha}}_k + 2^{-k}\gamma\,.
	\end{equation*}
	Otherwise, $c^*_k=0$ and
	\begin{equation*}
		c^*_k\Big(1-\log\Big(\frac{c^*_k}{q^{\ssup{\alpha}}_k}\Big)\Big) = 0.
	\end{equation*}	
	So noting that $q^{\ssup{\alpha}}_k>0 $ and $ \gamma>0$ and that the sum of $q^{\ssup{\alpha}}_k$ converges give us
	\begin{equation*}
		\sum_{k\in\N}c^*_k\Big(1-\log\Big(\frac{c^*_k}{q^{\ssup{\alpha}}_k}\Big)\Big) \leq \gamma + \sum_{k\in\N}q^{\ssup{\alpha}}_k < \infty.
	\end{equation*}
	
	Suppose, for contradiction, that $\sum_{k\in\N}c^*_k=\infty$. Then, in order for the left hand side of the above inequality to converge, we require $1-\log\Big(\frac{c^*_k}{q^{\ssup{\alpha}}_k}\Big)\rightarrow0$ as $k\rightarrow\infty$. Consequently, there exists a $K\geq1$ such that for $k\geq K$, $\frac{c^*_k}{q^\star_k}<3$, and hence
	\begin{equation*}
		\sum_{k\geq K}c^*_k \leq 3\sum_{k\geq K}q^{\ssup{\alpha}}_k\leq 3\sum_{k\in\N}q^{\ssup{\alpha}}_k<\infty\,.
	\end{equation*}
	We have a contradiction, and $\sum_{k\in\N}c^*_k<\infty$ as required.\qed
\end{proofsect}

%

\begin{lemma}\label{cumulant}
	The limit cumulant generating function exists and is given by
	\begin{equation*}
	\begin{aligned}
		\L\left(t\right) &=  \lim_{N\rightarrow\infty}\frac{1}{\beta\abs{\Lambda_N}}\log\mathbb{E}_{\nu_{N,\alpha}}\Big[\exp\Big(\beta\abs{\L_N}\langle t,\blambda_N\rangle\Big)\Big] \\ &		= \sum_{k\in\N}\frac{q^{\ssup{\alpha}}_k}{\beta}\left(\e^{\beta t_k}-1\right)<\infty\,,\qquad t\in\ell_\infty\left(\R\right)\,.
		\end{aligned}
	\end{equation*}
	Moreover, $\L$ is G{\^a}teaux differentiable, lower semicontinuous, and strictly convex.
\end{lemma}

\begin{proof}
	First, let us evaluate the logarithmic moment generating function. Recall, that our reference process is  a independent superposition of countably many independent marked Poisson point processes.  Denote the marginal law of $\blambda^{\ssup{k}}_N$ by $\nu^{(k)}_N$, then we have,
$$
\begin{aligned}
	\L\left(t\right) &=   \lim_{N\rightarrow\infty}\frac{1}{\beta\abs{\Lambda_N}}\log\mathbb{E}_{\nu_{N,\alpha}}\Big[\exp\Big(\beta\abs{\L_N}\langle t,\blambda_N\rangle\Big)\Big]= \lim_{N\rightarrow\infty}\frac{1}{\beta\abs{\Lambda_N}}\sum_{k\in\N}\log \mathbb{E}_{\nu_{N,\alpha}}\left[\exp\left(\beta\abs{\Lambda_N}t_k\blambda^{\ssup{k}}_N\right)\right]\\
	&=  \sum_{k\in\N}\frac{q^{\ssup{\alpha}}_k}{\beta}\left(\ex^{\beta t_k}-1\right).
	\end{aligned}
$$

	To see that $\L\left(t\right)$ is finite, note that $t\in\ell_\infty\left(\R\right)$ implies that $T:=\sup_{j\in\N}t_j$ is finite. Hence 
	\begin{equation*}
		\L\left(t\right)\leq\left(\ex^T-1\right)\bar{q}^{\ssup{\alpha}}<\infty.
	\end{equation*}
	
	To confirm G{\^a}teaux differentiability, let $t,s\in\ell_\infty\left(\R\right)$ and consider
	\begin{equation*}
		\frac{\d}{\d \epsilon}\L\left(t+\epsilon s\right) = \sum_{k\in\N}\frac{q^{\ssup{\alpha}}_k}{\beta}\, s_k \ex^{\beta(t_k+\epsilon s_k)}.
	\end{equation*}
	This sum is finite because $t$ and $s$ are bounded above and $q^{\ssup{\alpha}}\in\ell_1\left(\R\right)$ for $\alpha \le 0$. In particular, the derivative is defined at $\epsilon=0$, and hence $\L$ is G{\^a}teaux differentiable.
	
	Lower semicontinuity is an immediate consequence of Fatou's Lemma. For any sequence $t^{(n)}\rightarrow t$ in $\ell_\infty\left(\R\right)$,
	\begin{equation*}
		\liminf_{n\rightarrow\infty}\L\big(t^{(n)}\big) = \liminf_{n\rightarrow\infty}\sum_{k\in\N}\frac{q^{\ssup{\mu}}_k}{\beta}\left(\ex^{\beta t^{(n)}_k}-1\right) \geq \sum_{k\in\N}\frac{q^{\ssup{\mu}}_k}{\beta}\left(\ex^{\beta t_k}-1\right) = \L\left(t\right).
	\end{equation*}
	
	To show strict convexity, consider distinct $t,s\in\ell_\infty\left(\R\right)$ and $\lambda\in\left[0,1\right]$. Then
	\begin{align*}
		\L\left(\lambda s + \left(1-\lambda\right)t\right) &= \sum_{k\in\N}\frac{q^{\ssup{\mu}}_k}{\beta}\left(\ex^{\beta(\lambda s_k + \left(1-\lambda\right)t_k)}-1\right)\\
		&< \lambda\sum_{k\in\N}\frac{q^{\ssup{\alpha}}_k}{\beta} \ex^{\beta s_k} + \left(1-\lambda\right)\sum_{k\in\N}\frac{q^{\ssup{\alpha}}_k}{\beta}\ex^{\beta t_k} - \frac{\bar{q}^{\ssup{\alpha}}}{\beta}\\
		&= \lambda \L\left(s\right) + \left(1-\lambda\right)\L\left(t\right),
	\end{align*}
	where the inequality follows from the strict convexity of the exponential function.
\end{proof}

\begin{remark}
	If we do not have $\alpha\leq0$, then we do not have $\L\left(t\right)<\infty$ for all $t\in\ell_\infty\left(\R\right)$. To see this, let $t$ be a constant sequence $t_k=C>0$. Then $\L\left(t\right) = C\bar{q}^{\ssup{\alpha}} = \infty$ unless $\alpha\leq0$.\hfill $ \diamond $
\end{remark}

\begin{lemma}\label{legendre}
	For all $x\in\ell_1\left(\R\right)$, we have
	\begin{equation*}
		\L^*\left(x\right) := \sup_{t\in\ell_\infty\left(\R\right)}\left\{\left\langle t,x\right\rangle - \L\left(t\right)\right\} =  I_\alpha(x)\,.
	\end{equation*}
\end{lemma}

\begin{proof}
	Let $g_x\left(t\right)$ denote the functional we wish to maximise in the definition of $\L^*$, so
	\begin{equation*}
	g_x\left(t\right) = \sum^\infty_{k=1}\big[x_kt_k + \frac{1}{\beta}q^{\ssup{\alpha}}_k\big(1-\ex^{\beta t_k}\big)\big].
	\end{equation*}
	
	First let us consider $x\in\ell_1\left(\R\right)\setminus \ell_1\left(\R_+\right)$. Hence there exists an index $k^\prime$ such that $x_{k^\prime}<0$. Now let $t^{(T)}=-T\delta_{k^\prime}\in\ell_\infty\left(\R\right)$. Therefore
	\begin{equation*}
	\L^*\left(x\right)\geq g_x\big(t^{(T)}\big) = -Tx_{k^\prime} + \frac{1}{\beta}q^{\ssup{\alpha}}_{k^\prime}\big(1-\ex^{-\beta T}\big)
	\xrightarrow{T\rightarrow\infty}+\infty.
	\end{equation*}
	This means $\L^*\left(x\right) = +\infty =I_\alpha\left(x\right)$ for all $x\in\ell_1\left(\R\right)\setminus \ell_1\left(\R_+\right)$.
	
	To show the required inequality on $\ell_1\left(\R\right)$, let us now search for critical points of $g_x$. Taking the G{\^a}teaux derivative of $g_x$ gives us
	\begin{equation*}
	\d g_x\left(t;s\right) = \sum^\infty_{k=1}s_k\left(x_k - q^{\ssup{\alpha}}_k\ex^{\beta t_k}\right), \qquad\forall t,s\in\ell_\infty\left(\R\right).
	\end{equation*}
	Now $t$ is a critical point if and only if $\d g_x\left(t;s\right)=0$ $\forall s \in\ell_\infty\left(\R\right)$. This means that we want to investigate the sequence $\tilde{t}_k = \frac{1}{\beta}\log \frac{x_k}{q^{\ssup{\alpha}}_k}$. If $\tilde{t}\in\ell_\infty\left(\R\right)$, then this gives us the supremum, and a simple substitution tells us that $\L^*\left(x\right) = I_\alpha\left(x\right)$ for such $x$. Unfortunately, this is not necessarily the case.
	
	Nevertheless, these critical points will give us the supremum over all sequences in $\left(\R\cup\left\{-\infty\right\}\right)^\N$. Since $\ell_\infty\left(\R\right)\subset\left(\R\cup\left\{-\infty\right\}\right)^\N$, we have
	\begin{equation}
		\L^*\left(x\right) = \sup_{t\in\ell_\infty\left(\R\right)}g_x\left(t\right) \leq \sup_{t\in\left(\R\cup\left\{-\infty\right\}\right)^\N}g_x\left(t\right) = I_\alpha\left(x\right).
	\end{equation}
	
	To find the reverse inequality, let us consider
	\begin{equation*}
	t^{(K)}_k =
	\begin{cases}
	\mathds{1}\left\{k\leq K\right\} \beta^{-1}\log\frac{x_k}{q^{\ssup{\alpha}}_k} &,  x_k\ne0,\\
	-K\mathds{1}\left\{k\leq K\right\} &, x_k=0.
	\end{cases}
	\end{equation*}
	Since $t^{(K)}$ truncates, it is clearly in $\ell_\infty\left(\R\right)$ for all $K$. Now let us substitute it into $g_x$.
	\begin{align*}
	g_x\left(t^{(K)}\right) &= \sum_{k\leq K:x_k\ne0}\frac{1}{\beta}\Big(x_k\log\frac{x_k}{q^{\ssup{\alpha}}_k}-x_k+q^{\ssup{\alpha}}_k\Big) +  \sum_{k\leq K:x_k=0} \frac{1}{\beta}q^{\ssup{\alpha}}_k\left(1-\ex^{-\beta K}\right)\\
	&= \sum^K_{k=1}\frac{1}{\beta}\Big(x_k\Big(\log\frac{x_k}{q^{\ssup{\alpha}}_k}-1\Big)+q^{\ssup{\alpha}}_k\Big] - \ex^{-\beta K}\sum_{k\leq K\colon x_k=0} q^{\ssup{\alpha}}_k\\
	&\xrightarrow{K\rightarrow\infty} I_\alpha \left(x\right).
	\end{align*}
	
	This sequence $\left(t^{(K)}\right)_{K\in\N}$ shows that for $x\in\ell_1\left(\R_+\right)$,
	\begin{equation*}
	\L^*\left(x\right) = \sup_{t\in\ell_\infty\left(\R\right)}g_x\left(t\right)\geq I_\alpha\left(x\right),
	\end{equation*}
	as required.
\end{proof}

With Baldi's Theorem the proof of Theorem~\ref{THM-Ideal} is complete.

\qed

\end{proofsect}

\begin{proofsect}{Proof of Proposition~\ref{zeroIdeal}}
(a)  To find the zeroes of the ideal gas rate function, first let us find the critical points by setting the G{\^a}teaux derivative of the function to zero. That is, we find the set of points $\tilde{x}\in\ell_1\left(\R_+\right)$ such that
\begin{equation*}
    \d I_\alpha\left(\tilde{x};y\right)=0\quad\forall y\in\ell_1\left(\R\right).
\end{equation*}
This yields a single equation for each element of the sequence $\tilde{x}$. This set of equations has the unique solution $ \tilde{x}=\xi $ given by $\xi_k = q_k\e^{\beta \alpha k}$, for all $k\in\N$. Since the rate function $ I_\alpha $ is strictly convex where it is finite, this critical point is the unique global minimiser.

\noindent (b)  The existence of the thermodynamic limit and the explicit function follows from the large deviation rate function $ I_\alpha $.

\qed
\end{proofsect}

\begin{proofsect}{Proof of Proposition~\ref{P:ideal1}}
\noindent (a)   Clearly, $ p(\beta,0)=\frac{1}{\beta(4\pi\beta)^{\frac{d}{2}}}\sum_{k=1}^\infty\frac{1}{k^{1+d/2}}<\infty $ for all $ d\ge 1 $. Convexity follows from properties of the Bose functions, $g(1+\frac{d}{2},-\beta\mu) $, see \eqref{defBosefunctionapa} in Appendix~\ref{app-Bose}. 

\noindent (b) This follows from the Bose functions, $g(n,x)$, being differentiable for $x>0$, and
\begin{equation*}
		\frac{\d}{\d x}g(n,x) = -g(n-1,x), \qquad\forall x>0.
	\end{equation*}
	Then the first derivative follows from directly differentiating the representation \eqref{p}.

\noindent (c)  This follows by direct computation. The exponential term ensures that the derivative of the finite-volume pressure is increasing in $ \alpha $. As long as the box $ \L_N $ has finite volume one can give the average particle density any pre-assigned value by choosing a chemical potential.\\

\noindent (d) The limit in \eqref{critical} is obtained by direct calculation in conjunction with basic properties of the Bose function summarised in Appendix~\ref{app-Bose}. The convergence of the unique root is ensured as long as the expected particle density stays below the critical density which is finite only in dimensions $ d\ge 3 $.
\qed
\end{proofsect}

\begin{proofsect}{Proof of Proposition~\ref{IdealFreeEnergy}}
	Since $p\left(\beta,s\right)=+\infty$ for $s>0$, we only need to search $s\leq0$. On the interior of this region $p$ is differentiable, and we look for stationary points. If $\varrho\geq\varrho_\mathrm{c}$, then there are no stationary points for $s<0$ and $s\varrho-p\left(\beta,s\right)$ is increasing in $s$. Hence the supremum is achieved at $s = 0$. If $\varrho<\varrho_\mathrm{c}$, then there is a unique stationary point. This is also a local maximum and is given at $s = \alpha$ as required. This has the required limit as $\varrho\uparrow\varrho_\mathrm{c}$ implying the  continuity for $f$.    \qed
\end{proofsect}

\section{Lambert W function}\label{app-Lambert}

\begin{figure}
	\centering
	\begin{tikzpicture}[xscale=1.5]
	\draw[->] (-1.1,0) -- (3.1,0) node[below]{$x$};
	\draw[<-] (0,1.6) node[left]{$W\left(x\right)$} -- (0,-4.1);
	\draw[very thick] (-0.37,-1) to [out=90,in=225] (0,0) to [out=45,in=190] (3,1.5) node[right]{$W_0$};
	\draw[very thick, dashed] (-0.37,-1) to [out=270,in=90] (-0.05,-4) node[left]{$W_{-1}$};
	\draw[dashed] (-0.37,0) node[above]{$-1/e$} -- (-0.37,-1) -- (0,-1) node[right]{$-1$};
	\end{tikzpicture}
	\caption{\label{fig:lambertW} The two real branches of $W$: $W_0$ and $W_{-1}$.}
\end{figure}
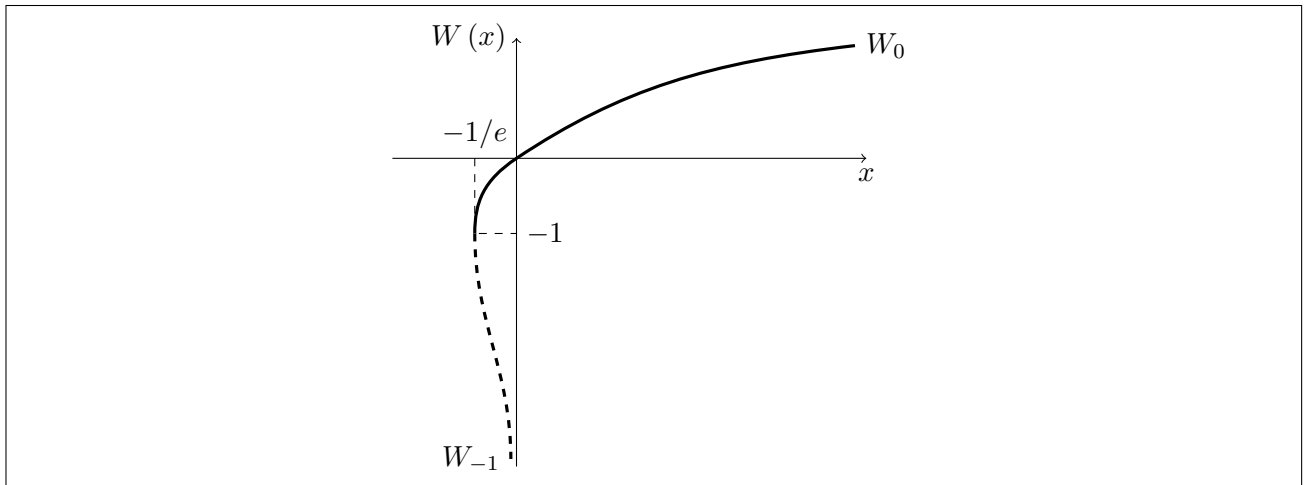

The Lambert W function (sometimes called elsewhere the \emph{Omega function}) is defined as the multi-valued inverse of the $\mathbb{C}\rightarrow\mathbb{C}$ function $w\mapsto w\ex^{w}$. We shall only be concerned with the two branches on $\R$. Figure \ref{fig:lambertW} shows these two real branches, denoted $W_0$ and $W_{-1}$. The $W_0$ branch is defined on $\left[-\ex^{-1},\infty\right)$, whereas the $W_{-1}$ branch is only defined on $\left[-\ex^{-1},0\right)$. Given a branch $W_l$ with $l\in\left\{0,-1\right\}$, we can find its (real) derivative $W'_l$ by differentiating the equation $W_l\left(x\right)\ex^{W_l\left(x\right)} = x$. This gives us
\begin{equation*}
W'_l\left(x\right) = \frac{1}{x}\frac{W_l\left(x\right)}{1+ W_l\left(x\right)}.
\end{equation*}
Taking further derivatives and applying induction shows that the branches are smooth on the interior of their respective domains, and gives expressions for each order of the derivative. We make use of some asymptotic expansions of $W_0$ and $W_1$:
\begin{equation*}
\begin{alignedat}{2}
W_0\left(x\right) &= x - x^2 + o\left(x^2\right)  & & \quad\text{as }x\to 0,\\
W_0\left(x\right) &= \log x - \log \left(\log x\right) + o\left(1\right) & &\quad\text{as }x\to+\infty,\\
W_{-1}\left(x\right) &= \log\left(-x\right) - \log\left(-\log\left(-x\right)\right) + o\left(1\right) & & \quad\text{as }x\uparrow 0.\\
\end{alignedat}
\end{equation*}
For more details, see \cite{CGHJK96}.

\end{appendices}

\end{document}